\documentclass[12pt]{article}
\usepackage{amsmath}
\usepackage{amssymb}
\usepackage{theorem}

\newcommand{\circledvartriangleleft}{\lessdot}

\newcommand{\eop}{\bigstar}

\newcommand{\card}[1]{{\vert #1 \vert} }
\newcommand{\norm}[1]{{\card{\card{#1}}}}

\newcommand{\proves}{\vdash}

\newcommand{\iso}{\thickapprox}

\newcommand{\Dom}{{\rm Dom}}

\newcommand{\llg}{{l\rm g}}

\newcommand{\Rang}{{\rm Rang}}

\newcommand{\cf}{{\rm cf}}

\newcommand{\initial}{\vartriangleleft}
\newcommand{\initialeq}{\vartrianglelefteq}

\newenvironment{proof}{\noindent{\bf Proof.}}{\par\bigskip}

\newenvironment{Proof}{\noindent{\bf Proof.}}{\par\bigskip} 

{\theorembodyfont{\rmfamily}\newtheorem{THEOREM}{Theorem}[section]}

{\theorembodyfont{\rmfamily}\newtheorem{Conclusion}[THEOREM]{Conclusion}}

{\theorembodyfont{\rmfamily}\newtheorem{LEMMA}[THEOREM]{Lemma}}

{\theorembodyfont{\rmfamily}\newtheorem{Main Theorem}[THEOREM]{Main Theorem}}
\newenvironment{main Theorem}{\begin{Main Theorem}} 
{\end{Main Theorem}}

{\theorembodyfont{\rmfamily}\newtheorem{Theorem}[THEOREM]{Theorem}}
\newenvironment{theorem}{\begin{Theorem}}{\end{Theorem}}

{\theorembodyfont{\rmfamily}\newtheorem{Definition}[THEOREM]{Definition}}

{\theorembodyfont{\rmfamily}\newtheorem{Conventions}[THEOREM]{Conventions}}

{\theorembodyfont{\rmfamily}\newtheorem{Main Definition}[THEOREM]{Main Definition}}
\newenvironment{main definition}{\begin{Main Definition}}
{\end{Main Definition}}

{\theorembodyfont{\rmfamily}\newtheorem{Lemma}[THEOREM]{Lemma}}

{\theorembodyfont{\rmfamily}\newtheorem{Notation}[THEOREM]{Notation}}

{\theorembodyfont{\rmfamily}\newtheorem{Convention}[THEOREM]{Convention}}

{\theorembodyfont{\rmfamily}\newtheorem{Note}[THEOREM]{Note}}

{\theorembodyfont{\rmfamily}\newtheorem{Observation}[THEOREM]{Observation}}

{\theorembodyfont{\rmfamily}\newtheorem{Remark}[THEOREM]{Remark}}

{\theorembodyfont{\rmfamily}\newtheorem{Main Fact}[THEOREM]{Main Fact}}
\newenvironment{main Fact}{\begin{Main Fact}}{\end{Main Fact}}

{\theorembodyfont{\rmfamily}\newtheorem{Fact}[THEOREM]{Fact}}

{\theorembodyfont{\rmfamily}\newtheorem{Subfact}[THEOREM]{Subfact}}

{\theorembodyfont{\rmfamily}\newtheorem{Claim}[THEOREM]{Claim}}
\newenvironment{claim}{\begin{Claim}}{\end{Claim}}

{\theorembodyfont{\rmfamily}\newtheorem{Main Claim}[THEOREM]{Main Claim}}
\newenvironment{main claim}{\begin{Main Claim}}{\end{Main Claim}}

{\theorembodyfont{\rmfamily}\newtheorem{Corrolary}[THEOREM]{Corrolary}}

{\theorembodyfont{\rmfamily}\newtheorem{Subclaim}[THEOREM]{Subclaim}}

{\theorembodyfont{\rmfamily}\newtheorem{Corollary}[THEOREM]{Corollary}}

{\theorembodyfont{\rmfamily}\newtheorem{Example}[THEOREM]{Example}}

{\theorembodyfont{\rmfamily}\newtheorem{Proposition}[THEOREM]{Proposition}}

{\theorembodyfont{\rmfamily}\newtheorem{Discussion}[THEOREM]{Discussion}}

{\theorembodyfont{\rmfamily}\newtheorem{Question}[THEOREM]{Question}}

\newenvironment{Proof of the Subfact}
{\noindent{\bf Proof of the Subfact.}}{\par\bigskip}

\newenvironment{Proof of the Theorem}
{\noindent{\bf Proof of the Theorem.}}{\par\bigskip}

\newenvironment{Proof of the Conclusion}
{\noindent{\bf Proof of the Conclusion.}}{\par\bigskip}

\newenvironment{Proof of the Observation}
{\noindent{\bf Proof of the Observation.}}{\par\bigskip}

\newenvironment{Proof of the Fact}
{\noindent{\bf Proof of the Fact.}}{\par\bigskip}

\newenvironment{Proof of the Lemma}
{\noindent{\bf Proof of the Lemma.}}{\par\bigskip}

\newenvironment{Proof of the Claim}
{\noindent{\bf Proof of the Claim.}}{\par\bigskip}

\newenvironment{Proof of the Subclaim}
{\noindent{\bf Proof of the Subclaim.}}{\par\medskip}

\newenvironment{Proof of the Main Claim}
{\noindent{\bf Proof of the Main Claim.}}{\par\bigskip}

\newcommand{\elementary}{\prec}

\newcommand{\vartrianglelefteq}{\trianglelefteq}

\newcommand{\tensor}{\otimes}

\newcommand{\into}{\rightarrow}

\newcommand{\rest}{\upharpoonright}  

\newcommand{\acc}{\mathop{\rm acc}}
\newcommand{\cl}{\mathop{\rm cl}}

\newcommand{\deq}{\buildrel{\rm def}\over =}

\newcommand{\DD}{{\cal D}}
\newcommand{\EE}{{\cal E}}

\newcommand{\HH}{{\cal H}}

\newcommand{\LL}{{\cal L}}

\title{On $\initial^\ast$-maximality}
\author{
Mirna D\v zamonja\\
School of Mathematics\\
University of East Anglia\\
Norwich, NR4 7TJ,UK\\
\scriptsize{M.Dzamonja@uea.ac.uk}\\
\scriptsize{http://www.mth.uea.ac.uk/people/md.html}\\
and\\
Saharon Shelah\\
Mathematics Department\\
Hebrew University of Jerusalem\\
91904 Givat Ram, Israel\\
\scriptsize{shelah@sunset.huji.ac.il}\\
\scriptsize{http://shelah.logic.at/}
}

\begin{document}

\baselineskip=16pt
\binoppenalty=10000
\relpenalty=10000
\raggedbottom

\maketitle

\begin{abstract} 
This paper investigates a connection between 
the semantic notion provided
by the ordering $\initial^\ast$ among theories in model theory and the
syntactic (N)SOP${}_n$ hierarchy of Shelah. It introduces two properties
which are natural extensions
of this hierarchy, called SOP${}_2$ and SOP${}_1$.
It is shown here that SOP${}_3$ implies SOP${}_2$ implies SOP${}_1$. In
\cite{Sh 500} it was
shown that SOP${}_3$ implies $\initial^\ast$-maximality and
we prove here that $\initial^\ast$-maximality in a model of GCH implies
a property called SOP${}''_2$. It has been subsequently shown by Shelah and
Usvyatsov that SOP${}''_2$ and SOP${}_2$ are equivalent, so obtaining an
implication between $\initial^\ast$-maximality and SOP${}_2$. It is not
known if SOP${}_2$ and SOP${}_3$ are equivalent.

Together with the known results
about the connection between the (N)SOP${}_n$ hierarchy and the existence
of universal models in the
absence of GCH, the paper provides a step toward the classification of
unstable theories without the strict order property.
\footnote{
This publication is numbered 692 in the list of publications
of Saharon Shelah. The authors
thank the United States-Israel Binational Science
Foundation for a grant supporting this research
and 
the NSF USA for their grant numbered NSF-DMS97-04477.
Mirna D\v zamonja thanks EPSRC for their support through grant number
GR/M71121, as well as the Royal Society for
their support through grant number SV/ISR/NVB.
We would also like to acknowledge the support of the Erd\"os Research Center
in Budapest, during the Workshop on Set-theoretic
Topology, July 1999, and support of the Hebrew University of Jerusalem and
the Academic Study Group during July 1999. Finally, warm thanks are
due to Alex Usvyatsov for his comments and improvements to the
manuscript. 


} 

{\em 
Changes from the published version:

In the published version of this paper
it is claimed that witnesses to being SOP${}_1$ can be chosen to be 
highly indiscernible, and this is justified by a certain notion of 1-fbti. The definition of this notion
(Definition 2.10) has a typo in a crucial place, and in addition Claim 2.11 for $t=2$ is incorrect and 
for $t=1$ the proof is incomplete. In this version we clarify these statements and proofs by introducting a new notion of indescernibility 3-fbti. The corrected statement is that witnesses to being SOP${}_1$ can be chosen to be 3-fbti.

That there are inconsistencies in the notions we used in the original paper was first observed by Lynn Scowl (September 2008),
Byunghan Kim (May 2009) and Enrique Casanovas and Martin Ziegler (July 2010). Whilst a Ph.D. student at UEA in 2008, Mark Wong also observed some incosistencies and made partial progress in rectifying them. There is a paper by Kim and Kim (to appear in APAL as of March 2011) which gives a different notion of 1-fbti and shows that witnesses can be chosen with that kind of indiscernibility. }

\end{abstract} 

\section{Introduction}\label{intro}

This paper investigates a connection between 
the ordering $\initial^\ast$ among theories in model theory and the
(N)SOP${}_n$ hierarchy of Shelah and as such provides a step toward the classification of
unstable theories without the strict order property. The thesis we pursue is that
the syntactic property SOP${}_2$ is closely related to the
semantic property of being maximal in the $\initial^\ast$-order.
We shall now give the relevant definitions and explain the motivation behind the paper
as well as noting our main results.
For the purpose of this introductory discussion we shall
limit ourselves to countable (complete first order) theories.

The following order among theories was introduced and investigated
by Keisler in \cite{Keisler}. 

\begin{Definition}\label{Keisl} (1) For any cardinal $\lambda$,
the Keisler order $\circledvartriangleleft_\lambda$
among theories is defined as follows: $T_0\circledvartriangleleft_\lambda
T_1$ if whenever $M_l (l<2)$ is a model of $T_0, T_1$ respectively
and $\DD$ is a regular ultrafilter over $\lambda$, \underline{then}
the $\lambda^+$-compactness of $M_1^\lambda/\DD$ implies the
$\lambda^+$-compactness of $M_0^\lambda/\DD$.

{\noindent (2)} We say $T_0\circledvartriangleleft T_1$ if for all
$\lambda$ we have $T_0\circledvartriangleleft_\lambda T_1$.
\end{Definition}

The relevance of this order to the project of 
classifying unstable theories without strict order property
lies in the two 
following theorems of Shelah (note that the second one implies the first).

\begin{Theorem} (Shelah \cite{Sh c}, VI4.3) Any (countable) theory 
with the strict order property is $\circledvartriangleleft$-maximal.
\end{Theorem}

As stated in \cite{Sh c}, pg xiv, Ch VI of \cite{Sh c} gives a rather complete
picture of Keisler's order and to complete it we should know more about
unstable theories without the strict order property. Paper \cite{Sh 500}
started a classification of such theories by introducing the hierarchy SOP${}_n$ for $n\ge 3$
and in particular it is stated there that being maximal in the Keisler order
is not a characterisation of theories with the strict order property,

\begin{Theorem}\label{dwa} (Shelah \cite{Sh 500}, see also \cite{ShUs xx}) Any theory 
with SOP${}_3$ is $\circledvartriangleleft$-maximal. 
\end{Theorem}

Details of the proof are given in \cite{ShUs xx}. 
Precise definitions of properties SOP${}_n$
for $n\ge 3$ will be repeated below in \S\ref{druga} but for the moment we note that it was
proved in \cite{Sh 500} that for $n\ge 3$
\[
\mbox{ strict order property}\implies SOP_{n+1}\implies SOP_n\implies\mbox{ not simple}
\]
and that all the implications are irreversible.
One may now wonder if having SOP${}_3$
is a characterisation of theories that are maximal in the Keisler order, giving us
a semantic equivalent to the syntactic notion of SOP${}_3$. 
This would be consistent with what is known about this order, see the Introduction
to Ch VI of \cite{Sh c}. This question remains open but instead one may attempt to give
a characterisation of SOP${}_3$ or SOP${}_n$ in terms of some other similarly defined order. 
This is suggested by \cite{Sh 500} which in fact gives a theorem stronger than \ref{dwa},
namely 

\begin{Theorem}\label{dwaprime} (Shelah \cite{Sh 500}, see also \cite{ShUs xx}) Any theory 
with SOP${}_3$ is $\vartriangleleft^\ast$-maximal. 
\end{Theorem}

The definition of this order will be recalled in 
\S\ref{first} where we shall also prove that being $\vartriangleleft^\ast$-maximal
implies being maximal in the Keisler order. Given this fact one may now ask if
being $\vartriangleleft^\ast$-maximal characterises theories with SOP${}_3$. To test
this claim it is natural to investigate a prototypical example of an NSOP${}_3$
theory that is still not simple, which is $T^\ast_{\rm feq}$. In \S\ref{first}
we shall recall the definition of this theory and show that in fact it is not
$\vartriangleleft^\ast$-maximal, as it is consistently 
strictly below the theory of a dense linear
order with no first or last element (all we need for the consistency is a partial GCH
assumption).

This naturally leads to the question of the possibility of refining the distinction
between simplicity and SOP${}_3$. Definition of the SOP${}_n$ hierarchy from \cite{Sh 500}
does not immediately give way to such a refinement as SOP${}_n$ is roughly speaking,
defined in terms of omitting loops of size $n$. However in \S\ref{druga} we introduce two properties
SOP${}_2$ and SOP${}_1$ that in fact satisfy 
\[
{\rm SOP}_3\implies {\rm SOP}_2 \implies {\rm SOP}_1\implies \mbox{ not simple}.
\]
We then ask if these properties in any way characterise the maximality in $\initial^\ast$.
To this end in \S\ref{ekvivalencija} we prove that any theory that is $\initial^\ast$-maximal
in a model of a sufficient amount of GCH must satisfy a syntactic property SOP${}''_2$.
Together with a subsequent result of Shelah and Usvyatsov in \cite{ShUs xx} that proved
that SOP${}''_2$ is equivalent to SOP${}_2$ we hence obtain that
$\initial^\ast$-maximality in any model of a sufficiently rich fragment of GCH implies
SOP${}_2$. (See \S3 for the definition of SOP${}''_2$ and the exact reference from \cite{ShUs xx}).
To summarise, our main result, appearing as Corollary \ref{zakljucak}(1) below is

\begin{Theorem}\label{zakljucakint} Suppose that $T$ is a theory that is
$\initial^\ast$-maximal in some universe of set theory in which $2^\lambda=\lambda^+$ holds
for all large enough regular $\lambda$. \underline{Then} $T$ has SOP${}_2$.
\end{Theorem}

Several questions remain open. The main one of course is if SOP${}_2$
is actually equivalent to $\initial^\ast$-maximality. Recall from the discussion above that
we know that SOP${}_3$ implies $\initial^\ast$-maximality. It is not known if SOP${}_3$
and SOP${}_2$ are actually equivalent.
We also note that 
Shelah and Usvyatsov have proved in \cite{ShUs xx} a local version of the implication
${\rm SOP}_2\implies\initial^\ast{\rm-maximality}$, see \S\ref{ekvivalencija} for a more
detailed discussion.

A burning question also is that we in fact do not know almost anything about the reverse
of other implications in the (consistent) diagram
\[
{\rm SOP}_3\implies\initial^\ast{\rm-maximality}\implies 
{\rm SOP}_2 \implies {\rm SOP}_1\implies \mbox{ not simple},
\]
apart that not all of them may be equivalences, as $T^\ast_{\rm feq}$ is not
simple but is NSOP${}_3$. In fact \cite{ShUs xx} proves that $T^\ast_{\rm feq}$ is not even
SOP${}_1$.

Before laying down the organisation of the paper let us also mention the connection
of the SOP${}_n$ hierarchy with another semantic property, which is the possibility of
having a universal model at $\lambda$ in some universe of set theory
where a sufficient amount of
GCH fails (under GCH every countable first order theory has a universal model in
every uncountable cardinal). The connection between this property and unstable
theories without the strict order property has been investigated in a series of papers,
notably in \cite{KjSh 409} where it is proved that if GCH fails sufficiently then there
are no universal dense linear orders. It was proved in \cite{Sh 500} that SOP${}_4$ is already
sufficient for such a negative universality result. The question of universality is 
interesting also for classes that are not elementary classes of models of a first order
theory, for example for classes without amalgamation the most interesting case is
the strong limit singular $\mu$ of cofinality $\aleph_0$. In \cite{GrSh 174} it is proved
that for such $\mu$ and $\lambda<\mu$ a strongly compact cardinal the class of models of
any $L_{\lambda,\mu}$-theory of cardinality $<\mu$ admits a universal model of
cardinality $\mu$.
A rather detailed description of what is known about the connection of unstable theories
without the strict order property and the universality problem may be found in the
introduction to \cite{DjSh 710}.

The paper is organised as follows.
In the first section we investigate the theory $T^\ast_{\rm feq}$.
This is simply the model completion
of the theory of infinitely many parametrised equivalence relations. We
show that
under a partial GCH assumption,
this theory is not maximal with respect to
$\initial^\ast_\lambda$,
as it is strictly below the theory of a
dense linear order. In the second section of the paper we
extend Shelah's NSOP${}_n$ hierarchy by introducing two further properties
SOP${}_1$ and SOP${}_2$, and we show that their names
are justified by their position in the hierarchy. 
Namely SOP${}_3\implies$ SOP${}_2\implies$ SOP${}_1$. Furthermore,
SOP${}_1$ theories are not simple. The last
section of the paper contains the main result showing that $\initial^\ast$-maximality
in a model of a sufficiently rich fragment of GCH implies SOP${}''_2$, and
hence SOP${}_2$ by Shelah-Usvyatsov.


The following conventions will be
used in the paper.

\begin{Convention} Unless specified otherwise, a ``theory" stands for a
first order complete theory. An unattributed $T$ stands for a theory.
We use $\tau(T)$ to denote the vocabulary of a theory $T$,
and $\LL(T)$ to denote the set of formulae of $T$.

By ${\frak C}={\frak C}_T$ we denote a $\bar{\kappa}$-saturated model of $T$,
for a large enough regular cardinal $\bar{\kappa}$
and we assume that any models of $T$ that we mention are elementary submodels
of ${\frak C}$.

$\lambda,\mu,\kappa$ stand for infinite cardinals.
\end{Convention}

\section{On the order $\initial_\lambda^\ast$}\label{first}

\begin{Definition}(1) \label{interpretations}
For (first order complete) theories $T_0$ and $T$ we say
that 
\[
\bar{\varphi}=\langle \varphi_R (\bar{x}_R):\, R\mbox{ a predicate of }
\tau(T_0)\mbox{ or a function symbol of }\tau(T_0)\mbox{ or }=\rangle,
\]
(where we have $\bar{x}_R=(x_0,\ldots x_{n(R)-1})$),
{\em interprets} $T_0$ in $T$,
or that $
\bar{\varphi}$ {\em is an interpretation of}
$T_0$ in $T$, or that
\[
T\proves``
\bar{\varphi}\mbox{ is a model of }T_0",
\]
\underline{if} each $\varphi_R(\bar{x}_R)\in \LL(T)$, and for any $M\models T$,
the model
$M^{[\bar{\varphi}]}$ described below
is a model of $T_0$. Here,
$N=M^{[\bar{\varphi}]}$ is a $\tau(T_0)$ model, whose set of elements is $\{a:\,
M\models\varphi_=(a,a)\}$ (so $M^{[\bar{\varphi}]}\subseteq M$) and
$R^N=\{\bar{a}:\,M\models\varphi_R[\bar{a}]\}$ for a predicate $R$ of $T_0$.

For any function symbol $f$ of $\tau(T_0)$ we have that
$N\models``f(\bar{a})=b"$ iff $M\models\varphi_f(\bar{a},b)$,
while
\[
M\models``\varphi_f(\bar{a},b)=\varphi_f(\bar{a},c) \implies
b=c"
\]
for all $\bar{a},b,c$.

{\noindent (2)} We say that the interpretation $\bar{\varphi}$ is
{\em trivial} if $\varphi_R(\bar{x}_R)=R(\bar{x}_R)$ for all
$R\in\tau(T_0)$, so $M^{[\bar{\varphi}]}=M\rest\tau(T_0)$, for any
model $M$ of $T$.
\end{Definition}

(The last clause in Definition \ref{interpretations}(1)
shows that we can in fact
restrict ourselves to vocabularies without function or constant
symbols.) 

We use the notion of interpretations to define a certain relation
among theories. This relation was introduced
by S. Shelah in \cite{Sh 500}, section \S2 and one can see \cite{ShUs xx} for a more
detailed exposition. The reason we are interested in this ordering is Shelah's 
Theorem \ref{dwa} quoted in the Introduction and we shall now start developing methods
for the proof of our main result \ref{zakljucak}.

\begin{Definition}\label{order} For (complete first order) theories $T_0,T_1$
we define:
\begin{description}

\item{(1)} A triple $(T, \bar{\varphi}_0, \bar{\varphi}_1)$
is called a $(T_0,T_1)$-{\em superior} iff $T$ is a theory
and $ \bar{\varphi}_l$ is an interpretation of $T_l$ in $T$,
for $l<2$.

\item{(2)} For a cardinal $\kappa$, a
$(T_0,T_1)$-superior $(T, \bar{\varphi}_0, \bar{\varphi}_1)$ is
called $\kappa$-{\em relevant} iff
$\card{T}<\kappa$.

\item{(3)} For regular cardinals $\lambda,\mu$ we say
 $T_0\initial^\ast_{\lambda,\mu} T_1$ if there
is a $\min(\mu,\lambda)$-relevant $(T_0,T_1)$-superior
triple $(T, \bar{\varphi}_0, \bar{\varphi}_1)$ such that
in every model $M$ of $T$ in which $M^{[\bar{\varphi}_1]}$ is
$\mu$-saturated, the model $M^{[\bar{\varphi}_0]}$ is
$\lambda$-saturated. If this happens, we call the triple a
{\rm witness} for $T_0\initial^\ast_{\lambda,\mu} T_1$.

\item{(4)} We say that $T_0\initial^\ast_{\lambda,\mu} T_1$ {\em over} $\theta$ if
$\theta\le\lambda,
\theta\le\mu$ and $T_0\initial^\ast_{\lambda,\mu} T_1$ as witnessed by a
$(T,\bar{\varphi_0},\bar{\varphi_1})$ with $\card{T}<\theta$.

\item{(5)} If $\lambda=\mu$, we write $\initial^\ast_\lambda$ in place of
$\initial^\ast_{\lambda,\mu}$.

\item{(6)} We say that $T_1\initial^\ast T_2$ iff $T_1\initial^\ast_\lambda T_2$
holds for all large enough regular $\lambda$.

\item{(6)} $T^\ast$ is $\initial^\ast_\lambda$-maximal iff $T\initial^\ast_\lambda T^\ast$
holds for all $T$. The notion of $\initial^\ast$-maximality is defined analogously.

\item{(7)} We say $T_0 \initial^\ast_{\lambda,\neq} T_1$ iff
$T_0\initial^\ast_\lambda T_1$ but $\neg (T_1\initial^\ast_\lambda T_0)$.

\end{description}
\end{Definition}

Although in this paper we do not consider this
in its own right, it is natural to define the local versions of the $\initial^\ast$-relation.
This is used by Shelah and Usvyatsov in \cite{ShUs xx} to obtain their local converse
to the implication $\initial^\ast$-maximality $\implies$ SOP${}_2$, see \S\ref{ekvivalencija}
for more discussion on this.

\begin{Definition}\label{localised}
Relations $\initial^{{\ast,\rm l}}_{\lambda,\mu}$ and
$\initial^{{\ast,\rm l}}_{\lambda}$ are the local
versions of $\initial^\ast_{\lambda,\mu}$ and
$\initial^\ast_{\lambda}$ respectively, where by a local version
we mean that in the definition of the relations, only types
of the form
\[
p\subseteq\{\pm \vartheta(x,\bar{a}):\,\bar{a}\in{}^{\llg(\bar{y})}M\}
\]
for some fixed $\vartheta(x,\bar{y})$ are considered.
\end{Definition}

\begin{Observation}\label{trivint} (0) If $T_0\initial^\ast_{\lambda,\mu} T_1$ and $l<2$,
\underline{then} there is a witness $(T, \bar{\varphi}^0,
\bar{\varphi}^1)$ such that $\bar{\varphi}^l$ is trivial, hence $T_l\subseteq
T$.

{\noindent (1)} $\initial^\ast_\lambda$ is a partial order among
theories (note that $T\initial^\ast_\lambda T$ for every complete $T$ of
size $<\lambda$,
and that the strict inequality is written as $T_1\initial^\ast_{\lambda,\neq}
T_2$).

{\noindent (2)} If $T_0\initial^\ast_{\lambda,\mu} T_1$ over $\theta$ and
$T_1\initial^\ast_{\mu,\kappa}T_2$ over $\theta$, \underline{then}
$T_0\initial^\ast_{\lambda,\kappa} T_2$ over $\theta$.
\end{Observation}

[Why? (0) Trivial.

(1) Suppose that $T_l\initial^\ast_\lambda T_{l+1}$ for $l<2$ over $\theta$, as exemplified
by $(T^\ast, \bar{\varphi}_0, \bar{\varphi}_1)$ and
$(T^{\ast\ast}, \bar{\psi}_1, \bar{\psi}_2)$ respectively. Without loss of
generality, 
$\bar{\varphi}_1$ is trivial (apply part (0)),
so as $T^\ast$ is complete we have $T_1\subseteq T^\ast$. Similarly,
without loss of generality, $\bar{\psi}_1$ is trivial
and so, as $T^{\ast\ast}$ is complete, we have $T_1\subseteq T^{\ast\ast}$.
As $T_1$ is complete, without loss of generality, $T^\ast$ and
$T^{\ast\ast}$ agree on the common part of their vocabularies, and
hence by Robinson Consistency Criterion, $T\deq T^\ast\cup T^{\ast\ast}$
is consistent. Also $\card{T^\ast}+\card{T^{\ast\ast}}<\theta$, hence
$\card{T}<\theta$. Clearly
$T$ interprets $T_0, T_1, T_2$ by $\bar{\varphi_0}$,
$\bar{\varphi_1}=\bar{\psi}_1$ and $\bar{\psi}_2$ respectively and
$T$ is complete. We now show
that the triple $(T, \bar{\varphi}_0, \bar{\psi}_2)$ is a $(T_0, T_2)$-superior
which witnesses $T_0\initial^\ast_{\lambda} T_2$ over $\theta$. So
suppose that $M$ is a model of $T$ in which $M^{[\bar{\psi}_2]}$
is $\lambda$-saturated. As $(T^{\ast\ast}, \bar{\psi}_1, \bar{\psi}_2)$
witnesses $T_1\initial^\ast_{\lambda} T_2$, 
we can conclude that $M^{[\bar{\varphi}_1]}=M^{[\bar{\psi}_1]}$
is $\lambda$-saturated. We can argue similarly that $M^{[\bar{\varphi}_0]}$
is $\lambda$-saturated.

(2) is proved similarly to (1).]

\medskip

In this section we consider an example of a theory
which is a prototypical
example of an NSOP${}_3$ theory that is not simple (see \cite{Sh 457}).
It is the model
completion of the theory of infinitely many (independent)
parametrised equivalence
relations, formally defined below. We shall prove that for
$\lambda$ such that $\lambda=\lambda^{<\lambda}$ and $2^\lambda=
\lambda^+$, this theory is 
strictly $\initial^\ast_{\lambda^+}$-below the theory
of a dense linear order with no first or last element.

\begin{Definition} (1) $T_{\rm feq}$ is the following theory
in $\{P, Q, E, R, F\}$
\begin{description}
\item{(a)}
Predicates $P$ and $Q$ are unary and disjoint, and
$(\forall x)\,[P(x)\vee Q(x)]$,
\item{(b)} $E$ is an equivalence relation 
on $Q$,
\item{(c)} $R$ is a binary relation on $Q\times P$ such that
\[
[x\,R\,z\,\,\&\,\,y\,R\,z\,\,\&\,\,x\,E\,y]\implies x=y.
\]
{\scriptsize{(so $R$ picks for each $z\in P$ (at most one)
representative of any $E$-equivalence class).}}
\item{(d)} $F$ is a (partial) binary function from $Q\times P$ to $Q$, which
satisfies
\[F(x,z)\in Q\,\,\&\,\, \,F(x,z)\,R\,z \,\,\&\,\,
x\, E\, F(x,z).
\]
{\scriptsize{(so for $x\in Q$ and $z\in P$, the function $F$ picks the
representative of the $E$-equivalence class of $x$ which is in the
relation $R$ with $z$).}}
\end{description}

{\noindent (2)} $T^+_{\rm feq}$ is $T_{\rm feq}$ with the requirement that
$F$ is total.

{\noindent (3)} For $n<\omega$, we let $T^n_{\rm feq}$ be
$T_{\rm feq}^+$ enriched by the sentence saying that over any
$n$
elements, any (not necessarily complete) 
quantifier free type consisting of basic (atomic and negations of the atomic)
formulae
with no direct contradictions, is realised.

\end{Definition}

\begin{Note} One may easily check that every model of $T_{\rm feq}$
can be extended to a model of
$T^+_{\rm feq}$ and that $T^+_{\rm feq}$ has the amalgamation property and the joint
embedding property. This theory also has a model completion, which can be
constructed directly, and which we denote by $T^\ast_{\rm feq}$. It follows
that
$T^\ast_{\rm feq}$ is a complete theory with infinite models, in which
$F$ is a full function.
\end{Note}

\begin{Remark} Notice that $T_{\rm feq}$ has been defined somewhat
differently than in [Sh 457, \S1], but the difference is non-essential, as the
following Claim \ref{modcom} shows that the two theories have the same model completion.
This claim also shows the origin of the name ``infinitely many 
independent equivalence relations" for $T^\ast_{\rm feq}$.
\end{Remark}

\begin{Claim}\label{modcom} Let $T$ be the theory defined (in [Sh457, 1]) by
\begin{description}
\item{(a)} $T$ has unary predicates $P$ and $Q$ and a three place relation $E$
writen as $y\,E_x z$,
\item{(b)} the universe of any model of $T$ is a disjoint union of $P$ and $Q$,
\item{(c)} $y \,E_x z\implies P(x)\,\,\&\,\,Q(y), Q(z)$,
\item{(d)} for any fixed $x\in P$ the relation $E_x$ is an equivalence relation on $Q$.
\end{description}
Then $T^\ast_{\rm feq}$ is the model completion of $T$.
\end{Claim}

\begin{Proof of the Claim} Let $M$ be a model of $T_{\rm feq}$, we shall extend $M$ to
a model of $T$ as follows. 
Each $E$-equivalence class $e=a/E$ gives
rise to an equivalence relation $E_e$ on $P$ given by:
\[
z_1E_e z_2\mbox{ iff }z_1,z_2\in P\mbox{ and } F(a,z_1)=F(a,z_2).
\]
This definition does not depend on $a$,
just on $a/E$. Let $P^N$ and $Q^N$ be $Q^M$ and $P^M$ respectively. Define $y\,E^N_x z$
iff $y\,E_e z$ where $e=x/E^M$. Clearly $N$ is a model of $T$.

Now suppose that we have a model $M$ of $T$ and we shall extend it to a model $N$ of
$T_{\rm feq}$. Let $P^N$ and $Q^N$ be $Q^M$ and $P^M$ respectively. Define $x\,E^N x'$
iff for every $y,z$ we have $y\,E_x z$ iff $y\,E_{x'} z$. Choose a representative
of each $E$-equivalence class and for any $z\in Q^N$ and such a representative $x$
let $F(x,z)=x$. Then for $x'\in Q^N$ which has not been chosen as a representative
of any equivalence class, find $x$ which has been chosen as its representative
and define $F(x',z)=F(x,z)$ for all $z\in P^N$.

This shows that $T_{\rm feq}$ and $T$ are cotheories (\cite{ChKe}, 3.5.6(2)). Being
the model completion of $T_{\rm feq}$, $T_{\rm feq}^\ast$ is its cotheory, and hence
a cotheory of $T$. Hence $T_{\rm feq}^\ast$ is a model companion of $T$. In order to
prove that it is the model completion of $T$ it suffices to show that $T$ has
the amalgamation property (\cite{ChKe}, 3.5.18) which is easily seen directly.
$\eop_{\ref{modcom}}$
\end{Proof of the Claim}

\begin{Observation}\label{eliminacija} $T^\ast_{\rm feq}$ has elimination of
quantifiers and for any $n$, any model of $T^\ast_{\rm feq}$ is a model of $T^n_{\rm feq}$.
\end{Observation}

\begin{Notation} $T_{\rm ord}$ stands for the theory of a
dense linear order with no first or last element.
\end{Notation}

The following convention will make the notation used in this
section simpler.

\begin{Convention}\label{simplification}
Whenever considering
$(T_{\rm ord},T^\ast_{\rm feq})$-superiors
$(T, \bar{\varphi}, \bar{\psi})$ we shall abuse the notation
and assume $\bar{\varphi}=(I, <_0)$ and
$\bar{\psi}=(P,Q,E,R, F)$. In such a case we may also write $P^M$
in place of $P^{M^{[\bar{\psi}]}}$ etc., and we may simply say
that $T$ is a $(T_{\rm ord},T^\ast_{\rm feq})$-superior.

\end{Convention}

We intend to prove that for $\lambda$ satisfying $\lambda^{<\lambda}$ and
$2^\lambda=\lambda^+$ the theory $T^\ast_{\rm feq}$ is strictly
$\initial^\ast_{\lambda^+}$-below 
$T_{\rm ord}$ (Theorem \ref{notmax} below). This will be done by a diagonalisation
argument where for a given $\lambda$-relevant
$(T_{\rm ord}, T_{\rm feq}^\ast)$-superior $T$ we inductively construct a model
of $T$ that is saturated for $T_{\rm feq}^\ast$ but not for $T_{\rm ord}$. Main Claim
\ref{main} provides one step in the required induction. In Stage A of its proof
we use the elimination of quantifiers in $T_{\rm feq}^\ast$ to
reduce the situation to $T_{\rm feq}$-types of four prescribed kinds, and then we
show that we may in fact work only with three of them. Stage B contains the main
point of the proof, which is the construction of a certain tree of models and
embeddings. Once this is done in Stage C we use the analysis from Stage A to show that
the $T^\ast_{\rm feq}$-type defined by the union of the embeddings is consistent.
In Stage D we take $N\elementary \frak{C}$ of size
$\lambda$ that realises this type and show that
$N$ must omit most of the Dedekind cuts induced by the tree of embeddings, and that
most of these cuts are not definable over $N$. After an application of an appropriate
automorphism of $\mathfrak C$ this finishes the proof of the Main Claim. The
proof of the theorem then follows by induction. The cardinal arithmetic assumptions
are used in Stage D and in the inductive proof of the theorem.

\begin{Definition} For a $\lambda$-relevant $(T_{\rm
ord},T^\ast_{\rm feq})$- superior $T$, {\em the statement}
\[
\ast[M, \bar{a},\bar{b}]=\ast[M, \bar{a},\bar{b}, T, \lambda]
\]
means:
\begin{description}
\item{(i)}
$M$ is a model of $T$ of size $\lambda$,
\item{(ii)} $\bar{a}=\langle a_i:\,i<\lambda
\rangle$, $\bar{b}=\langle b_i:\,i<\lambda
\rangle$, are sequences of elements of $I^{M^{[\bar{\varphi}]}}$ such that
\[
i<j<\lambda\implies a_i<_0 a_j<_0 b_j<_0 b_i,
\]
\item{(iii)}
there is no $x\in M^{[\bar{\varphi}]}$ such that for all $i$ we have $a_i<_0
 x<_0 b_i$,
\item{(iv)}
the Dedekind cut $\{x:\,\bigvee_{i<\lambda}x<_0 a_i\}$ is not
definable
by any formula of $\LL(M)$ with parameters in $M$.
\end{description}
\end{Definition}

\begin{Main Claim}\label{main} Assume $\lambda^{<\lambda}=\lambda$ and
$(T,\bar{\varphi}, \bar{\psi})$ is a $\lambda$-relevant
$(T_{\rm ord}, T_{\rm feq}^\ast)$-superior. Further assume
that
$\ast[M, \bar{a},\bar{b}]$ holds, and $p=p(z)$ is a
(consistent) $T_{\rm feq}^\ast$-type  over $M^{[\bar{\psi}]}$.
\underline{Then}
there is $N\models T$
with $M\elementary N$, such that $p(z)$ is realised in
$N^{[\bar{\psi}]}$ and $\ast[N, \bar{a},\bar{b}]$ holds.
\end{Main Claim}

\begin{Proof of the Main Claim} 

{\bf Stage} {\bf A}. Without loss of generality,
$p$ is complete in the
$T_{\rm feq}^\ast$-language over $M^{[\bar{\psi}]}$.
(By Convention \ref{simplification}, we can 
consider $p$ to be a type over $M$ (rather than $M^{[\bar{\psi}]}$).
We shall use this Convention throughout the proof). If $p$ is
realised in $M$, our conclusion follows by taking $N=M$, so let us assume that
this is not the case.
Using the elimination of quantifiers for $T^\ast_{\rm
feq}$,
we can without loss of generality assume that $p(z)$ consists
of quantifier free formulae with parameters in $M$. This means
that one of the following four cases must happen:

\underline{Case 1}. (This will be the main case)
$p(z)$ implies that $z\in P$
and it determines which elements of $Q^M$ are $R$-connected
to $z$. Hence
for some function $f:\,Q^{M}\into Q^{M}$
which respects $E$, i.e.
\[
a\,E\,b\implies f(a)=f(b),
\]
and
\[
f(a)\in a/E^{M};
\]
we have
\[
p(z)=\{P(z)\}\cup\{\,b \,R\,z:\,b\in \Rang(f)\}
\]
and no $a\in P^M$ satisfies $p$.

\underline{Case 1A}. Like Case 1, but $f$ is a partial function and
\begin{equation*}
\begin{split}
p(z)=\{P(z)\}&\cup\{\,f(b) \,R\,z:\,b\in \Dom(f)\}\\&\cup\{\,\neg (b
\,R\,z):\,(b/E^M \cap \Rang(f))=\emptyset\}.\\
\end{split}
\end{equation*}
(This Case will be reduced to Cases 1-3 in Subclaim \ref{addition}).

\underline{Case 2}. $p(z)$ determines that $z\in Q$ and that it is
$E$-equivalent to some $a^\ast\in Q^{M}$, but not
equal to any ``old" element. Note that
in this case if $b^\ast$ realises $p(z)$, we cannot have
$b^\ast R c$ for any $c\in P^{M}$, as this would
imply $F(a^\ast,c)=b^\ast\in M^{[\bar{\psi}]}$ (and we have assumed
that $p$ is not realised in $M^{[\bar{\psi}]}$). Hence, the complete
$M$-information is given by
\[
p(z)=\{Q(z)\}\cup \{a^\ast\,E\,z\}\cup\{a\neq z:\,a\in a^\ast/E^M\}.
\]

\underline{Case 3}.
$p(z)$ determines that $z\in Q$, but it has a different $E$-equivalence class
than any of
the elements of $Q^M$. As $p$ is complete, it must determine
for which $c\in P^M$ we have $z\,R\,c$, and for which $c,d\in
P^M$ we have $F(z,c)=F(z,d)$. Hence,
for some $f:\,P^M\into\{\mbox{yes, no}\}$
and some $\EE$, an equivalence relation on $P^M$ such that
$c\EE d\implies f(c)=f(d)$, we have
\begin{equation*}
\begin{split}
p= \{Q(z)\}\cup\{
\neg (a\,\,E\,\,z):\,a\in Q^M\}
\cup\{(z\,\,R\,\,b)^{f(b)}:\,
b\in P^M\}\\
\cup\{\left(F(z,c)=F(z,d)\right)^{{\rm if }c\EE d}:\,
c,d\in P^M\}.\\
\end{split}
\end{equation*}

In the above description, we have used

\begin{Notation} For a formula $\vartheta$ we let $\vartheta^{\rm
yes}\equiv\vartheta$ and $\vartheta^{\rm no}\equiv\neg\vartheta$.
\end{Notation}

\begin{Subclaim}\label{addition} It suffices to deal with Cases 1,2,3, ignoring
the Case 1A.
\end{Subclaim}

\begin{Proof of the Subclaim}
Suppose that we are in the Case 1A. Let \[
\{d_i/E^M:\, i<i^\ast\le\lambda\}
\]
list the $d/E^M$ for $d\in Q^M$ such that $d'\in d/E^M\implies \neg (d'\,R\,z)
\in p(z)$. 
We choose by induction on $i\le i^\ast$ a pair $(M_i, c_i)$ such that
\begin{description}
\item{(a)} $M_0=M$, $\norm{M_i}=\lambda$, 
\item{(b)} $\langle M_i:\,i\le i^\ast \rangle$ is an increasing continuous
elementary chain,
\item{(c)} $\ast[M_i,\bar{a},\bar{b}]$
\item{(d)} $c_i\in (d_i/E^{M_{i+1}})\setminus M_i$, for $i<i^\ast$.
\end{description}
For $i$ limit or $i=0$, the
choice is trivial. For the situation when $i$ is a successor, we use Case 2.

Let $\langle c_i/E^{M_{i^\ast}}:\,i\in [i^\ast, i^{\ast\ast})\rangle$ list
without repetitions the
$c/E^{M_{i^\ast}}$ which are disjoint to $M$. Note that $\card{i^{\ast\ast}}\le
\lambda$. Let
\[
p^+(z)\deq p(z)\cup\{c_i\,R\,z:\,i<i^{\ast\ast}\}.
\]
Then $p^+(z)$ is a complete type (for $M^{[\bar{\psi}]}_{i^\ast}$), and $\ast[M_{i^\ast},
\bar{a},\bar{b}]$ holds by (c). If $p^+(z)$ is realised in $M_{i^\ast}$, we
can let $N=M_{i^\ast}$ and we are done. Otherwise, $p^+(z)$ is not realised in
$M_{i^\ast}$ and is a type of the form required in Case 1, so we can proceed
to deal with it using the assumptions on Case 1.
$\eop_{\ref{addition}}$
\end{Proof of the Subclaim}

{\bf Stage} {\bf B}. Let us assume that $p$ is a type as in one of the Cases
1,2 or 3, which we can do by Subclaim \ref{addition}.
We shall define $\langle M_\alpha:\,\alpha<\lambda\rangle$, an
$\elementary$-increasing continuous sequence of elementary submodels of $M$, each of size $<\lambda$,
and with union $M$,
such that:
\begin{description}
\item{(a)} In Case 1, each $M_\alpha$ is closed under $f$,
\item{(b)} In Case 2, $a^\ast\in M_0$,
\item{(c)} For every $\alpha<\lambda$,
\[
(M_\alpha, \{(a_j, b_j):\,j<\lambda\}\cap M_\alpha)\elementary (M, \{(a_j,
b_j):\,j<\lambda\}),
\]
\end{description}

Hence,
for some club $C$ of $\lambda$ consisting of limit ordinals $\delta$,
we have that for all $\delta\in C$,
\[
a_j\in M_\delta\iff b_j\in M_\delta\iff j<\delta,
\]
\[
(\forall c\in I^{M_\delta})(\exists j<\delta)\,[c<_0 a_j\vee b_j<_0c].
\]
Let $C=\{\delta_i:\,i<\lambda\}$ be an
increasing continuous enumeration.

Now we come to the {\bf main point} of the proof.

By induction on $i=\llg(\eta)<\lambda$ we shall choose
$\bar{h}=\langle h_\eta:\,\eta\in
{}^{\lambda>}2\rangle$, a sequence such that
\begin{description}
\item{($\alpha$)} $h_\eta$ is an elementary embedding of
$M_{\delta_{\llg(\eta)}}$
into
${\mathfrak C}_T$, whose range will be denoted by $N_\eta$.
\item{($\beta$)} $\nu\initial \eta\implies h_\nu\subseteq h_\eta$.
\item{($\gamma$)} If $\eta_l\in{}^{\lambda>}2$
for $l=0,1$ and
$\eta_0\cap\eta_1=\eta$, \underline{then}:
\begin{description}
\item{(i)}
$N_{\eta_0}\cap
 N_{\eta_1}=N_\eta$.
\item{(ii)} 
In addition, if $a_l\in Q^{N_{\eta_l}}$ for $l=0,1$ and
$a_0\, E^{{\mathfrak C}_T} a_1$, \underline{then} for some
$a\in Q^{N_{\eta}}$ we have $a_l E^{{\mathfrak C}_T} a$
for $l=0,1$. (Equivalently, if $a_l\in Q^{N_{\eta_l}}$
and $\neg(\exists a\in N_\eta)(\wedge_{l<2} a_l\, E\, a)$,
then $\neg (a_0\, E\, a_1)$).
\end{description}

\item{($\delta$)} $\models``h_{\eta\frown\langle 0\rangle}(b_{
\delta_{\llg (\eta)}})
 <_0 h_{\eta\frown\langle 1\rangle}(a_{\delta_{\llg (\eta)}})"$
(see Convention \ref{simplification} on $<_0$).
 \end{description}
 
Note that the requirement of $h_\eta$ being elementary
and onto $N_\eta$ in particular implies that

\begin{description} 
\item{($\delta'$)} If for some $l<2$ and $\eta\in {}^{\lambda>}2$ we have
$a\in
N_{\eta\frown\langle l\rangle} \setminus N_\eta$ and $b\in N_\eta$, \underline{then}
$a E^{{\mathfrak C}_T} b$ iff $a=h_{\eta\frown\langle l\rangle}(a')$
for some $a'$ such that $a' E^{{\mathfrak C}_T} h_\eta^{-1}(b)$.
\end{description}

We now describe the inductive choice of $h_\eta$ for
$\eta\in {}^{\lambda>} 2$, the induction being on $i=\llg(\eta)$.
Let $h_{\langle\rangle}={\rm id}_{M_0}$.
If $i$ is a limit ordinal, we just let
$h_\eta\deq\bigcup_{j<\llg(\eta)}h_{\eta\rest j}$. Hence, the
point is to handle the successor case.

Fixing
$i<\lambda$, let $\langle \eta_{i,\alpha}:\,\alpha
<\alpha^\ast\le\lambda\rangle$ list $ {}^{i+1} 2$, in such a manner that
$\eta_{i,2\alpha}\rest i=
\eta_{i,2\alpha+1}\rest i$ and $\eta_{ i,2\alpha+l}
(i)=l$ for $l<2$ (we are using the assumption
$\lambda^{<\lambda}=\lambda$). Now we choose
$h_{\eta_{i,2\alpha+l}}$ by induction on $\alpha$. Hence, coming to $\alpha$, let us denote by
$\eta_l$ the sequence $\eta_{i,2\alpha+l}$, and let
$\eta_0\cap\eta_1=\eta$
(so $\eta_0\rest i=\eta_1\rest i=\eta$).
Let $M_{\delta_{i+1}}\setminus
M_{\delta_i}=\{d^i_j:\,j<j_i^\ast\}$,
so that $d_0^i=a_{\delta_i}$ and $d_1^i=b_{\delta_i}$. We consider the type
$\Gamma$, which is the union of
\begin{description}
\item{(a)} 
\begin{equation*}
\Gamma_0\deq \left\{
\begin{split}\varphi(x^0_{j_0},\ldots,x^0_{j_{n-1}};
h_{\eta}(\bar{c})):\, & n<\omega\,\,\&\,\,\bar{c}\subseteq M_{\delta_i}
\,\,\&\,\, j_0,\ldots j_{n-1}<j^\ast_i
\,\,\& \\
&M_{\delta_{i+1}}\models
\varphi(d^i_{j_0},\ldots,d^i_{j_{n-1}};\bar{c})\\
\end{split}
\right\},
\end{equation*}
{\scriptsize{(taking care of one ``side" (for $\eta_0$ or
$\eta_1$) of the part $(\alpha)$
above)}}

\item{(b)} $\Gamma_1$, defined analogously to $\Gamma_0$, but with
$x^0_{j_0},\ldots,x^0_{j_{n-1}}$ replaced everywhere by
$x^1_{j_0},\ldots,x^1_{j_{n-1}}$,

{\scriptsize{(taking care of the remaining ``side" of $(\alpha)$ above,
interchanging $\eta_0$ and $\eta_1$)}}

\item{(c)} $\Gamma_2=\{ (x^0_0, x^0_1)_I\cap (x^1_0, x^1_1)_I=\emptyset\},
$

{\scriptsize{(this says that
the above intervals in $<_0$
are disjoint, which after the right choice of
$h_{\eta_0}(d^i_j)=$ a realisation of $x^0_j$ or
$h_{\eta_0}(d^i_j)=$ a realisation of $x^1_j$ ($j<2$),
and similarly for $h_{\eta_1}$, will take care of the
part $(\delta)$ above.)}}

\item{(d)} $\Gamma_3=\Gamma_3^0\cup \Gamma_3^1$, where for $l<2$
\[
\Gamma_3^l=\left\{x_j^l\neq c:\, j<j_i^\ast, c\in
 \cup\{\Rang(h_\rho):\,
h_\rho\mbox{ already defined}\}
\right\}.
\]

\item{(e)} 
\[
\Gamma_4=\{x^0_{j_1}\neq x^1_{j_2}:\,j_0,j_1<j_i^\ast\}
\]
{\scriptsize{((d)+ (e) are taking care of ($\gamma$) above, part (i)).}}

\item{(f)} 
\begin{equation*}
\Gamma_5=
\left\{
\begin{split}
\neg(x^0_{j_0}E x^1_{j_1}):\,&\mbox{if }j_0,j_1<j_i^\ast\\
&\mbox{but there is no
} a\in M_{\delta_i}\mbox{ with }[d^i_{j_0} E a
\,\,\,\vee\,\,d^i_{j_1} E a]\\
\end{split}
\right\}.
\end{equation*}
{\scriptsize{(together with $\Gamma_6$ below, taking care of part $(\gamma)$(ii),
see below. Note that the type is defined using $\vee$ rather than
$\wedge$, but this will turn out to be sufficient.) }}

\item{(g)} $\Gamma_6=\Gamma_6^0\cup \Gamma^1_6$, where
\begin{equation*}
\Gamma_6^l=
\left\{
\begin{split}
\neg(x^l_{j}E b):\,&j<j_i^\ast
\mbox{ and }b\mbox{ is an element of the set}\\
&\cup\{\mbox{Rang}(h_\rho):\,h_\rho\mbox{ already defined and  }
\neg(\exists c\in N_\eta)[b\,E\,c]\}\\
\end{split}
\right\}.
\end{equation*}

\end{description}
First note that requiring $\Gamma_5\cup \Gamma_6$ throughout the construction
indeed guarantees that $(\gamma)$(ii) can be satisfied. Namely, suppose that
the realisations of $x^0_{j_0}$ and $x^1_{j_1}$ are $E$-equivalent. Then by $\Gamma_5$
we must have that for some $l<2$ and $a\in M_{\delta_i}$ we have that
$d^i{j_l}\,E a$. By transitivity then the realisation of $x^{1-l}_{j_{1-l}}$
would have to be $E$-equivalent to $h_\eta(a)$, which might only be precluded by
$d^i_{j_{1-l}}$ being $E$ equivalent to some $c$ such that $a$ and $c$ are
not $E$-equivalent. This cannot happen by $\Gamma_6$.


We conclude that,
if $\Gamma$ is consistent, as $\mathfrak
C$ is $\bar{\kappa}$-saturated, the functions $h_{\eta_l}$ can be defined. 
Namely, for a realisation $\{c^l_j:\,j<j_i^\ast,l<2\}$ of $\Gamma$,
we can define $g_l(d^i_j)=c^l_j$, and then we let
$h_{\eta_0}=g_0$ if $c^0_1<_0 c^1_0$, and $g_1$ otherwise. We let
$h_{\eta_1}=g_{1-l}$ if $h_{\eta_0}=g_l$.
This guarantees that $(\delta)$ above is satisfied.

Let us then show that $\Gamma$ is consistent. Suppose
for contradiction that this is not so,
so we can find a finite $\Gamma'
\subseteq \Gamma$ which is inconsistent. 
Let $\{j_0, \ldots, j_{n-1}\}$ be an increasing enumeration
of a set including all $j<j_i^\ast$ such that $x^l_j$ is mentioned
in $\Gamma'$ for some $l<2$ and let $\bar{d}=\langle
d_{j_0}^i,\ldots d_{j_{n-1}}^i\rangle$. 
Without loss of generality, 0 and 1 appear in the list
$\{j_0, \ldots, j_{n-1}\}$ and hence
$j_0=0$ while $j_1=1$. By closing under
conjunctions and increasing $\Gamma'$ (retaining that $\Gamma'\subseteq
\Gamma$ is finite) if necessary, we may
assume that for some formula $\sigma(x_0,\ldots x_{n-1};\bar{c})
\in{\rm tp}(\bar{d}/M_{\delta_i})$, we have
\[
\Gamma'\cap\Gamma_l=\{\vartheta_l(x_{j_0}^l,\ldots x_{j_{n-1}}
^l;h_\eta(\bar{c}))\} \]
for $l<2$, where $\vartheta_l(x_{j_0}^l,\ldots x_{j_{n-1}}
^l;h_\eta(\bar{c}))$ is the formula obtained from $\sigma$ by replacing $x_k$
by $x^l_{j_k}$ and $\bar{c}$ by $h_\eta(\bar{c})$.

Let $\vartheta_2$ be the formula comprising $\Gamma_2$ and
$\vartheta^l_3(\bar{x}^l;\bar{c}^l_3)=\bigwedge (\Gamma^l_3\cap \Gamma')$,
while $\vartheta_4=\bigwedge (\Gamma_4\cap \Gamma')$ and
$\vartheta_5=\bigwedge (\Gamma_5\cap \Gamma')$. Let
$\vartheta_3=\vartheta_3^0\wedge\vartheta_3^1$ and $\vartheta=\bigwedge_{k<6, k\neq
2,3}\vartheta_k$.
Without loss of generality, $\vartheta$ includes 
statements $x_0^l\neq\ldots\neq x_{n-1}^l$ and $x^l_0<_0x^l_1$
for $l<2$. We may also assume that $(x^0_0, x^0_1)_I\cap
(x^1_0, x^1_1)_I=\emptyset$ is included in $\Gamma'$. The choice
of $n$ may be assumed to have been such that for some $c^l_0,\ldots
c^l_{n-1}$ (for $l<2$) from $\bigcup\{\Rang(h_\rho):\,
h_\rho\mbox{ already defined}\}$, we have
\[
\Gamma'\cap\Gamma_3=\{x^l_{j_m}\neq c^l_{k}:\,l<2, k<n, m<n\},
\]
and finally that 
\begin{equation*}
\begin{split}
\Gamma'\cap\Gamma_5=\{\neg(x^0_{j_k} E x^1_{j_m}):\,& k,m<n
\,\,\&\,\\
&\neg (\exists a\in M_{\delta_i})
[d^i_{j_k}E a\,\,\vee\,\, d^i_{j_m} E a]\}.\\
\end{split}
\end{equation*}
By extending $h_\eta$ to an automorphism $\hat{h}_\eta$ of $\mathfrak C$,
and applying $(\hat{h}_\eta)^{-1}$, we may assume that
$h_\eta={\rm id}_{M_{\delta_i}}$. We can also assume that
no $c^l_{k}$ is an element of $M_{\delta_i}$, as otherwise
the relevant inequalities can be absorbed by $\sigma$. 

We shall use the following general

\begin{Fact}\label{cinjenica} Suppose that $N\elementary \mathfrak C$ and
$\bar{e}\in 
{}^m\mathfrak C$ is disjoint from $N$, while $N\subseteq A$.

\underline{Then}
\begin{equation*}
\begin{split}
r(\bar{x})\deq &{\rm tp}(\bar{e}, N)\cup\{x_k\neq d:\,d\in A\setminus N, k<m
\}\\ &\cup\{\neg(x_k E d):\,d\in A\,\,\&\,\, (d/E)\cap N=\emptyset, k<m\},\\
\end{split}
\end{equation*}
is consistent (and in fact, every finite subset of it is realised in $N$).
\end{Fact}

\begin{Proof of the Fact} Otherwise, there is a finite
$r'(\bar{x})\subseteq r(\bar{x})$ which is inconsistent. Without
loss of generality, $r'(\bar{x})$ is the union of sets of the following
form (we have a representative type of the sets for each clause)
\begin{itemize}
\item $\{\varrho(\bar{x},\bar{c})\}$ for some $\bar{c}\subseteq N$
and $\varrho$ such that $\models\varrho[\bar{e},\bar{c}]$.
\item $\{x_k\neq \hat{c}_k^s:\,k<m \}$ for some $\hat{c}_0^s,\ldots
\hat{c}^s_{m-1}\in A\setminus N$ and $s<s^\ast<\omega$,
\item $\{\neg (x_k\,E\, \hat{d}_k^t):\,k<m\}$ for some $\hat{d}_0^t,\ldots
\hat{d}^t_{m-1}\in A\setminus N$ and $t<t^\ast<\omega$ such that $(\hat{d}^t_k/E)\cap
N=\emptyset$. 
\end{itemize}
By the elementarity of $N$, there is $\bar{e}'\in N$ with $N\models
\varrho[\bar{e}',\bar{c}]$. By the choice of the rest of the formulae
in $\bar{r}'(\bar{x})$, we see that $\bar{e}'$ satisfies them as well,
which is a contradiction.
$\eop_{\ref{cinjenica}}$
\end{Proof of the Fact}

Let $\bar{x}^l=(x^l_{0},\ldots, x^l_{{n-1}})$ for $l<2$. Let
\[
\Phi_0\deq\{\varphi(\bar{x}^0):\,\varphi(x^0_{j_0},\ldots,x^0_{j_{n-1}})\in
\Gamma'\cap (\Gamma_0\cup\Gamma_3^0\cup \Gamma_6^0)\}.
\]
Applying the last phrase in the above Fact
to $\Phi_0(\bar{x}^0)$, the model
$M_{\delta_i}$ and $\bar{d}$, we obtain a sequence
$\bar{e}^0=(e^0_0, \ldots
e^0_{n-1})\in M_{\delta_i}$ which realises $\Phi_0(\bar{x}^0)$. For
$k,m$ such that $\neg(x^0_{j_k} E x^1_{j_m})\in \Gamma_5$ we have
$\neg(\exists a\in M_{\delta_i})(a E
d^i_{j_k}\vee a E d^i_{j_m})$. So
\[
\neg(x_{k} E e^0_m)\wedge\neg (e^0_m\,E\,x_m)
\in {\rm tp}(\bar{d}/
M_{\delta_i}). 
\]
Let now 
\begin{equation*}
\begin{split}
\Phi_1(\bar{x}^1)=\{\neg(x_{k}^1\,E\,e^0_m)\wedge
\neg(e^0_{k}\,E\,x^1_{m}):\,\neg(x^0_{j_k}\,E\,x^1_{j_m})
\in
\Gamma_5\}\\\cup\{x_{k}^1\neq e^0_m:\,k,m<n\}\cup
\{\varphi(\bar{x}^1):\,\varphi(x^1_{j_0},\ldots, x^1_{j_{n-1}})\in
\Gamma'\cap(\Gamma_1\cup\Gamma^1_3\cup \Gamma_6^1)\}.\\ \end{split}
\end{equation*}
$\Phi_1(\bar{x}^1)$ is a finite set
to which we can apply the last phrase of Fact \ref{cinjenica}. In this way
we find $\bar{e}^1=(e^1_0,\ldots e^1_{n-1})\in M_{\delta_i}$
realising $\Phi_1(\bar{x}^1)$. Now we show that
$\bar{e}^0\frown\bar{e}^1$ realises $\Gamma'\setminus\Gamma_2$. So suppose
$\neg(x^0_{j_k} E x^1_{j_m})\in \Gamma'\cap \Gamma_5$,
then $\neg(x_k^1\,E\,e^0_m)\in \Phi_1$, hence
$\neg (e^1_k E
e^0_m)$. Also $\wedge_{k,m<n} (e^1_k\neq e^0_m)$ holds, by the choice of
$\Phi_1$, so $\bar{e}^0\frown\bar{e}^1$ realises $\Gamma'\cap\Gamma_4$.
Now we need to deal with $\Gamma_2$.
Let 
\[
\DD\deq\{(\bar{u}^0,\bar{u}^1):\,(\bar{u}^0,\bar{u}^1)
\mbox{ satisfies } \vartheta\}.
\]
So $\DD$ is first order definable with parameters in $M_{\delta_i}$
and we have just shown that $\DD\cap M_{\delta_i}\neq\emptyset$. 
Also if $\bar{e}^0\frown\bar{e}^1\subseteq M_{\delta_i}$ satisfies $\vartheta$,
it necessarily realises $\Gamma'\setminus
\Gamma_2$ (as no $c^l_k\in M_{\delta_i}$, see the definition). As
$\Gamma'$ is presumed to be inconsistent, no $(\bar{u}^0,\bar{u}^1)
\in\DD\cap M_{\delta_i}$ can realise $\Gamma'$, i.e. satisfy $\vartheta_2$, and
hence for no $(\bar{u}^0,\bar{u}^1)\in \DD\cap M_{\delta_i}$ are
the intervals $(u^0_0, u^0_1)_I$
and $(u^1_0, u^1_1)_I$ disjoint. Now we
claim that if $(\bar{u}^0,\bar{u}^1)
\in M_{\delta_i}\cap\DD$, then $(u^0_0, u^0_1)_I\cap (a_{\delta_i},
b_{\delta_i})_I \neq\emptyset$. 

Indeed, suppose otherwise, say $u^0_1<_0 d^i_0=a_{\delta_i}$, so
$u^0_1<_0 x_0\in {\rm tp}(\bar{d}/M_{\delta_i})$.
Arguing as above, with $\bar{u}^0$ in place of $\bar{e}^0$ and
$\Phi_1(\bar{x})\cup\{u^0_1<_0 x_0^1\}$ in place of
$\Phi_1(\bar{x}^1)$, we can find 
$\bar{u}\in M_{\delta_i}$ satisfying
$
(u^0_1<_0 u_0)$ and such that $(\bar{u}^0,\bar{u})$
satisfies $\vartheta$. So $(\bar{u}^0, \bar{u})\in
\DD\cap M_{\delta_i}$ and the intervals $(u^0_0,u^0_1)_I$ and $(u_0,u_1)_I$ are
disjoint, a contradiction. A similar contradiction can be derived from the assumption
$b_{\delta_i}=d^{\delta_i}_1<_0 u^0_0$. Now note that $(\bar{u}^0,\bar{u}^1)
\in\DD\implies(\bar{u}^1,\bar{u}^0)
\in\DD$, so if $(\bar{u}^0,\bar{u}^1)
\in\DD\cap M_{\delta_i}$ we also have $(u^1_0,u^1_1)_I\cap(a_{\delta_i},
b_{\delta_i})_I\neq\emptyset.$

By the choice of $C$, there is no $x\in M_{\delta_i}$ with $d^{\delta_i}_0\le
_0 x\le _0 d^{\delta_i}_1$,
hence
 \begin{equation*}
\mbox{if }(\bar{u}^0,\bar{u}^1)
\in M_{\delta_i}\cap\DD\mbox{ and }l<2,\mbox{ we have }u^l_0<_0d^{\delta_i}_0
<_0d^{\delta_i}_1<_0 u^l_1.
\tag*{$(\ast)$}
\end{equation*}
Let $\sigma^\ast(\bar{y})$
be $\exists \bar{x} ((\bar{x},\bar{y})\in \DD)$. Hence if
\[
\varrho_0(z)=(\exists \bar{y})[\sigma^\ast(\bar{y})\wedge z\le_0 y_0]
\]
and
\[
\varrho_1(z)=(\exists \bar{y})[\sigma^\ast(\bar{y})\wedge y_1\le_0 z],
\]
then 
\[
M_{\delta_i}\models(\forall z_0,z_1)[
\varrho_0(z_0)\wedge 
\varrho_1(z_1)\implies z_0<_0 z_1]. \]
Of course, this holds in ${\mathfrak C}$ as well, so
\begin{description}
\item{(a)} $\varrho_0(z)$ defines an initial segment of $I$,
\item{(b)} $\varrho_1(z)$ defines an end segment of $I$,
\item{(c)} the segments defined by $\varrho_0(z)$ and $\varrho_1(z)$ are disjoint,
\item{(d)} $\varrho_0(M_{\delta_i})\cup \varrho_1(M_{\delta_i})=I\cap M_{\delta_i}$.

[Why? Note that $(\bar{e}^0,
\bar{d})\in\DD$. Hence
$\sigma^\ast(\bar{d})$ holds. As for every $a\in I\cap M_{\delta_i}$
we have $a<_Ia_{\delta_i}$ or $a>_I b_{\delta_i}$, the conclusions follows.]

\item{(e)} $\varrho_0(a_{\delta_i})$ and $\varrho_1(b_{\delta_i})$ hold.

[Why? Again because $\sigma^\ast(\bar{d})$ holds.]
\end{description}

The above arguments show that $\{x:\,(\exists \bar{y})[(\sigma^\ast(
\bar{y})\wedge x<_0 y_0)]\}$ defines the Dedekind cut $\{x:\,x<_0a_{\delta_i}
\}$ over $M_{\delta_i}$, which contradicts the choice of $C$ and the fact that
the Dedekind cut induced by $(\bar{a},\bar{b})$ is not definable (which is a
part of the definition of $\ast[M,\bar{a},\bar{b}]$).

{\bf Stage }{\bf C}.
Now we have shown that
the trees $\langle N_\eta:\,\eta\in {}^{\lambda>}2\rangle$,
$\langle h_\eta:\,\eta\in {}^{\lambda>}2\rangle$
of models and embeddings can be defined as required, and we consider
\[
p^\ast\deq\cup_{\eta\in {}^{\lambda>}2} h_\eta(p\rest M_{\delta_{\llg(\eta)}}).
\]
We shall show that $p^\ast$
is finitely satisfiable, hence satisfiable.
Let $\Gamma'\subseteq p^\ast$ be finite.
Recalling the analysis of $p$ from Stage A, we consider each of the cases
by which $p$ could have been defined (ignoring Case 1A, as justified by
Subclaim \ref{addition}.)

\underline{Case 1}. In this case there is a function $f:\,Q^M
\into Q^M$ respecting $E$, with $a E f(a)$ for all $a\in Q^M$, and
without loss of generality
there are some $\eta_0,\ldots \eta_{m-1}\in {}^{\lambda>}2$
and $\{b^j_i:\,j<m, i<n_j\}
\subseteq \Rang(f)$ such that 
\[
\Gamma'=\{P(z)\}\cup\bigcup_{j<m}\{h_{\eta_j}(b^j_i) R\, z:\,i<n_j\},
\]
and for each $j$ we have $\{ b^j_i:\,i<n_j\}\subseteq M_{\delta_{\llg(\eta_j
)}}$. Let $n\deq\Sigma_{j<m} n_j$, hence $\Gamma'$ is a quantifier
free (partial) type over $n$ variables in $\mathfrak C^{[\bar{\psi}]}$.
By Observation \ref{eliminacija}, we only need to check that in
$\Gamma'$ there are no direct contradictions with the axioms of
$T_{\rm feq}^+$.

The only possibility for such a contradiction is that for some $j_0, j_1$ and
$b^{j_0}_i$, $b^{j_1}_k$ we have
\[
h_{\eta_{j_0}}(b^{j_0}_i)\neq h_{\eta_{j_1}}(b^{j_1}_k) \wedge
h_{\eta_{j_0}}(b^{j_0}_i) E\, h_{\eta_{j_1}}(b^{j_1}_k)
\]
and $h_{j_0}(b^{j_0}_i)R\, z, h_{j_1}(b^{j_1}_k)R\, z \in \Gamma'$.
In such a case, any $c$ which would realise $\Gamma'$ would contradict
part (c) of the definition of $T_{\rm feq}^+$. Suppose that $b_i^{j_0},
b^{j_1}_k$ and $\eta_0,\eta_1$ are as above. Let $\eta_l\deq\eta_{j_l}$ for
$l<2$ and let $\eta=\eta_0\cap\eta_1$. By part $(\gamma)$(ii) in the definition
of $\bar{h}$, there is $\hat{b}\in N_\eta$ such that $h_{\eta_0}(b^{j_0}_i)
E\, \hat{b}$ and $h_{\eta_1}(b^{j_1}_k) E\, \hat{b}$. For some $b\in 
M_{\delta_{\llg(\eta)}}$ we have
\[
h_{\eta_0}(b)=h_{\eta_1}(b)=h_\eta(b)=\hat{b},
\]
so applying the elementarity of the maps, we obtain
\[
b^{j_0}_i E\, b E\, b^{j_1}_k.
\]
On the other hand,
by the definition of $p^\ast$ we have $b^{j_0}_i R\,z\in p(z)$ and
$b^{j_1}_k R\, z\in p(z)$. By the the demands on $p$
this implies that $b^{j_0}_i=b^{j_1}_k\notin M_{\delta_{\llg(\eta)}}$ and
$f(b)= b^{j_0}_i$, contradicting the fact that $M_{\delta_{\llg(\eta)}}$ is closed
under $f$.

\underline{Case 2}. For a fixed $a^\ast\in M_0$ we have
\[
p(z)=\{Q(z)\}\cup\{a^\ast E z\}
\cup\{z\neq c:\,c\in a^\ast/E^M\},
\]
so without loss of generality
\[
\Gamma'=\{Q(z)\}\cup\{a^\ast\,E\,z\}\cup\{z\neq h_{\eta_j}(c_j):\,j<m\}
\]
for some $c_0,\ldots,c_{m-1}\in a^\ast/E^M$ and $\eta_0,\ldots,\eta_{m-1}
\in {}^{\lambda>}2$,
as $h_{\langle\rangle}=
{\rm id}_{M_0}$. As $a^\ast/E$ is infinite in any model of $T^\ast_{\rm feq}$,
the set $\Gamma'$ is consistent.

\underline{Case 3}. We may assume that for some
equivalence relation ${\cal E}$ on $P^M$, a function $f$ from
$P^M$ into $\{{\rm yes}, {\rm no}\}$, sequences $\eta_0,\ldots
\eta_{n-1} \in {}^{\lambda>}2$, and $\{a^k_i:\,i<m,k<n\}\subseteq Q^M$ and
$\{b^k_i,c^k_i,d^k_i:\,i<m,k<n\}\subseteq P^M$ we have $e_1\,\EE\,e_2\implies
f(e_1)=f(e_2)$ and
\begin{equation*}
\begin{split}
\Gamma'(z)=&\{Q(z)\}\cup\bigcup_{k<n}\{\neg (h_{\eta_k}(a^k_i) E z):\,i<m\}
\cup\bigcup_{k<n}\{(zR\,h_{\eta_k}(b^k_i))^{f(b^k_i)}:\,i<m\}\\
&\cup \bigcup_{k<n}\{[F(z,h_{\eta_k}(c_i^k))=F(z,h_{\eta_k}(d_i^k))]^{{\rm if
}c_i^k\EE d_i^k}: \,i<m\}.\\
\end{split}
\end{equation*}
We could have a contradiction if for some $k_1,k_2, i_1, i_2$ we had
$f(b^{k_1}_{i_1})=$yes, $f(b^{k_2}_{i_2})=$no, but
$h_{\eta_{k_1}}(b^{k_1}_{i_1})=
h_{\eta_{k_2}}(b^{k_2}_{i_2})$, which cannot happen by $\gamma$(i) and the fact
that each $h_\eta$ is 1-1. Another
possibility is that for some $b^{k_1}_{i_1},b^{k_2}_{i_2}$ we have
$f(b^{k_1}_{i_1})
=f(b^{k_2}_{i_2})$=yes, but $h_{\eta_{k_1}}(b^{k_1}_{i_1})\neq
h_{\eta_{k_2}}(b^{k_2}_{i_2})$ while $h_{\eta_{k_1}}(b^{k_1}_{i_1})\,E\,
h_{\eta_{k_2}}(b^{k_2}_{i_2})$. To see that this cannot happen, we distinguish
various possibilities for $b^{k_1}_{i_1},b^{k_2}_{i_2}$ and use part
$(\gamma)$(ii) in the choice of $\bar{h}$.

Yet another possible source of
contradiction could come from a similar consideration involving the last clause
in the definition of $\Gamma'(z)$, which cannot happen for similar reasons.

{\bf Stage} {\bf D}.
Now we can conclude, using $\lambda=\lambda^{<\lambda}$ and $\card{T}<\lambda$,
that there is a model $N^\ast\elementary \mathfrak C$ of size $\lambda$ with
$\bigcup_{\eta\in {}^{\lambda>}2} N_\eta\subseteq N^\ast$, such that $p^\ast$ is
realised in $N^\ast$. For $\nu\in {}^{\lambda}2$, let
 $h_\nu\deq\cup_{i
<\lambda}h_{\nu\rest i}$, and let $N_\nu\deq\Rang(h_\nu)$, while $p_\nu
\deq h_\nu(p)$.

For such $\nu$, let
\[
q_\nu(x)\deq\{I(x)\} \cup\{h_\nu(a_i)<_0 x<_0 h_\nu(b_i):\,i<\lambda\}.
\]
Hence we have that for $\nu\neq \rho$
from ${}^{\lambda} 2$, the types $q_\nu$ and $q_\rho$ are contradictory,
by $(\delta)$ above. As $\norm{N^\ast} +\card{L(T)}\le\lambda$,
there are only $\le\lambda$ definable Dedekind
cuts of $<_0$ over $N^\ast$, and only
$\le \lambda$ types $q_\nu$ are realised in $N^\ast$.
Hence there is $\nu\in {}^{\lambda} 2$
(actually $2^\lambda$ many) such that the Dedekind cut
$\{x:\,\vee_{i<\lambda} x<_0 h_\nu(a_i)\}$ is not definable over
$N^\ast$ and $q_\nu$ is not realised in $N^\ast$. So $N^\ast$ omits
$q_\nu$ and realises $p_\nu$. We let $N=
h(N^\ast)$, where $h$ is an automorphism of $\mathfrak C$ extending
$h^{-1}_\nu$. 
$\eop_{\ref{main}}$
\end{Proof of the Main Claim}

\begin{Theorem}\label{notmax} Assume that $\lambda^{<\lambda}=\lambda$
and $2^\lambda=\lambda^+$.

{\noindent (1)} For any $\lambda$-relevant $(T_{\rm ord}, T^\ast_{\rm feq})$
-superior $(T,\bar{\varphi},\bar{\psi})$, the theory $T$
has a model $M^\ast$ of cardinality $\lambda^+$
such that
\begin{description}
\item{(i)} $\bar{\varphi}^{M^\ast}$ is not $\lambda^+$-saturated,
\item{(ii)} $\bar{\psi}^{M^\ast}$ is $\lambda^+$-saturated.
\end{description}

{\noindent (2)} We can strengthen the claims in (i) and (ii) to include any
interpretations
of a dense linear order and $T_{\rm feq}^\ast$-respectively in $M^\ast$,
even with parameters.
\end{Theorem}

\begin{Proof} We prove (1), and (2) is proved similarly.
Using the Main Claim
\ref{main}, we can construct $M^\ast$ of size $\lambda^+$,
by an $\elementary$-increasing continuous sequence $\langle M_i^\ast:\,i\le
\lambda^+ \rangle$, with $\norm{M_i^\ast}=\lambda$
satisfying
$\ast[M_i,\bar{a},\bar{b}]$ for each $i\le\lambda^+$,
and letting
$M^\ast=M_{\lambda^+}$. The Main Claim
\ref{main} is used in the successor steps. To assure that $M^\ast$ is 
$\lambda^+$-saturated for $T^\ast_{\rm feq}$,
we use the assumption $2^\lambda=\lambda^+$, to do the bookkeeping of all
$T_{\rm feq}^\ast$-types involved.
$\eop_{\ref{notmax}}$
\end{Proof}

\begin{Conclusion} Under the assumptions of Theorem \ref{notmax},
the theory
$T_{\rm feq}^\ast$ is $\initial^\ast_{\lambda^+}$ strictly below
the theory of a dense linear order with no first or last elements.
\end{Conclusion}

[Why? It is below by Shelah's Theorem \ref{dwaprime} above.]

\medskip

We recall that our motivation for studying $\initial^\ast$
is to try to characterise SOP${}_3$ (or SOP${}_2$) theories by the $\initial^\ast$-maximality.
As we explained in the Introduction this has origins in the connection between the
maximality in the Keisler order and having the strict order property, so we should
show here what is the connection between the maximality in Keisler's order and the
maximality in the order
$\initial^\ast$. The following Claim \ref{connection} does that for countable theories.

\begin{Claim}\label{connection} Suppose that $T$ is a countable theory that is
$\vartriangleleft_{\lambda^\ast}^\ast$-maximal.
\underline{Then} it is
maximal in the
Keisler order $\circledvartriangleleft_\lambda$. 
\end{Claim} 

\begin{Proof of the Claim} Suppose otherwise and let $T_1$ be a theory that is
$\vartriangleleft^\ast$-maximal but not
maximal in the Keisler order $\circledvartriangleleft_\lambda$. In particular we
have $T_{\rm ord}\vartriangleleft_\lambda^\ast T_1$, so there is a $\lambda$-relevant
$(T_{\rm ord}, T_1)$-superior triple $(T,\bar{\varphi_0}, \bar{\varphi_1})$-exemplifying
this. By Observation \ref{trivint}(0) we may assume that the interpretation $\bar{\varphi_1}$
is trivial, so $T_1\subseteq T$- for simplicity.
 
Since $T$ is not maximal in the
Keisler order $\circledvartriangleleft_\lambda$, by \cite{Sh c} 4.2 (1) there
is a regular ultrafilter $\DD$ which is not good and a model $M$ of $T$ such that
$M^\lambda/\DD$ is nevertheless $\lambda^+$-compact. We can extend $M$ to a model $N$ of
$T$ and consider $N^\ast=N^\lambda/\DD$. This is a model of $T$ and by the Extension Theorem
for ultrafilters we have that $[N^\ast]^{\bar{\varphi_1}}=M^\lambda/\DD$, so it is
$\lambda^+$-compact and hence it is $\lambda^+$-saturated. Again by the Extension Theorem
we have that $[N^\ast]^{\bar{\varphi_1}}=(N^{\bar{\varphi_1}})^\lambda/\DD$. Now on the one hand
we have by the $\vartriangleleft_{\lambda^\ast}^\ast$-maximality of $T_1$ that
$(N^{\bar{\varphi_1}})^\lambda/\DD$ must be $\lambda^+$-saturated, hence
$\lambda^+$-compact. But on the other hand $(N^{\bar{\varphi_1}})^\lambda/\DD$ cannot
be $\lambda^+$-compact because $\DD$ is not a good ultrafilter and $T_{\rm ord}$ is
maximal in the Keisler order, contradicting \cite{Sh c} 4.2 (1).
$\eop_{\ref{connection}}$
\end{Proof of the Claim}



\section{On the properties SOP${}_2$ and SOP${}_1$}\label{druga}
In his paper \cite{Sh 500}, S. Shelah investigated a hierarchy
of properties unstable theories without strong order property may have.
This hierarchy is named NSOP${}_n$ for $3\le n<\omega$, where the
acronym NSOP stands for ``not strong order property". The
negation of NSOP${}_n$ is denoted by SOP${}_n$. It was shown in
\cite{Sh 500} that SOP${}_{n+1}\implies$ SOP${}_n$,
that the implication is strict and that
SOP${}_3$ theories are not simple. 
In this section we investigate two further notions, which with the
intention of furthering the above hierarchy, we name
SOP${}_2$ and SOP${}_1$. The original definitions of SOP${}_n$ for
$n\ge 3$ do not immediately lend themselves to extending the hierarchy for $n=1,2$,
but the properties we define nevertheless fulfill that role.
In section \ref{ekvivalencija}, a connection between this
hierarchy and $\initial^\ast_\lambda$-maximality will be established.

Recall from \cite{Sh 500} one of the equivalent definitions
of SOP${}_3$. (The equivalence is established in Claim 2.19
of \cite{Sh 500}).

\begin{Definition}\label{sop3} (1) A (complete) theory $T$ has SOP${}_3$
iff there is an indiscernible sequence
$\langle \bar{a}_i:\,i<\omega\rangle$ and formulae
$\varphi(\bar{x},\bar{y})$, $\psi(\bar{x},\bar{y})$ such that
\begin{description}
\item{(a)} $\{\varphi(\bar{x},\bar{y}), \psi(\bar{x},\bar{y})\}$ is contradictory,
\item{(b)} for some sequence $\langle \bar{b}_j:\,j<\omega\rangle$ we have
\[
i\le j\implies \models\varphi [\bar{b}_j,\bar{a}_i]\mbox{ and }
i> j\implies \models\psi [\bar{b}_j,\bar{a}_i],
\]
\item{(c)} for $i<j$, the set $\{\varphi(\bar{x},\bar{a}_j),
\psi(\bar{x},\bar{a}_i)\}$ is contradictory.
\end{description}

{\noindent (2)} NSOP${}_3$ stands for the negation of SOP${}_3$.

\end{Definition}


\begin{Definition}\label{nsop21}
{\noindent (1)} $T$ has SOP${}_2$ if there is a formula $\varphi(\bar{x},
\bar{y})$ which exemplifies this property in
$\mathfrak{C}=\mathfrak{C}_T$, which means:

There are $\bar{a}_\eta\in {\mathfrak C}$ for $\eta\in {}^{\omega>}2$
such that
\begin{description}
\item{(a)} for every
$\rho\in {}^{\omega}2$, the set
$\{\varphi(\bar{x}, \bar{a}_{\rho\rest n}):\,n<\omega\}$ is
consistent,
\item{(b)} if $\eta,\nu\in {}^{\omega>}2$ are incomparable,
$\{\varphi(\bar{x},\bar{a}_\eta),
\varphi(\bar{x},\bar{a}_\nu)\}$ is inconsistent.
\end{description}

{\noindent (2)} $T$ has SOP${}_1$ if there is a formula $\varphi(\bar{x},
\bar{y})$ which exemplifies this in ${\mathfrak C}$, which means:

There are $\bar{a}_\eta\in {\mathfrak C}$, for $\eta\in {}^{\omega>}2$ such that:

\begin{description}
\item{(a)} for $\rho\in {}^\omega 2$ the set
$\{\varphi(\bar{x},\bar{a}_{\rho\rest n}):
\,n<\omega\}$ is consistent.

\item{(b)} if $\nu\frown\langle 0\rangle\initialeq \eta\in {}^{\omega>}2$,
\underline{then} $\{\varphi(\bar{x},\bar{a}_\eta), \varphi(\bar{x},
\bar{a}_{\nu\frown\langle 1\rangle})\}$ is inconsistent.

\end{description}

{\noindent (3)} NSOP${}_2$ and NSOP${}_1$ are the negations of SOP${}_2$ and
SOP${}_1$ respectively.

\end{Definition}

The following Claim establishes the relative position of the
properties introduced in Definition \ref{nsop21} within the (N)SOP hierarchy.

\begin{Claim}\label{position} For any complete first order theory $T$, we have
\[
\mbox{SOP${}_3\implies$SOP${}_2\implies$SOP${}_1$}.
\]
\end{Claim}

\begin{Proof of the Claim} Suppose that $T$ is SOP${}_3$, as exemplified
by $\langle \bar{a}_i:\,i<\omega\rangle$, $\langle
\bar{b}_j:\,j<\omega\rangle$ and formulae $\varphi(\bar{x},\bar{y})$ and
$\psi(\bar{x},\bar{y})$ (see Definition \ref{sop3}), and
we shall show that $T$ satisfies SOP${}_2$. We define \[
\vartheta(\bar{x},\bar{y}^0\frown
\,\bar{y}^1)\equiv\varphi(\bar{x},\bar{y}^0)\wedge
\psi(\bar{x},\bar{y}^1), \mbox{ where }\llg(\bar{y}^0)= \llg(\bar{y}^1).\]
Let us first prove the consistency of
\begin{equation*}
\Gamma\deq 
\begin{split}T&\cup\{\neg(\exists\bar{x})[\vartheta(\bar{x},
\bar{y}_\eta)\wedge \vartheta(\bar{x}, \bar{y}_\nu)]:\,\eta\perp\nu\mbox{ in
}{}^{\omega>}2\} \cup\\
&\cup\bigcup_{n<\omega}\{(\exists\bar{x})
[\bigwedge_{k\le n}\vartheta(\bar{x},\bar{y}_{\eta\rest k})]:\, \eta\in {}^n
2\}.\\ \end{split}
\end{equation*}
Suppose for contradiction
that $\Gamma$ is not consistent, then for some $n<\omega$, the following set
is inconsistent:
\begin{equation*}
\Gamma'\deq 
\begin{split}
T &\cup\{\neg(\exists\bar{x})[\vartheta(\bar{x}, \bar{y}_\eta)\wedge
\vartheta(\bar{x}, \bar{y}_\nu)]:\,\eta,\nu\mbox{ incomparable in
}{}^{n\ge}2\}\\
& \cup\{(\exists\bar{x})[\bigwedge_{k\le n}\vartheta
(\bar{x},\bar{y}_{\eta\rest k})]:\, \eta\in {}^n 2\}.\\
\end{split}
\end{equation*}
Fix such $n$. We pick ordinals $\alpha_\eta,\beta_\eta<\omega$ for
$\eta\in {}^{n\ge} 2$ so that
\begin{description}
\item{(i)} $\nu\initial\eta\implies \alpha_\nu<\alpha_\eta<\alpha_{\eta}+1
<\beta_\eta<\beta_\nu$,
\item{(ii)} $\beta_{\eta\frown\langle 0\rangle}<\alpha_{\eta\,\frown\,
\langle 1\rangle}$. \end{description}
For $\eta\in {}^{n\ge }2$ let $\bar{a}^\ast_\eta\deq\bar{a}_{\alpha_\eta}\frown
\bar{a}_{\beta_\eta}$. We show that $\mathfrak C$ and
$\{\bar{a}^\ast_\eta:\,\eta\in {}^{n\ge} 2\}$ exemplify that $\Gamma'$ is consistent. So, if $\eta\in {}^{n\ge} 2$ then we have
$\bigwedge_{k\le n}\vartheta[\bar{b}_{\alpha_\eta+1}, \bar{a}^\ast_{\eta\rest
k}]$ as for every $k\le n$ we have $\alpha_{\eta\rest k}<\alpha_{\eta}+1$, so
$\varphi[\bar{b}_{\alpha_\eta+1}, \bar{a}_{\alpha_{\eta\rest k}}]$ holds, but
also for all $k\le n$, as $\eta\rest k\initialeq\eta$, we have
$\beta_{\eta\rest k}>\alpha_\eta+1$, so
$\psi[\bar{b}_{\alpha_\eta+1}, \bar{a}_{\beta_{\eta\rest k}}]$ holds.
Hence $(\exists \bar{x})[\bigwedge_{k\le n}\vartheta(\bar{x},
\bar{a}^\ast_{\eta\rest k})]$.
On the other hand, if $\nu\frown\langle l\rangle\initialeq \eta_l$ for $l<2$,
then $\{\vartheta(\bar{x}, \bar{a}^\ast_{\eta_0}), \vartheta(\bar{x},
\bar{a}^\ast_{\eta_1})\}$ is contradictory as the conjunction implies
$\psi(\bar{x}, \bar{a}_{\beta_{\eta_0}})\wedge\varphi(\bar{x}, \bar{a}_{\alpha_{\eta_1}})$,
which is contradictory by $\beta_{\eta_0}<\alpha_{\eta_1}$ and (c) of
Definition \ref{sop3}. This shows that $\Gamma'$ is consistent, hence we have
also shown that $\Gamma$ is consistent.

Having shown that $\Gamma$ is consistent, we can find witnesses
$\{\bar{a}^\ast_\eta:\,\eta\in {}^{\omega>}2\}$ in $\mathfrak C$ realising $\Gamma$. Now we just
need to show that $\{\vartheta(\bar{x}, \bar{a}^\ast_{\eta\rest
n}):\,n<\omega\}$ is consistent for every $\eta\in  {}^\omega 2$. This follows by the compactness theorem and
the definition of $\Gamma$. Hence we have shown that SOP${}_3\implies
{\rm SOP}_2$.

The second part of the claim is obvious (and the witnesses for SOP${}_2$
can be used for SOP${}_1$ as well).
$\eop_{\ref{position}}$
\end{Proof of the Claim}

\begin{Question} Are the implications from Claim \ref{position}
reversible?
\end{Question}

\begin{Claim}\label{infrank}
If $T$ satisfies SOP${}_1$, then $T$ is not simple. In fact,
if $\varphi(\bar{x},\bar{y})$ exemplifies SOP${}_1$ of $T$, then 
the same formula exemplifies that $T$ has the tree property.
\end{Claim}

\begin{Proof of the Claim} 
Let $\varphi(\bar{x},\bar{y})$ and $\{\bar{a}_\eta:\,\eta\in 
{}^{\omega>}2\}$ exemplify SOP${}_1$. Then
\[
\Gamma_\eta\deq\{\varphi(\bar{x},\bar{a}_{\eta\frown\langle
0\rangle_n\frown\langle 1\rangle}):\,n<\omega\}
\]
for $\eta\in {}^{\omega >}2$ consists of pairwise
contradictory formulae. (Here $\langle
0\rangle_n$ denotes a sequence consisting of $n$ zeroes.) For
$n<\omega$ and $\nu\in
{}^n\omega$ let
\[
\rho_\nu\deq\langle 0\rangle_{\nu(0)+1}\frown\langle
1\rangle\frown\langle 0\rangle_{\nu(1)+1}\ldots
\frown\langle 0\rangle_{\nu(n-1)+1}\frown\langle
1\rangle,
\]
so $\rho_\nu
\in {}^{\omega>} 2$ and $\nu\initialeq \eta\implies
\rho_\nu\initialeq\rho_\eta$.
For $\nu\in {}^n\omega$ let $\bar{b}_\nu=\bar{a}_{\rho_\nu}$.  We observe first
that $\{\varphi(\bar{x},\bar{b}_{\nu\,\hat{}\, \langle k\rangle}):\,k<\omega\}$ is
a set of pairwise contradictory formulae, for $\nu\in
{}^n\omega$; namely, if $k_0\neq k_1$, then
$\varphi(\bar{x},\bar{b}_{\nu\frown \langle k_l\rangle})$ for $l<2$ are two
different elements of $\Gamma_{\rho_\nu}$.
On the other hand, $\{\varphi(\bar{x},\bar{b}_{\nu\rest n}):\,n<\omega\}$ is
consistent for every $\nu\in {}^{\omega}\omega$. Hence
$\varphi(\bar{x},\bar{y})$ and $\{\bar{b}_\nu:\,\nu\in {}^{\omega>}\omega\}$
exemplify that $T$ has the tree property, and so $T$ is not simple.
$\eop_{\ref{infrank}}$
\end{Proof of the Claim}

This ends the discussion of the properties of SOP${}_1$ and SOP${}_2$ that are
directly relevant to the main thesis of the paper-the reader
only interested in the connection with the order $\initial^\ast$ can now
turn directly to \S\ref{ekvivalencija}. The rest of this section 
however contains some further syntactic developments of these properties
which are of interest if one wishes to understand the type theory induced
by them. The indescernibility results we have here were recently
used by Shelah and Usvyatsov \cite{ShUs xx} to define a rank function on
NSOP${}_1$ theories (see Theorem \ref{ranksucess})). 

The definition of when a theory has SOP${}_1$ can be made in another
equivalent fashion.

\begin{Definition}\label{another}
Let $\varphi(\bar{x},\bar{y})$ be a formula  of
$\LL(T)$. We say $\varphi(\bar{x},\bar{y})$ has SOP${}'_1$ iff
there is $\langle\bar{a}_\eta:\,\eta\in {}^{\omega>} 2
\rangle$ in ${\mathfrak C}_T$ such that
\begin{description}
\item{(a)} $\{\varphi(\bar{x},\bar{a}_{\rho\rest n})^{\rho(n)}:\,
n<\omega\}$ is consistent for every $\rho\in {}^\omega 2$, where we use the
notation
\[
\varphi^l=\left\{
\begin{array}{cc}
\varphi &\mbox{ if }l=1,\\
\neg\varphi &\mbox{ if }l=0
\end{array}
\right.
\]
for $l<2$.
\item{(b)} If $\nu\frown\langle 0\rangle\initialeq\eta\in {}^{\omega>}2$,
\underline{then} $\{\varphi(\bar{x},\bar{a}_\eta), \varphi(\bar{x},
\bar{a}_\nu)\}$ is inconsistent.
\end{description}

We say that $T$ has property SOP${}'_1$ iff some formula of $\LL(T)$ has it.
\end{Definition}

\begin{Claim}\label{equivalent} 
\begin{description}
\item{(1)} If $\varphi(\bar{x},\bar{y})$
exemplifies SOP${}_1$ of $T$ \underline{then} $\varphi(\bar{x},\bar{y})$ (hence
$T$) 
has property SOP${}'_1$.

\item{(2)} If $T$ has property SOP${}'_1$ \underline{then} $T$ has SOP${}_1$.
\end{description}
\end{Claim}

\begin{Proof of the Claim} (1)
Let $\{\bar{a}_\eta:\,\eta\in {}^{\omega>}2\}$ and
$\varphi(
\bar{x}, \bar{y})$ exemplify that
$T$ has SOP${}_1$. For $\eta\in {}^{\omega>}
2$ we define $\bar{b}_{\eta}\deq\bar{a}_{\eta\frown\langle 1\rangle}$. We shall
show that $\varphi(\bar{x}, \bar{y})$ and $\{\bar{b}_{\eta}:\,\eta\in
{}^{\omega>}2
\}$ exemplify SOP${}'_1$.

Given $\hat{\eta}\in {}^{\omega}2$.
Let $\bar{c}$ exemplify that item (a) from Definition \ref{nsop21}(2) holds
for $\hat{\eta}$.
Given $n<\omega$, we consider $\varphi[\bar{c}, \bar{b}_{\hat{\eta}\rest
n}]^{\hat{\eta}(n)}$. If $\hat{\eta}
(n)=1$, then, as  $\bar{b}_{\hat{\eta}\rest n}=
\bar{a}_{\hat{\eta}\rest n\frown\langle 1\rangle}=\bar{a}_{\hat{\eta}\rest
(n+1)}$, we have that $\varphi[\bar{c}, \bar{b}_{\hat{\eta}\rest
n}]^{\hat{\eta}(n)} =\varphi[\bar{c}, \bar{a}_{\hat{\eta}\rest (n+1)}]$
holds.
If $\hat{\eta}(n)=0$, then 
\[
(\hat{\eta}\rest n)\frown\langle 0\rangle
= \hat{\eta}\rest(n+1).
\]
As $\varphi[\bar{c},\bar{a}_{\hat{\eta}\rest(n+1)}]$ holds,
by (b) of Definition
\ref{nsop21}(2),
we have that
$\varphi[\bar{c}, \bar{a}_{\hat{\eta}\rest n\frown\langle 1\rangle}]$
cannot hold, showing again that, 
$\varphi[\bar{c}, \bar{b}_{\hat{\eta}\rest
n}]^{\hat{\eta}(n)}=\neg\varphi[\bar{c}, \bar{a}_{\hat{\eta}\rest n
\frown\langle 1\rangle}]$ holds.
This shows that $\{\varphi(\bar{x}, \bar{b}_{\hat{\eta}\rest n})^{\hat{\eta}(n)}:\,
n<\omega\}$
is consistent, as exemplified by $\bar{c}$.

Suppose $\nu\frown\langle 0\rangle\initialeq \eta\in {}^{\omega>}2$
and that $\varphi[\bar{d},\bar{b}_\eta]\wedge \varphi[\bar{d},\bar{b}_\nu]$
holds. So both $\varphi[\bar{d},\bar{a}_{\eta\frown\langle 1\rangle}]$ and
$\varphi[\bar{d},
\bar{a}_{\nu\frown\langle 1\rangle}]$ hold. On the other hand, as
$\nu\frown\langle
0\rangle\initialeq\eta$, clearly $\nu\frown\langle
0\rangle\initialeq\eta\,\frown\,\langle
1\rangle$, and so (b) of Definition \ref{nsop21}(2) implies that
$\{\varphi(\bar{x},\bar{a}_{\eta\frown\langle 1\rangle}),
\varphi(\bar{x},
\bar{a}_{\nu\frown\langle 1\rangle})\}$ is contradictory, a contradiction.
Hence  the set $\{\varphi(\bar{x},\bar{b}_{\eta}),
\varphi(\bar{x},
\bar{b}_{\nu})\}$ is contradictory

{\noindent (2)} Define first for $\eta\in {}^{\omega\ge}2$ an element
$\rho_\eta \in {}^{\omega\ge}2$ by letting
\[
\rho_\eta (3k)=\eta(k),
\]
\[
\rho_\eta (3k+1)=0,
\]
\[
\rho_\eta (3k+2)=1,
\]
and if $\llg(\eta)=m<\omega$, then $\llg(\rho_\eta)=3m$. Note that for
$\eta\in {}^\omega 2$
and $k<\omega$ we have $\rho_{\eta\rest k}=\rho_\eta\rest (3k)$. 

Let $\varphi(\bar{x},\bar{y})$
and $\{\bar{a}_\eta:\,\eta\in {}^{\omega>}2\}$ exemplify property SOP${}'_1$.
We
pick $c_0\neq c_1$
and define for $\eta\in {}^{\omega>}2$
\[
\bar{b}_{\eta\frown\langle 1\rangle}\deq \bar{a}_{\rho_\eta}\frown
\bar{a}_{\rho_\eta\frown\langle 1\rangle}
\frown \langle c_0,c_1\rangle,
\]
\[
\bar{b}_{\eta\frown\langle 0\rangle}\deq \bar{a}_{\rho_\eta\frown\langle
0,0\rangle}\,\frown\, 
\bar{a}_{\rho_\eta}
\,\frown\, \langle c_0,c_1\rangle,
\]
\[
\bar{b}_{\langle\rangle}\deq\langle c_0\rangle_{2n+2},
\]
where $\langle c\rangle_k$ stands for the sequence of $k$ entries, each of
which is
$c$, and $n=\llg(\bar{y})$ in $\varphi(\bar{x},\bar{y})$. We define
\begin{equation*}
\begin{split}
\psi(\bar{x},\bar{z})\equiv \psi(\bar{x},\bar{z}^0\,\frown\,\bar{z}^1\,
\frown\,w^0 \,\frown\,w^1)
&\equiv\\
& [w^0=w^1]\vee [
\varphi(\bar{x},\bar{z}^0) \wedge \neg \varphi(\bar{x},\bar{z}^1)],\\
\end{split}
\end{equation*}
where $\bar{z}=\bar{z}^0\frown\bar{z}^1\frown\langle w^0, w^1\rangle$
and $\llg(\bar{z}^0)=\llg(\bar{z}^1)=\llg(\bar{y})$.
We now claim that $\psi(\bar{x},\bar{z})$ and $\{\bar{b}_\eta:\,\eta\in
{}^{\omega>}2\}$
exemplify that SOP${}_1$ holds for $T$. Before we start checking this, note
that for
$\eta\in {}^{\omega>}2$ we have:
\begin{description}
\item{$\bullet_1$} $\psi(\bar{d}, \bar{b}_{\langle \rangle})$ holds for any
$\bar{d}$,

\item{$\bullet_2$} $\psi(\bar{d}, \bar{b}_{\eta\frown\langle 0 \rangle})$
holds iff $\varphi(\bar{d}, \bar{a}_{\rho_\eta\frown\langle 0,0\rangle})\wedge\neg
\varphi(\bar{d}, \bar{a}_{\rho_\eta})$ holds,

\item{$\bullet_3$} $\psi(\bar{d}, \bar{b}_{\eta\frown\langle 1 \rangle})$
holds iff $\neg \varphi(\bar{d}, \bar{a}_{\rho_\eta\frown\langle 1\rangle})\wedge
\varphi(\bar{d}, \bar{a}_{\rho_\eta})$ holds.
\end{description}
Let us verify \ref{nsop21}(2)(a), so let $\eta\in {}^\omega 2$. Pick
$\bar{c}$ such that $\varphi[\bar{c},\bar{a}_{\rho_\eta\rest
n}]^{\rho_\eta(n)}
$ holds for all $n<\omega$. We claim that
\begin{equation*}
\psi[\bar{c}, \bar{b}_{\eta\rest
n}]\mbox{ holds for all }n.\tag{$\ast$}
\end{equation*}
The proof is by a case analysis of $n$.

If \underline{$n=0$}, this is trivially true. If \underline{$n=k+1$} and
\underline{$\eta(k)=0$}, then we need to verify that 
$\varphi[\bar{c},\bar{a}_{\rho_{\eta\rest k}\frown\langle 0,0\rangle}]$ holds
and $\neg \varphi[\bar{c},\bar{a}_{\rho_{\eta\rest k}}]$ holds. We have
\[
\rho_{\eta\rest k}\frown\langle 0,0\rangle=\rho_\eta\rest (3k+2),
\]
and $\rho_\eta (3 k+2)=1$. Hence $\varphi[\bar{c},
\bar{a}_{\rho_{\eta\rest k}\frown\langle 0,0\rangle}]$ holds by the choice of
$\bar{c}$. On the other hand, we have $\rho_{\eta\rest k}=\rho_\eta\rest
(3k)$, and $\rho_\eta(3k)=\eta(k)=0$, hence
$\neg \varphi[\bar{c},\bar{a}_{\rho_{\eta\rest k}}]$ holds.

If \underline{$n=k+1$} and
\underline{$\eta(k)=1$}, then we need to verify that 
$\varphi[\bar{c},\bar{a}_{\rho_{\eta\rest k}}]$ holds
while $\varphi[\bar{c},\bar{a}_{\rho_\eta\rest (3k)\frown\langle 1\rangle}]$
does not. As $\rho_{\eta\rest k}=\rho_\eta\rest
(3k)$, and $\rho_\eta(3k)=\eta(k)=1$, we have that
$\varphi[\bar{c},\bar{a}_{\rho_\eta\rest k}]$ holds. Note that
$\varphi[\bar{c},\bar{a}_{\rho_\eta\rest (3k+2)}]$ holds as
$\rho_\eta (3k+2)=1$. We also have $(\rho_\eta\rest (3k+1))\frown
\langle 0\rangle\initialeq \rho_\eta\rest(3k+2)$. Hence
$\neg\varphi[\bar{c},\bar{a}_{\rho_\eta\rest(3k+1)}]$ by part
(b) in Definition \ref{another}. But
\[
\neg \varphi[\bar{c},\bar{a}_{\rho_\eta\rest (3k+1)}]\equiv
\neg \varphi[\bar{c},\bar{a}_{\rho_\eta\rest (3k)\frown\langle 1\rangle}]\equiv
\neg \varphi[\bar{c},\bar{a}_{\rho_{\eta\rest k}\frown\langle 1\rangle}]
\]
holds, so we are done proving $(\ast)$.

Let us now verify \ref{nsop21}(2)(b). So suppose $\nu\frown\langle 0\rangle
\initialeq\eta$ and consider $\{\psi(\bar{x},\bar{b}_{\nu\frown\langle
1\rangle}),\psi(\bar{x},\bar{b}_\eta)\}$. Let $\sigma$ and $l$ be such that $\eta=\sigma \frown
\langle l\rangle$.

\underline{Case 1}. $\nu=\sigma$.

Hence $l=0$. So $\psi(\bar{x},\bar{b}_\eta)\implies \neg 
\varphi(\bar{x},\bar{a}_{\rho_\nu})$
and
$\psi(\bar{x},\bar{b}_{\nu\frown\langle 1\rangle})
\implies \varphi(\bar{x},\bar{a}_{\rho_\nu})$, by $\bullet_2$
and $\bullet_3$ respectively,
showing that $\{\psi(\bar{x},\bar{b}_\eta),
\psi(\bar{x},\bar{b}_{\nu\frown\langle 1\rangle})\}$ is inconsistent.

\underline{Case 2}. $\nu\initial\sigma$ and $l=0$.

Hence $\nu\frown\langle 0\rangle\initialeq\sigma$. Clearly $\rho_\nu\frown
\langle 0\rangle \initialeq \rho_\sigma\frown\langle 0,0\rangle$,
as 
\[
\rho_\sigma(\llg(\rho_\nu))=\sigma(\llg(\nu))=0.
\]
We have $\psi(\bar{x},\bar{b}_{\nu\frown\langle 1\rangle})
\implies \varphi(\bar{x},\bar{a}_{\rho_\nu})$ by $\bullet_3$ and
$\psi(\bar{x},\bar{b}_\eta)=\psi(\bar{x},\bar{b}_\sigma\frown\langle 0\rangle)$ implies
$\varphi(\bar{x},\bar{a}_{\rho_\sigma\frown\langle 0,0\rangle})$ by
$\bullet_2$, while the
two formulae being implied are contradictory, by (b) in the definition of
SOP${}'_{1}$.

\underline{Case 3}. $\nu\initial\sigma$ and $l=1$.

Observe that $\psi(\bar{x},\bar{b}_\eta)\implies
\varphi(\bar{x},\bar{a}_{\rho_\sigma})$ by $\bullet_3$ and
$\psi(\bar{x},\bar{b}_{\nu\frown\langle 1\rangle})\implies
\varphi(\bar{x},\bar{a}_{\rho_\nu})$. As above, using $\nu\frown\langle 0
\rangle \initialeq\sigma$, we
show that the set $\{\varphi
(\bar{x},\bar{a}_{\rho_\nu}), \varphi
(\bar{x},\bar{a}_{\rho_\sigma})\}$ is inconsistent.
$\eop_{\ref{equivalent}}$
\end{Proof of the Claim}

\begin{Conclusion}\label{zakljucnica} $T$ has SOP${}_1$ iff $T$ has property
SOP${}'_1$ 
from Claim \ref{equivalent}.
\end{Conclusion}

\begin{Question} Is the conclusion of \ref{zakljucnica} true when
the theory $T$
is replaced by a formula $\varphi$?
\end{Question}

{\em Start changes}

It turns out that witnesses to being SOP${}_1$ can be chosen to be 
highly indiscernible.

\begin{Definition}\label{treeindisc}
(1) Given an ordinal $\alpha$ and sequences $\bar{\eta}_l=
\langle \eta^l_0,\eta^l_1,\ldots,\eta_{n_l}^l\rangle$ for $l=0,1$
of members of ${}^{\alpha>} 2$, we say that $\bar{\eta}_0\iso_1\bar{\eta}_1$
iff
\begin{description}
\item{(a)} $n_0=n_1$,
\item{(b)} the truth values of
\begin{itemize}
\item $\eta^l_k=\langle\rangle$,
\item $\eta^l_{k_1}\cap \eta^l_{k_2}\initialeq \eta^l_{k_3}\cap\eta^l_{k_4}$
\end{itemize}
do not depend on $l$,
\item{(c)} 
$\eta^l_{k_1}\ntrianglelefteq\eta^l_{k_2}\implies \eta^0_{k_1} (\llg(\eta^0_{k_1}\cap \eta^0_{k_2}))=
\eta^1_{k_1} (\llg(\eta^1_{k_1}\cap \eta^1_{k_2}))$
\end{description}
for $k_1,k_2, k_3, k_4\le n_0$.

{\noindent (2)}
We say that the sequence
$\langle\bar{a}_\eta:\, \eta\in {}^{\alpha>}2\rangle $ of elements of
$\mathfrak C$ (for an
ordinal $\alpha$) is {\em 1-fully binary tree indiscernible (1-fbti)}
iff whenever $\bar{\eta}_0\iso_1\bar{\eta}_1$ are sequences of elements
of ${}^{\alpha>}2$, \underline{then}
\[
\bar{a}_{\bar{\eta}_0}\deq\bar{a}_{\eta^0_0}\frown\ldots \frown\bar{a}_{\eta^0_{n_0}}
\]
and the similarly defined $\bar{a}_{\bar{\eta}_1}$,
realise the same type in $\mathfrak C$.

{\noindent (3)} Suppose that $\delta$ is a limit ordinal $>0$. Define
$h^\ast=h^\ast_\delta:\,{}^{\delta>}2\to {}^{\delta>}2$ by letting for $\eta\in {}^{\delta>}2$
\begin{itemize}
\item $\llg(h^\ast(\eta))=2\llg(\eta)+1$,
\item $i<\llg(h^\ast(\eta))\implies h^\ast(\eta)(2i)=0, h^\ast(\eta)(2i+1)=\eta(i)$,
\item $h^\ast(\eta)(2\llg(\eta))=1$.
\end{itemize}
For $n<\omega$ and $\bar{\eta}\in {}^n({}^{\delta>}2)$ we define $h^\ast(\bar{\eta})=
\langle h^\ast (\eta_l):\,l<n\rangle$.

We say $\bar{\eta}\iso_2\bar{\nu}$ iff $h^\ast (\bar{\eta})\iso_1 h^\ast(\bar{\nu})$. We define 2-fbti
like 1-fbti but using $\iso_2$ in place of $\iso_1$.

\end{Definition} 

{\bf Observation 2.10 A} The following can be easily checked:

{\noindent (1)} Let $\bar{\eta}, \bar{\nu}\in {}^n({}^{\alpha>}2)$ and let $\bar{\eta}'$ and $\bar{\nu}'$
be the closures of $\bar{\eta}, \bar{\nu}$, respectively, under intersections. \underline{Then} $\bar{\eta}\iso_1 \bar{\nu}$ iff $\bar{\eta}'\iso_1 \bar{\nu}'$.

{\noindent (2)} If $\langle \bar{a}_{\bar{\eta}}:\,\eta\in {}^{\delta>} 2\rangle $ is 1-fbti then 
$\langle \bar{a}_{h^\ast(\bar{\eta})}:\,\eta\in {}^{\delta>} 2\rangle $ is 2-fbti.

{\noindent (3)} $h^\ast(\eta)$ is never $\langle,\rangle$ and $h^\ast(\eta_0)$ is never a strict initial segment of
$h^\ast(\eta_1)$.

\begin{Claim}\label{thinning} If 
$t\in \{1,2\}$ and $\langle \bar{b}_\eta:\,\eta\in {}^{\omega>}
2 \rangle$ are of given constant length, and $\delta\ge \omega$ is a (limit for $t=2$) ordinal, \underline{then} we can find
$\langle \bar{a}_\eta:\,\eta\in {}^{\delta>}2\rangle$ such that
\begin{description}
\item{(a)} $\langle \bar{a}_\eta:\,\eta\in {}^{\delta>}2\rangle$ is $t$-fbti,
\item{(b)} if $\bar{\eta}=\langle \eta_m:\,m<n\rangle$, where each $\eta_m
\in {}^{\delta>}2$, is given, \underline{then}
we can find $\nu_m\in {}^{\omega>}2\,(m<n)$ such that with $\bar{\nu}\deq
\langle \nu_m:\,m<n\rangle$, we have $\bar{\nu}\iso_t\bar{\eta}$
and sequences $\bar{a}_{\bar{\eta}}$
and $\bar{b}_{\bar{\nu}}$ realise the same type in $\mathfrak C$. 
\end{description}
\end{Claim}

\begin{Proof of the Claim}\footnote{Note that the definition of $\iso_1,\iso_2$ has changed from the one given in the published version of this paper, but the following proof is basically the same as the one there.}
Let us first deal with $t=1$. By Observation 2.10 A (1) above, we may reduce to checking clause (b) only for tuples
$\bar{b}_\eta$ where $\eta$ is closed under intersections.
By Compactness Theorem it suffices
to assume $\delta=\omega$.
The proof goes through a series of steps through which we obtain 
increasing degrees of indiscernibility. We shall need some auxiliary
definitions. Let $\alpha$ be an infinite ordinal.

\begin{Definition}\label{semiindis} (1) Given $\bar{\eta}=\langle
\eta_0,\ldots,\eta_{k-1}\rangle$, a sequence of elements of ${}^{\alpha>}
2$, and an ordinal $\gamma$. We define $\bar{\eta}'={\rm cl}_\gamma
(\bar{\eta})$ as follows:
\[
\bar{\eta}'=\langle\langle\rangle, \eta_0,\eta_0\rest\gamma,\eta_1,\eta_1\rest\gamma,\eta_0
\cap\eta_1,\eta_2,\eta_2\rest\gamma,
\eta_0\cap\eta_2,\eta_1\cap\eta_2\ldots\rangle. \]
We also define $u_\gamma[\bar{\eta}]=\{\eta_l\in \bar{\eta}:\,\llg(\eta_l)>\gamma\}$.

{\noindent (2)} We say that $\bar{\eta}\iso_{\gamma,n}\bar{\nu}$ iff
$\bar{\eta}'\deq{\rm cl}_\gamma
(\bar{\eta})$ and $\bar{\nu}'\deq{\rm cl}_\gamma
(\bar{\nu})$ satisfy
\begin{description}
\item{(i)} $\bar{\eta}'=\langle \eta'_l:\,l<m\rangle$ and
$\bar{\nu}'=\langle \nu'_l:\,l<m\rangle$ are both in ${}^m({}^{\alpha>}2)$
for some $m$,
\item{(ii)} for $l< m$ we have $\eta'_l\in {}^{\gamma\ge} 2
\iff \nu'_l\in {}^{\gamma\ge} 2$, and for such $l$ we have $\eta'_l
=\nu'_l$,
\item{(iii)} $n\ge\card{u_\gamma[\bar{\eta}]}$, 

\item{(iv)} $\eta'_l, \eta'_k\in u_\gamma[\bar{\eta}]\implies [\llg(\eta'_l)<\llg(\eta'_k)
\iff \llg(\nu'_l)<\llg(\nu'_k)$],

\item{(v)} $\eta'_{l_1}\initialeq \eta'_{l_2}\iff\nu'_{l_1}
\initialeq \nu'_{l_2}$,
and the same holds for the equality,
\item{(vi)} if $\eta'_{l_1}$ is not an initial segment of $\eta'_{l_2}$, then
$\eta'_{l_1}(\llg(\eta'_{l_1}\cap \eta'_{l_2}))=\nu'_{l_1}(\llg(\nu'_{l_1}\cap \nu'_{l_2}))$.
\end{description}

{\noindent (3)} $\langle \bar{a}_\eta:\,\eta\in {}^{\alpha> }2\rangle$
is $(\gamma, n)$-{\em indiscernible} iff
for every $k$, for every $\bar{\eta},\bar{\nu}
\in {}^k({}^{\alpha>} 2)$ with $\bar{\eta}\iso_{\gamma,n}\bar{\nu}$,
the tuples $\bar{a}_{\bar{\eta}}$ and $\bar{a}_{\bar{\nu}}$ realise the
same type.

{\noindent (4)} $(\le \gamma, n)$-{\em indiscernibility} is the
conjunction of $( \beta, n)$-indiscernibility for all $\beta\le
\gamma$.

{\noindent (5)} We say that $\langle \bar{a}_\eta:\,\eta\in
 {}^{\alpha> }2\rangle$ is 0-{\em fbti} iff it is $(\gamma, n)$-indiscernible
 for all $\gamma$ and $n$.

 \end{Definition}
 
 {\bf Note 2.12 A} (1) ${\rm cl}_0(\bar{\eta})$ is simply the closure of $\bar{\eta}$ under intersections,
 joined with $\langle\rangle$ in appropriate places.
 
 {\noindent (2)} $\bar{\eta}\iso_{\gamma ,n}\bar{\nu}$ iff
${\rm cl}_\gamma(\bar{\eta}) \iso_{\gamma ,n}
{\rm cl}_\gamma
(\bar{\nu})$.
 
\begin{Subclaim}\label{pomocnis} If $\bar{a}_\eta\in {}^k{\mathfrak C}$ for
$\eta\in {}^{\omega>}2$ are tuples of constant length and closed under intersections, \underline{then}

for any $\alpha\ge\omega$ we can find $\bar{a}'=
\langle \bar{a}'_\eta:\,\eta\in
 {}^{\alpha> }2\rangle$ such that
 \begin{description}
 \item{($x$)} $\bar{a}'$ is 0-fbti,
 \item{($xx$)} for every $m$ and a finite set $\Delta$ of formulae, we can find
 $h:\,{}^{m\ge}2\into {}^{\omega>}2$ such that
 \begin{description}
 \item{$(\alpha)$} $\langle \bar{a}'_\eta:\,\eta\in
 {}^{m\ge  }2\rangle$ and $\langle \bar{a}_{h(\eta)}:\,\eta\in
 {}^{m\ge }2\rangle$ realise the same $\Delta$-type,
 \item{$(\beta)$} $h$ satisfies $h(\eta)\,\hat{}\,\langle l\rangle
 \initialeq h(\eta\,\hat{}\,\langle l\rangle)$ for $\eta \in
 {}^{m>} 2$ and $l<2$, and
 \[
 \llg(\eta_1)=\llg(\eta_2)\implies \llg(h(\eta_1))=
 \llg(h(\eta_2)).
 \]
 \end{description}
 \end{description}
 \end{Subclaim}
 
 \begin{Proof of the Subclaim} By Compactness Theorem it suffices
to work with $\alpha=\omega$.

Let $(\ast)_{\gamma,n}$ be the conjunction of
the statement $(x)_{\gamma,n}$ given by
 \begin{equation*}
 \bar{a}'\mbox{ is }(\le\gamma,n)\mbox{-indiscernible},
 \end{equation*}
 and $(xx)$ above. We prove by induction on $n$ and then $\gamma$ that
 for any $\gamma\le\omega$ we can find $\bar{a}'$ for which $(\ast)_{\gamma,n}$
 holds.  
 
 \underline{$n=0$}. We use $\bar{a}'_\eta=\bar{a}_\eta$.
 
 \underline{$n+1$}. By
 induction on $\gamma\le\omega$, we prove that there is $\bar{a}'$ for which
 $(\ast)_{\gamma,n+1}+(\ast)_{\omega,n}+(xx)$ holds.
  \begin{description}
\item{\underline{$\gamma<\omega$}}. 
\end{description}
Without loss of
generality, the sequence $\langle \bar{a}_\eta:\,\eta\in {}^{\omega>}2\rangle$
is $(\le\omega,n)$-indiscernible, as $(xx)$ as a relation between
$\langle \bar{a}_\eta:\,\eta
\in {}^{\omega>}2\rangle$ and $\langle \bar{a}'_\eta:\,\eta\in
{}^{\omega>}2\rangle$ is transitive.
Suppose we are given $\bar{\eta}^\ast,\bar{\nu}^\ast$ 
satisfying $\bar{\eta}^\ast\iso_{\gamma,n+1}\bar{\nu}^\ast$. By Note 2.1.2 A, we may assume
$\bar{\eta}^\ast,\bar{\nu}^\ast$ to be the same as their cl${}_\gamma$ closures
and the same will hold for any $\bar{\eta}, \bar{\nu}$ that we
mention in this context.

If $\card{u_\gamma[\bar{\eta^\ast}]}\le n$, the conclusion follows by the assumptions.
We shall assume $\card{u_\gamma[\bar{\eta^\ast}]}>n$, so
$\card{u_\gamma[\bar{\eta^\ast}]}=n+1$. Moreover, if
$\min(u_\gamma[\bar{\eta^\ast}])=\min(u_\gamma[\bar{\nu^\ast}])$ and for any $l$ with
$\llg(\eta^\ast_l)=\min(u_\gamma[\bar{\eta^\ast}])$ we have $\eta^\ast_l=\nu^\ast_l$, then using
$(x)_{\min(u_\gamma[\bar{\eta^\ast}]),n}$, we get that $\bar{a}_{\bar{\eta}^\ast}$
and $\bar{a}_{\bar{\nu}^\ast}$ realise the same type.
By the same argument,
fixing a finite set $\Delta$ of formulae,
for every $\bar{\eta}$, we get that the
${\rm{tp}}_\Delta(\bar{a}_{\bar{\eta}}) $ depends just on the

\[
\bar{\eta}/\iso_{\gamma,n+1}\deq\Upsilon\mbox{ and }
\{\eta_l:\,l<\llg(\bar{\eta})\}\cap {}^{\min(u_\gamma[\bar{\eta}])}2=\{\eta_l:\,l\in v^\Upsilon\}
\]
for some
$v^\Upsilon\subseteq\llg(\bar{\eta})$. Let us define $F_{\Upsilon,\Delta}^0$
by $F_{\Upsilon,\Delta}^0(\langle\eta_l:\,l\in v^\Upsilon\rangle)={\rm
tp}_\Delta(\bar{a}_{\bar{\eta}})$. By the closure properties of $\bar{\eta}$
and the definition of $\iso_{\gamma, n+1}$, we get
that for $l_1\neq l_2\in v^\Upsilon$ the truth value of
$\eta_{l_1}\rest(\gamma+1)=
 \eta_{l_2}\rest(\gamma+1)$ depends only on $\Upsilon$. We can hence replace
 $v^\Upsilon$  by  a  set  $v_\ast^\Upsilon\subseteq v^\Upsilon$ such that
 $\langle\eta_l:\,l\in  v_\ast^\Upsilon\rangle$ are the  representatives under
the equality
 of the restrictions to $\gamma+1$.
 
 As we have fixed $\Delta$, there is a finite set $A$ of $\Upsilon$s
that can be used as representatives for the values of $F^0_{\Upsilon,
\Delta}$. Let $r$ be the size of the range of $F^0_{\Upsilon,
\Delta}$.  Let
 $k^\ast=2^{\gamma+1}$  (so  finite)  and 
 let $\{\mu^\ast_k:\,k<k^\ast\}$ list ${}^{\gamma+1}2$. We define
a function $F_{\Upsilon,\Delta}$ on ${}^{k^\ast}({}^{\omega>}2)$ by letting
 \begin{equation*}
 \begin{split}
 F_{\Upsilon,\Delta}(x_0,\ldots,x_k,\ldots)_{k<k^\ast}\deq&
 F_{\Upsilon,\Delta}^0(\langle\eta_l:\,l\in v^\Upsilon_\ast\rangle),\\
 & \mbox{ 
 where }
 \eta_l\rest(\gamma+1)=\mu^\ast_k\implies\eta_l=\mu^\ast_k\frown x_k.\\
 \end{split}
 \end{equation*}
 Define a function $F$ with arity $k^\ast$ so that
 $F((\ldots,x_k,\ldots)_{k<k^\ast})$  is  defined  iff for some $m<\omega$ we
 have
 $\{x_k:\,k<k^\ast\}\subseteq {}^m 2$ and then
 \[
 F((\ldots,x_k,\ldots)_{k<k^\ast})=\langle
 F_{\Upsilon,\Delta}((\ldots,x_k,\ldots)_{k<k^\ast}):\,\Upsilon\in A
 \rangle.
 \]
 Therefore $F$ is a function from $\bigcup_{m<\omega} \prod_{k<k^\ast}  {\rm lev}_m ({}^{\omega>}2)$ into a set of
 size $r$.
We recall the following definition and restatement of the Halpern-Lauchli theorem \cite{HaLa},
due to Laver and Pincus and presented in \cite{PiHa}. 

{\bf Definition 2.13 A} (1) A tree $S$ is {\em strongly embedded} in a tree $T$ if  there is a
strictly increasing embedding $f^\ast$ of $S$ as a suborder of $T$
such that
\begin{itemize}
\item any nonmaximal node in $f(S)$ has the same number of immediate successors in $T$ and in $f(S)$,
and
\item all nodes on any common level of $S$ are mapped by $f$ to a common level of $T$.

{\noindent (2)} A nonempty subtree of  ${}^{\omega>}\omega$ is {\em well-behaved} if it is finitely branching and has no maximal nodes (hence it has $\omega$ levels).
\end{itemize}

{\bf Halpern-Lauchli theorem} 
Let $r,d<\omega$. Suppose that $\langle T_i:\,i<d\rangle$ are well-behaved trees and that
$c$ is a colouring of $\bigcup_{n<\omega}\prod_{i<d} 
{\rm lev}_{n}(T_i)$ into $r$ colours. 
\underline{Then} there are $f^\ast$,  $\langle S_i:\,i<d\rangle$ and $\langle h_i:\,i<d\rangle$
such that
\begin{itemize}
\item $f^\ast:\,\omega\into\omega$ is a strictly increasing function,
\item each $S_i$ is a well-behaved tree,
\item $h_i$ is
a strong embedding of $S_i\into T_i$, 
\item for each $n<\omega$ and $i<d$,
the common height in $T_i$ of elements of $h_i``{\rm lev}_n(S_i)$ is $f^\ast(n)$, and
\item
$\bigcup_{n\in u}\prod_{i<d} 
h_i``{\rm lev}_{n}(S_i)$ is $c$-monochromatic.
\end{itemize}
Moreover, in the case that all $T_i$ are the same tree, we can assume that all $h_i$ are contained in a common function $h$.

Therefore we can apply the Halpern-Lauchli theorem to $F$. We get a sequence $\langle S_k:\,k<k^\ast\rangle$ of well-behaved trees exemplify the conclusion of the Halpern-Lauchli theorem with $h_k=h\rest S_k$ and $f^\ast(n)={\rm ht}[h``{\rm lev}_n(S_k)]$. Since the only well-behaved subtree of ${}^{\omega>}2$
is ${}^{\omega>}2$ itself, we 
can conclude that there is $h:\,{}^{\omega>}2\into {}^{\omega>}2$ such that
 \begin{itemize}
 \item $h\rest {}^{\gamma\ge} 2$ is the identity,
 \item $\llg(h(\eta))$ depends just on $\llg(\eta)$ (not on
 $\eta$),
 \item $h(\eta)\frown\langle l\rangle\initial h(\eta\frown\langle
 l\rangle)$
 for $l=0,1$,
 \item for some $c$ we have that for all $m<\omega$
 \[
 \{\eta_k:\,k<k^\ast\}\subseteq{}^m2\implies
 F((h(\eta_0),h(\eta_1),\ldots,h(\eta_k),\ldots  )_{k<k^\ast})=c. 
 \]
 \end{itemize}
 Let $\bar{a}'_\eta$ for $\eta\in {}^{\omega>}2$ be defined to be
 $\bar{a}_\eta$ if $\eta\in {}^{\gamma>} 2$, and otherwise $\bar{a}_{h(\nu)}$ for the unique $\nu$ such that $\eta\rest(\gamma+1)=\mu^\ast_k$ and
$\eta=\mu^\ast_k\frown\nu$.

 We have obtained the desired conclusion, but localized to $\Delta$.
 The induction
 step ends  by an application of the compactness theorem.
 
 \begin{description}
\item{\underline{$\gamma=\omega$}} This is vacuously true.
\end{description}
Having carried the induction, the
conclusion of the Subclaim follows from $\bigwedge_n (\ast)_{0,n}$. 
$\eop_{\ref{pomocnis}}$
\end{Proof
of the Subclaim}   
Now we go {\em back to the proof of the Claim}. Given $\langle\bar{b}_\eta:\,\eta\in {}^{\omega>}2\rangle$ as in the
assumptions, by the Subclaim we can assume that they are 0-fbti.
We choose by induction on $n$ a function $h_n:\,{}^{n\ge }2\into
{}^{\omega>}2$ as follows.
Let $h_0(\langle\rangle)=\langle\rangle$. If $h_n$ is defined, let
\[
k_n\deq\max\{\llg(h_n(\eta))+1:\,\eta\in {}^{n\ge}2\}
\]
and let
\[
h_{n+1}(\langle\rangle)=\langle\rangle,\quad
h_{n+1}(\langle 1\rangle\,\hat{}\,\nu)=\langle 1\rangle\,\hat{}\,h_n(\nu),\quad
h_{n+1}(\langle 0\rangle\,\hat{}\,\nu)=\langle 0,\ldots,0\rangle\hat{}h_n(\nu),
\]
where the sequence of 0s in the last part of the definition has length $k_n$.
The point of the definition of $h_n$ is that if $\bar{\eta}^l=\langle
\eta^l_0,\ldots,\eta^l_{n_l}\rangle$ for $l=0,1$ are given and
$n^\ast=\llg(\cl_0(\bar{\eta}^0))$, then
\[
\bar{\eta}^0\iso_1\bar{\eta}^1\implies\langle h_{n^\ast}(\eta^0_0),\ldots,
h_{n^\ast}(\eta^0_{n_0})\rangle\iso_{0,n^\ast}
\langle h_{n^\ast}(\eta^1_0),\ldots,
h_{n^\ast}(\eta^1_{n_1})\rangle.
\]
To check this, we verify the six relevant items of the definition of $\iso_{0,
n^\ast}$.

\begin{description}

\item{(i)} Follows because $n_0=n_1$ by the definition of $\iso_1$.

\item{(ii)} If $h_{n^\ast}(\eta^0_i)\cap h_{n^\ast}(\eta^0_j)=\langle\rangle$
then $\eta^0_i\cap \eta^0_j=\langle\rangle$ so $\eta^1_i\cap \eta^1_j
=\langle\rangle$ by the definition of $\iso_1$, and hence
$h_{n^\ast}(\eta^1_i)\cap h_{n^\ast}(\eta^1_j)=\langle\rangle$. The opposite implication
holds by symmetry.

\item{(iii)} Follows by the definition of $n^\ast$.

\item{(iv)} Suppose 
\[
0<\llg(h_{n^\ast}(\eta^0_i)\cap h_{n^\ast}(\eta^0_j))<\llg(h_{n^\ast}(\eta^0_k)\cap h_{n^\ast}(\eta^0_s)).
\]
Let $m\le n^\ast$ be the first such that
\[
0<\llg(h_{n^\ast}(\eta^0_i\rest m)\cap h_{n^\ast}(\eta^0_j\rest m))<\llg(h_{n^\ast}(\eta^0_k
\rest m)\cap h_{n^\ast}(\eta^0_s\rest m)).
\]
Clearly, $m>0$. To simplify the notation, let us assume that $m=n^\ast$.
Let $\eta^0_t=\langle l_t\rangle\frown\nu^0_t$ for $t\in \{i,j,k,s\}$
and for some $l_t\in \{0,1\}$ depending on $t$. The situation we
describe can happen iff $l_i=l_j=1$ and $l_k=l_s=0$, by the definition
of $h_n$. By the definition of $\iso_1$ this can be recognised by the
$\iso_1$-type of $\bar{\eta}^0$.

\item{(v), (vi)} Follow because the corresponding properties are
preserved by $h_n$.

\end{description}

Fix an $n<\omega$ and define $\bar{a}_\eta=\bar{b}_{h_n(\eta)}$ for
$\eta\in {}^{n\ge 2}$. By the above argument it follows
that $\langle \bar{a}_{\eta}:\,\eta\in {}^{n\ge}2\rangle$ are
1-fbti. As $n$ was arbitrary, we can finish by compactness.

For $t=2$, we use exactly the same proof.
$\eop_{\ref{thinning}}$
\end{Proof of the Claim}

The following Theorem \ref{3indiswitness} will finally tell us that witnesses for SOP${}_1$ can be chosen with a certain degree of
indiscernability. We need to introduce a new notion of indiscernability:

\begin{Definition}\label{3treeindisc}
(1) Given an ordinal $\alpha$ and sequences $\bar{\eta}_l=
\langle \eta^l_0,\eta^l_1,\ldots,\eta_{n_l}^l\rangle$ for $l=0,1$
of members of ${}^{\alpha>} 2$, we say that $\bar{\eta}_0\iso_3\bar{\eta}_1$
iff
\begin{description}
\item{(a)} $n_0=n_1$,
\item{(b)} the truth values of
\begin{itemize}
\item $\eta^l_k=\langle\rangle$,
\item $\eta^l_{k_1}\cap \eta^l_{k_2}\initialeq \eta^l_{k_3}\cap\eta^l_{k_4}$
\end{itemize}
do not depend on $l$,
\item{(c)} 
$\eta^l_{k_1}\ntrianglelefteq\eta^l_{k_2}\implies \eta^0_{k_1} (\llg(\eta^0_{k_1}\cap \eta^0_{k_2}))=
\eta^1_{k_1} (\llg(\eta^1_{k_1}\cap \eta^1_{k_2}))$,
\item{(d)}  $\eta^l_{k_1}\ntrianglelefteq\eta^l_{k_2}\implies \eta^0_{k_1}=(\eta^0_{k_1}\cap \eta^0_{k_2})
\frown\langle 1\rangle$ iff $\eta^1_{k_1}=(\eta^1_{k_1}\cap \eta^1_{k_2})\frown\langle 1\rangle$.
\end{description}
for $k_1,k_2, k_3, k_4\le n_0$.

{\noindent (2)}
We say that the sequence
$\langle\bar{a}_\eta:\, \eta\in {}^{\alpha>}2\rangle $ of elements of
$\mathfrak C$ (for an
ordinal $\alpha$) is {\em 3-fully binary tree indiscernible (3-fbti)}
iff whenever $\bar{\eta}_0\iso_1\bar{\eta}_1$ are sequences of elements
of ${}^{\alpha>}2$, \underline{then}
\[
\bar{a}_{\bar{\eta}_0}\deq\bar{a}_{\eta^0_0}\frown\ldots \frown\bar{a}_{\eta^0_{n_0}}
\]
and the similarly defined $\bar{a}_{\bar{\eta}_1}$,
realise the same type in $\mathfrak C$.
\end{Definition} 

\begin{theorem}\label{3indiswitness} Suppose that $T$ has SOP${}_1$ as witnessed by
$\varphi(\bar{x}, \bar{y})$ and a sequence $\bar{a}=\langle \bar{a}_\eta:\,\eta\in {}^{\omega>}2\rangle$. 
\underline{Then} there is $\bar{d}=\langle \bar{d}_\eta:\,\eta\in {}^{\omega>}2\rangle$ exemplifying that
$\varphi(\bar{x}, \bar{y})$ has SOP${}_1$ and $\langle \bar{d}_\eta:\,\eta\in {}^{\omega>}2\setminus
\{\langle\rangle\}\rangle$ is 3-fbti.
\end{theorem}

\begin{proof} Let $k^\ast=\llg(\bar{y})$. First define $\bar{b}_\eta$ for $\eta\in {}^{\omega>}2$ by
$\bar{b}_\eta=\bar{a}_{\eta\frown\langle 0\rangle}\frown \bar{a}_{\eta\frown\langle 1\rangle}$.
Let for any $\bar{z_0}, \bar{z_1}$ of length $k^\ast$ and $l\in \{0,1\}$,
$\psi_l(\bar{x}, \bar{z}_0\frown \bar{z}_1)\equiv \varphi(\bar{x}, \bar{z}_l)$. Now we use
Claim \ref{thinning} applied to $\langle \bar{b}_\eta:\,\eta\in {}^{\omega>}2\rangle$. Therefore we can
find $\bar{c}=\langle \bar{c}_\eta:\,\eta\in {}^{\omega>}2\rangle$ such that
\begin{description}
\item{(a)} $\bar{c}$ is 1-fbti,
\item{(b)} for any finite $n$ and $\bar{\eta}\in {}^n({}^{\omega>}2)$ there is 
$\bar{\nu}\in {}^n({}^{\omega>}2)$ such that $\bar{\nu} \iso_1\bar{\eta}$ and $\bar{b}_{\eta}$ and $\bar{c}_{\nu}$ realise the same type in $\mathfrak C$.
\end{description}
Let $\bar{d}_{\eta}$ for $\eta\in {}^{k^\ast} ({}^{\omega>}2)$ be defined by induction on the length of
$\eta$ so that $\bar{d}_{\eta\frown\langle 0\rangle}\frown \bar{d}_{\eta\frown\langle 1\rangle}=
\bar{c}_\eta$ and $\bar{d}_{\langle\rangle}=\bar{c}_{\langle\rangle}$. This is possible by the choice of 
$\bar{b}$ and $\bar{c}$.

\begin{claim}\label{prvisubclaim} If $\nu\frown\langle 0\rangle\initialeq\eta$ then 
$\varphi(\bar{x}, \bar{d}_\eta)$ and $\varphi(\bar{x}, \bar{d}_{\nu\frown\langle 1\rangle})$ are incompatible.
\end{claim}

\begin{Proof of the Claim} Let $\eta=\rho\frown\langle l\rangle$ for some $l\in \{0,1\}$. Consider
$\{\psi_l(\bar{x},\bar{c}_\rho), \psi_1(\bar{x}, \bar{c}_\nu)\}$, we claim that this set is inconsistent.
We know that 
\[
\psi_l(\bar{x},\bar{c}_\rho)\equiv\varphi(\bar{x}, \bar{d}_{\rho\frown \langle l\rangle})\equiv 
\varphi(\bar{x}, \bar{d}_\eta), \quad \psi_1(\bar{x}, \bar{c}_\nu)\equiv \varphi(\bar{x}, \bar{d}_{\nu\frown
\langle 1\rangle}).
\]
By the 1-fbti property of $\bar{c}$ and the choice of $\bar{c}$ with respect to $\bar{b}$ it suffices to check that $\{\psi_l(\bar{x},\bar{b}_\rho), \psi_1(\bar{x}, \bar{b}_\nu)\}$ is inconsistent. This means that 
$\{\varphi(\bar{x},\bar{a}_\eta), \varphi(\bar{x}, \bar{a}_{\nu\frown \langle 1\rangle}\}$ is inconsistent, which is true by the choice of $\bar{a}$.
$\eop_{\ref{prvisubclaim}}$
\end{Proof of the Claim}

\begin{claim}\label{drugiclaim} For any $\rho\in {}^\omega 2$, $\{\varphi(\bar{x}, \bar{d_{\rho\rest n}}):\,n<\omega\}$ is consistent.
\end{claim}

\begin{Proof of the Claim} It suffices to show that for any
\[
\langle\rangle\initial \eta_0\initial \eta_1\initial\ldots \eta_k
\]
the set $\{\varphi(\bar{x}, \bar{d}_{\eta_{l+1}\rest\llg(\eta_l)\frown \eta_{\l+1}(\llg(\eta_l)}):\,l>k\}\cup
\{\varphi(\bar{x},\bar{d}_{\langle\rangle})\}$ is consistent. This means $\{\psi_{\eta_{\l+1}(\llg(\eta_l))}
(\bar{x}, \bar{c}_{\eta_{l+1}\rest\llg(\eta_l)}):\,l<k\}\cup \{\varphi (\bar{x}, \bar{c}_{\langle\rangle})\}$
is consistent. By the choice of $\bar{b}$ and $\bar{c}$ this is to say
$\{\psi_{\eta_{\l+1}(\llg(\eta_l))}
(\bar{x}, \bar{b}_{\eta_{l+1}\rest\llg(\eta_l)}):\,l<k\}\cup \{\varphi(\bar{x}, \bar{a}_{\langle\rangle})\}$ or
$\{\varphi(\bar{x}, \bar{a}_{\eta_{l+1}\rest\llg(\eta_l)}):\,l<k\}\cup \{\varphi(\bar{x}, \bar{a}_{\langle\rangle})\}$
is consistent, but this is true by the choice of $\bar{a}$.
$\eop_{\ref{drugiclaim}}$
\end{Proof of the Claim}

\begin{claim}\label{treciclaim} $\langle d_{\eta}:\,\eta\in {}^\omega 2\setminus \{0\}\rangle$ is 3-fbti.
\end{claim}

\begin{Proof of the Claim} Suppose that $\bar{\eta}_0\iso_3\bar{\eta}_1$ and consider
$\bar{d}_{\bar{\eta_0}}$ and $\bar{d}_{\bar{\eta_1}}$. For each $\eta_l^k$ let $\nu^l_k$ be such that
$\eta^l_k= \nu^l_k\frown \langle m^l_k\rangle$ for some $m^l_k\in \{0,1\}$ and let 
$\bar{\nu}_0,\bar{\nu}_1$ be defined from $\nu^l_k (l\in \{0,1\}, k<\llg(\bar{\eta}_0))$. Then
$\bar{\eta}_0\iso_3\bar{\eta}_1\implies \bar{\nu}_0\iso_3\bar{\nu}_1$, hence $\bar{c}_{ \bar{\nu}_0}$
and $\bar{c}_{ \bar{\nu}_1}$ realise the same type, which implies that $\bar{d}_{ \bar{\eta}_0}$
and $\bar{d}_{ \bar{\eta}_1}$ do.
$\eop_{\ref{treciclaim}}$
\end{Proof of the Claim}

$\eop_{\ref{3indiswitness}}$
\end{proof}

{\em End changes.}

As we mentioned before, it would be really interesting to know if
SOP${}_2$ and SOP${}_1$ are equivalent. A step towards understanding this
question is provided by
the next claim which shows that
in the case of theories which are SOP${}_1$ and NSOP${}_2$, the witnesses to
being SOP${}_1$
can be chosen to be particularly nice. {\em note a change here to 3-fbti from the old version}

\begin{Claim}\label{difference} Suppose that
$\varphi(\bar{x},\bar{y})$ satisfies SOP${}_1$,
but for no $n$ does the formula
$\varphi_n(\bar{x},\bar{y}_0,\ldots,\bar{y}_{n-1})\equiv
\wedge_{k<n}\varphi(\bar{x},\bar{y}_k)$ satisfy
SOP${}_2$. \underline{Then} there are witnesses
$\langle \bar{a}_\eta:\,\eta\in {}^{\omega>}2\rangle$ for
$\varphi(\bar{x},\bar{y})$
satisfying SOP${}_1$
which in addition satisfy:
\begin{description}
\item{(c)} if $X\subseteq {}^{\omega>}2$, and there are no $\eta,\nu\in
X$ such that $\eta\frown\langle 0\rangle\initialeq \nu$, \underline{then}
$\{\varphi(\bar{x},\bar{a}_\eta):\,\eta\in X\}$ is consistent.
\item{(d)} $\langle \bar{a}_\eta:\,\eta\in {}^{\omega>}2\rangle$ is 
3-fbti.
\end{description}

(In particular, such a formula and
witnesses can be found for any theory satisfying
SOP${}_1$ and NSOP${}_2$.)
\end{Claim}

\begin{Proof of the Claim} We shall be using the following colouring theorem,
for which we could not find a specific reference and so we include
a proof of it.

\begin{Lemma}\label{naknadni}
Suppose $\cf(\kappa)=\kappa$ and we colour ${}^{\kappa>}2$ by
$\theta<\kappa$ colours.
\underline{Then} there is an embedding $h:\,{}^{\omega>}2\into {}^{\kappa>}2$
such that $h(\eta)\hat{}\langle l\rangle\initialeq h(\eta\hat{}\langle
l\rangle)$ and Rang$(h)$ is monochromatic.
\end{Lemma}

\begin{Proof of the Lemma} Let $c$ be a colouring as in the assumptions and
let $\{a_i:\,i<\theta\}$ list $\Rang(c)$. We claim that there is $\nu^\ast
\in {}^{\kappa>}2$ and $j<\theta$ such that for every $\nu\in {}^{\kappa>}
2$ satisfying $\nu^\ast\initialeq\nu$ there is $\rho\in {}^{\kappa>}2$
with $\nu\initialeq\rho$ and $c(\rho)=j$. For otherwise, we can choose
by induction on $i\le\theta$ a member $\eta_i\in {}^{\kappa>} 2$ with
$i<j\implies \eta_i\initialeq\eta_j$ such that for no $\rho\in
{}^\kappa 2$ do we have $\eta_{i+1}\initialeq\rho$ and $c(\rho)=i$, using
$\theta<\cf(\kappa)$. As $\theta<\kappa$, we obtain a contradiction.

Having found such $\nu^\ast, j$ we define $h(\eta)$ for $\eta\in {}^n 2$
by induction on $n<\omega$. For $n=0$ we choose $h(\langle\rangle)$
to satisfy $\nu^\ast\initialeq h(\langle\rangle)$ and $c(h(\langle\rangle)
=j$, which is possible by the choice of $\nu^\ast$ and $j$. For $n+1$,
for any $\eta\in {}^{n+1} 2$ we choose for $l=0,1$ a member $h(\eta
\frown\langle l\rangle)$ of ${}^{\kappa>}2$ which is above
$h(\eta)\frown \langle l\rangle$ and on which $c$ is $j$, which
again is possible by the choice of $\nu^\ast$ and $j$.
$\eop_{\ref{naknadni}}$

\end{Proof of the Lemma}

Let $\varphi(\bar{x},\bar{y})$ be a SOP${}_1$
formula which is not SOP${}_2$, and moreover assume that for no
$n$ does the formula $\varphi_n$ defined as above satisfy SOP${}_2$. By Theorem
\ref{3indiswitness}, we can find witnesses $\langle \bar{a}_\eta:\,\eta\in {}^{\omega>}2\rangle$
which are 3-fbti. By the compactness theorem, we can
assume that we have a 1-fbti
sequence $\langle \bar{a}_\eta:\,\eta\in {}^{\omega_1>}2\rangle$ with the
properties corresponding to (a) and (b) of Definition \ref{nsop21}(2), namely
\begin{description}
\item{(a)} for every $\eta\in {}^{\omega_1}2$, the set
$\{\varphi(\bar{x},\bar{a}_{\eta\rest\alpha}):\,\alpha<\omega_1\}$ is
consistent,
\item{(b)} if $\nu\frown\langle 0\rangle\initialeq\eta\in {}^{\omega_1>}2$,
\underline{then} $\{\varphi(\bar{x},\bar{a}_{\nu\frown\langle 1\rangle}), 
\varphi(\bar{x},\bar{a}_{\eta})\}$
is inconsistent.
\end{description}

We shall now attempt to choose $\nu_\eta$ and $w_\eta$ for $\eta\in
{}^{\omega_1>}2$, by induction on $\llg(\eta)=\alpha<\omega_1$ so that:
\begin{description}
\item{(i)} $\nu_{\eta}\in {}^{\omega_1>}2$,
\item{(ii)} $\beta<\alpha\implies \nu_{\eta\rest\beta}\initial\nu_\eta$,
\item{(iii)} $\beta<\alpha\implies \nu_\eta(\llg(\nu_{\eta\rest\beta}))=
\eta(\beta)$,
\item{(iv)} $w_\eta\subseteq {}^{\omega_1>}2$ is finite and
$\nu\in w_\eta\implies \llg(\nu)<\llg(\nu_\eta)$,
\item{(v)} if $\llg(\eta)$ is a limit ordinal $>0$, then
$w_\eta=\emptyset$,
\item{(vi)} if $\eta\in {}^\beta 2$ and $l<2$, \underline{then}
$w_{\eta\frown\langle l\rangle}\subseteq\{\rho\in
{}^{\omega_1>}2:\,\nu_\eta\frown\langle l\rangle\initialeq\rho\}$
and $\max\{\llg(\rho):\,\rho\in w_{\eta\frown\langle l\rangle}
\}<\llg(\nu_{\eta\frown\langle l\rangle})$,
\item{(vii)} for each $\eta$ there is $\rho^\ast=\rho^\ast_\eta$
such that 
\begin{description}
\item{$(\alpha)$} $\nu_\eta\initial\rho^\ast\in {}^{\omega_1}2$,
\item{($\beta)$} $\card{\{\alpha<\omega_1:\,\rho^\ast(\alpha)=1\}}=\aleph_1$,
\item{$(\gamma)$} letting
\[
p_\eta(\bar{x})\deq\{\varphi(\bar{x}, \bar{a}_\Upsilon):\,\Upsilon\in
w_{\eta\rest\gamma} \mbox{ for some }\gamma\le\llg(\eta)\},
\]
we have that for all large enough $\beta^\ast$, the set
\[
p_\eta(\bar{x})\cup\{\varphi(\bar{x},\bar{a}_{\rho^\ast\rest\beta
}):\,\beta>\beta^\ast\,\,\wedge\,\,\rho^\ast(\beta)= 1\}
\]
is consistent,
\end{description}
\item{(viii)} 
$p_\eta(\bar{x})\cup\{\varphi(\bar{x},\bar{a}_\rho):\,\rho\in
w_{\eta\frown\langle 0\rangle}\cup w_{\eta\frown\langle 1\rangle}\}$
is inconsistent.
\end{description}
Before proceeding, we make several remarks about this definition. Firstly,
requirements (vii) and (viii) taken together imply that for each $\eta\in
{}^{\omega_1>}2$ we have that $w_{\eta\frown\langle 0\rangle}\cup
w_{\eta\frown\langle 1\rangle}\neq\emptyset$. Secondly, the definition of
$w_{\eta\frown\langle l\rangle}$ for $l\in \{0,1\}$ implies that 
\[
\wedge_{l=0,1}\rho_l\in w_{\eta\frown\langle l\rangle}\implies \rho_0\perp
\rho_1. \]
Thirdly, in (vii), any $\rho^\ast$ which satisfies that
$\nu_\eta\initial\rho^\ast$ and $\card{\{\gamma:\,\rho^\ast(\gamma)=1\}}=\aleph_1$ can be chosen as
$\rho^\ast_\eta$, by indiscernibility.

Now let us assume that a choice as above is possible, and we have made it.
Hence for each $\eta\in {}^{\omega_1>}2$ there is a finite $q_\eta\subseteq
p_\eta$ such that 
\begin{equation*}
q_\eta(\bar{x})\cup\{\varphi(\bar{x},\bar{a}_\rho):\,\rho\in
w_{\eta\frown\langle 0\rangle}\cup w_{\eta\frown\langle
1\rangle}\}\tag{$\ast$}
\end{equation*}
is inconsistent. Notice that there are $q$ and $\eta^\ast\in {}^{\omega_1}2$
such that \[
(\forall\eta_1)[\eta^\ast\initialeq\eta_1\in {}^{\omega_1>}2\implies
(\exists\eta_2\in {}^{\omega_1>}2)(\eta_1\initialeq\eta_2\,\,\wedge\,\,
q_{\eta_2}=q)].
\]
Namely, otherwise, we would have the following: each $p_\eta$ is countable,
hence for every $\eta$ there is
$g(\eta)$ with $\eta\initial g(\eta)\in {}^{\omega_1>}2$ and
\[
g(\eta)\initialeq \eta_1\implies q_{\eta_1}\nsubseteq p_\eta.
\]
Let $\eta_0\deq\langle\rangle$, and for $n<\omega$ let $\eta_{n+1}
=g(\eta_n)$. Let $\eta\deq\cup_{n<\omega}\eta_n$, hence $p_\eta=
\cup_{n<\omega} p_{\eta_n}$ (as $w_\eta=\emptyset$), and so $q_\eta
\subseteq p_{\eta_n}$ for some $n$, a contradiction.

Having found such $q,\eta^\ast$, by renaming
and using Lemma \ref{naknadni}, we can assume that 
$\eta^\ast\deq\langle\rangle$ and that for all $\eta\in {}^\omega 2$ we have
$q_\eta=p_{\langle\rangle}=
q$ (as $\eta\initialeq\nu\implies p_\eta\subseteq p_\nu$).
For $\eta\in {}^{\omega>}2$ let $\bar{\tau}_\eta$ list $w_\eta$. Without loss
of generality, by thinning and renaming, we have that for all
$\eta_1,\eta_2$,
\[
\langle\nu_{\eta_1}\rangle\frown\bar{\tau}_{\eta_1\frown\langle 0\rangle}\frown
\bar{\tau}_{\eta_1\frown\langle 1\rangle}
\iso_1
\langle\nu_{\eta_2}\rangle\frown\bar{\tau}_{\eta_2\frown\langle 0\rangle}\frown
\bar{\tau}_{\eta_2\frown\langle 1\rangle}.
\]
Similarly to the proof of Claim \ref{equivalent}, we can define a formula
$\psi(\bar{x},\bar{y})$ and $\{\bar{b}_\eta:\,\eta\in {}^{\omega>}2\}$
such that
\[
\psi(\bar{x}, \bar{b}_\eta)\equiv \bigwedge q\wedge
\bigwedge\{\varphi(\bar{x}, \bar{a}_\rho):\,\rho\in w_\eta\}.
\]

We claim that $\psi(\bar{x},\bar{y})$ and $\langle
\bar{b}_\eta:\,\eta\in {}^{\omega>}2\rangle$ exemplify SOP${}_2$ of $T$,
which is then a contradiction
(noting that $\psi$ is a formula of the form $\varphi_n$ for some $n$,
where $\varphi_n$ was defined in the statement of the Claim). We
check the two properties from Definition \ref{nsop21}(1). 

To see (a), let $\eta\in {}^\omega 2$ be given. We have that $p_\eta$ is
consistent, and $q\subseteq p_\eta$. For $n<\omega$, we have
\[
\psi(\bar{x}, \bar{b}_{\eta\rest n})\equiv \bigwedge q\wedge
\bigwedge\{\varphi(\bar{x}, \bar{a}_\rho):\,\rho\in w_{\eta\rest
n}\}. \]
As this is a conjunction of a set of formulae each of which is from  $p_\eta$,
we have that  $\{\psi(\bar{x}, \bar{b}_{\eta\rest n}):\,n<\omega\}$ is consistent. To check
(b), suppose $\eta\perp \nu\in {}^{\omega>}2$. Let $n$ be such that
$\eta\rest n=\nu\rest n$ but $\eta(n)\neq \nu(n)$. Hence
\[
\psi(\bar{x},\bar{b}_\eta)\equiv \bigwedge q \wedge
\bigwedge\{\varphi(\bar{x},\bar{a}_\rho):\,\rho\in w_{\eta\rest n\frown
\eta(n)}\}
\]
and
\[
\psi(\bar{x},\bar{b}_\nu)\equiv \bigwedge q \wedge
\bigwedge\{\varphi(\bar{x},\bar{a}_\rho):\,\rho\in w_{\eta\rest n\frown
\nu(n)}\},
\]
so taken together, the two are contradictory by the choice of $q$.

We conclude that the choice of $\nu_\eta$ and $w_\eta$ cannot be carried
throughout $\eta\in {}^{\omega_1>}2$. So, there is $\alpha<\omega_1$ and
$\eta\in {}^\alpha 2$ such that $\nu_\eta, w_{\eta\frown\langle l\rangle},
\nu_{\eta\frown\langle l\rangle}$
for $l<2$ cannot be chosen, and $\alpha$ is the first ordinal for which there
is such $\eta$.
Let $\nu^0_\eta\in {}^{\omega_1>}2\triangleright
\cup_{\beta<\alpha}\nu_{\eta\rest\beta} \frown
\langle \eta (\alpha-1)\rangle$ if the latter part is defined,
otherwise let $\nu^0_\eta\triangleright
\cup_{\beta<\alpha}\nu_{\eta\rest\beta}$. This choice
of $\nu_\eta=\rho$ for any $\rho\unrhd\nu_\eta^0$
with $\rho\in {}^{\omega_1}2$ satisfies items (i)-(iii)
above. We conclude that $w_{\eta\frown\langle l\rangle}$
for $l<2$ using any $\rho\unrhd\nu^0_\eta$ 
with $\rho\in {}^{\omega_1>}2$ for $\nu_\eta$ could not have been
chosen, and examine why this is so. Note that $p_\eta$ is already defined. Let
\begin{equation*}
\Theta\deq \left\{
\begin{split}(\rho,\gamma,w):\,&\nu^0_\eta\initial\rho\in {}^{\omega_1}
2,\\ &\llg(\nu^0_\eta)\le\gamma<\omega_1,\\
&(\exists^{\aleph_1}\beta <\omega_1)(\rho(\beta)=1),\\
&w\subseteq\{\Upsilon\in {}^{\omega_1>}2:\,\rho\rest\gamma\initialeq\Upsilon\}
\mbox{ is finite and}\\
&\mbox{ for some }\beta_\rho<\omega_1\mbox{ the set }\\
&p_\eta\cup\{\varphi(\bar{x},\bar{a}_{\rho\rest\beta}):\,\rho(\beta)=1
\,\,\&\,\,\beta\in [\beta_\rho,\omega_1)\}\cup\\
&\cup\{\varphi(\bar{x},\bar{a}_{\Upsilon}):\,\Upsilon\in w\}\\
&\mbox{ is consistent}\\
\end{split}
\right\}.
\end{equation*}
We make several {\bf observations}:

{\noindent (0)} If $(\rho,\gamma,w)
\in \Theta$ and
$w\subseteq w'$ with $w'$ finite and $w'\setminus
w$ is contained in $\{\rho\rest\beta:\,\beta_\rho
\le\beta\,\,\wedge\,\,\rho(\beta)=1\}$, \underline{then}
$(\rho,\gamma,w') \in \Theta$.

[This is obvious.]

{\noindent (1)} If $(\rho_l,\gamma,w_l)\in \Theta$ and for some
$\sigma\in {}^{\omega_1>}2$ with $\nu^0_\eta\initialeq \sigma$ we have
$\sigma\frown\langle l\rangle\initial \rho_l\rest\gamma$ for $l<2$,
while $\rho_0$ and $\rho_1$ are eventually equal,
\underline{then} $(\rho_l,\llg(\sigma), w_0\cup w_1)\in \Theta$.

[Why? We have $w_l\subseteq\{\Upsilon\in
{}^{\omega_1>}2:\,\rho_l\rest\gamma\initialeq\Upsilon\}$ is finite, so clearly
$w_0\cup w_1\subseteq\{\Upsilon\in
{}^{\omega_1>}2:\,\sigma\initialeq\Upsilon\}$ is finite. By the assumption, we
have that for some $\beta_l<\omega_1$ for $l<2$
\[
p_\eta\cup\{\varphi(\bar{x},\bar{a}_{\rho_l\rest\beta}):\,
\beta>\beta_l\,\,\wedge\,\,\rho_l(\beta)=1 
\} \cup\{\varphi(\bar{x},\bar{a}_{\Upsilon}):\,\Upsilon\in w_l\}
\]
is consistent. Suppose that (1) is not true
with $l=0$ and let $\beta^\ast\ge
\max\{\beta_0,\beta_1\}$ be such that $\beta^\ast<\omega_1$ and for
$\beta>\beta^\ast$ the equality $\rho_0(\beta)=\rho_1(\beta)$ holds. Hence
we have that
\[
p_\eta\cup\{\varphi(\bar{x},\bar{a}_{\rho_0\rest\beta}):\,
\beta>\beta^\ast\,\,\wedge\,\,\rho_0(\beta)=1\}
\cup\{\varphi(\bar{x},\bar{a}_{\Upsilon}):\,\Upsilon\in w_0\cup w_1\}
\]
is inconsistent. By increasing $w_0$ if necessary, (0) implies that
\[
p_\eta\cup\{\varphi(\bar{x},\bar{a}_\Upsilon):\,\Upsilon\in w_0\cup w_1\}
\]
is inconsistent. Let
$\nu_\eta\deq\sigma$, for $l<2$ let $w_{\eta\frown\langle l\rangle}=w_l$, and let $\nu_{\eta\frown\langle
l\rangle}\deq\rho_l\rest\beta^\ast_l$ for a large enough $\beta^\ast_l$ so that
$\beta^\ast<\beta^\ast_l$ and $\max(\{\llg(\Upsilon):\,\Upsilon\in
w_{\eta\frown\langle
l\rangle}\})<\beta^\ast_l$. This choice shows that we could have chosen
$\nu_\eta, w_{\eta\frown\langle
l\rangle}$ as required, contradicting the choice of $\eta$.]

{\noindent (2)} If $\nu^0_\eta\initial\rho\in {}^{\omega_1}2$ for some $\rho$
such that there are $\aleph_1$ many $\beta<\omega_1$ with $\rho(\beta)=1$, and
$\llg(\nu^0_\eta)\le\gamma<\omega_1$, \underline{then} $(\rho,\gamma,\emptyset)
\in \Theta$.

[Why? By the choice of $p_\eta$ and the remark about the freedom in the choice
of $\rho^\ast$ that we made earlier.]

Now we use the choice of $\eta$ to define witnesses to $T$ being SOP${}_1$
which also satisfy the requirements of the Claim. For $\tau\in {}^{\omega>}2$,
let $\bar{b}_\tau\deq\bar{a}_{\nu^0_\eta\frown\tau}$. Let us check the
required properties. Properties (a),(b) and (d) follow from the choice of
$\{\bar{a}_\sigma:\,\sigma\in {}^{\omega_1>}2\}$. Let $X^\ast\subseteq 
{}^{\omega>}2$ be such that there are no $\sigma,\nu\in X^\ast$ with
$\sigma\frown \langle 0\rangle\initialeq\nu$, we need to show that
$\{\varphi(\bar{x},\bar{b}_\tau):\,\tau\in X^\ast\}$ is consistent. It
suffices to show the same holds when $X^\ast$ replaced by an arbitrary 
finite $X\subseteq X^\ast$. Fix such an $X$. Clearly, it suffices to show that
for some $\rho, \gamma$, letting $w=\{\nu^0_\eta\frown\tau:\,\tau\in X\}$, we
have $(\rho,\gamma,w)\in \Theta$. 

Let $\rho^\ast\in {}^{\omega_1}2$ be such that $\nu^0_\eta\initial\rho^\ast$
and $\rho^\ast(\beta)=1$ for $\aleph_1$ many $\beta$.
By induction on 
$n\deq\card{X}$ we show:

there is $\rho\in{}^{\omega_1}2$ such that for some
$\gamma\ge\max\{\llg(\sigma):\,\sigma\in w\}$, we have
$(\rho,\gamma,w)\in \Theta$ and $\beta>\gamma\implies
\rho(\beta)=\rho^\ast(\beta)$, while
$\rho(\gamma)=1$.

\underline{$n=0$}. Follows by observation (2) above.

\underline{$n=1$}. Let $X=\{\tau\}$ and $\gamma=\llg(\tau)+\llg(\nu^0_\tau)$.
Let $\rho\in {}^{\omega_1}2$ be such
that $\rho\rest\gamma=\nu^0_\eta\frown\tau$, $\rho(\gamma)=1$ and
$\beta>\gamma\implies \rho(\beta)=\rho^\ast(\beta)$. By observation (2) above,
we have that $(\rho,\gamma,\emptyset)\in 
\Theta$. Then, by observation (0), we have
$(\rho,\gamma,w)\in \Theta$.

\underline{$n=k+1\ge 2$}. 
{\em Case 1}. $w$ is linearly ordered by $\initial$.

Let $\tau\in w$ be of maximal length, so clearly $\sigma\in
w\setminus\{\tau\}\implies\sigma
\frown\langle 1\rangle\initialeq\tau$. Let $\rho\in {}^{\omega_1}2$ be such that
$\tau\frown\langle 1\rangle\initial\rho$ and $\beta>\llg(\tau)$, while
$\rho(\beta)
=\rho^\ast(\beta)$. Now continue as in the case $n=1$.

{\em Case 2}. Not Case 1. 

Let $\sigma\in {}^{\omega_1>}2$ be $\initial$-maximal
such that $(\forall\tau)(\tau\in w\implies\sigma\initialeq\tau)$. This is well
defined,
as $w\neq\emptyset$ is finite. Let $w_l\deq\{\tau\in w:\,\sigma\frown\langle l\rangle
\initialeq\tau\}$, so $w_0\cap w_1=\emptyset$ but neither of $w_0,w_1$ is
empty. Now we have that $\sigma\notin w$, as otherwise we could choose
$\tau\in w_0$ such that $\sigma\frown\langle 0\rangle\initialeq\tau$, obtaining
an easy contradiction with our assumptions on $X$. Hence $w=w_0\cup w_1$.
We can now use observation (1) and the inductive hypothesis.
$\eop_{\ref{difference}}$
\end{Proof of the Claim}

To complete this discussion of the syntactic properties (N)SOP${1,2}$
we shall quote a result from \cite{ShUs xx} in which the understanding of SOP${}'_1$
and the
witnesses for SOP${}_1$ developed here was used to show that NSOP${}_1$
theories admit a rank function.

\begin{Definition}\label{rank} Given (partial) types $p(\bar{x}), q(\bar{y})$
and a formula $\varphi(\bar{x},\bar{y})$.
By induction on $n<\omega$ we define when
\[
\mbox{rk}^1_{\varphi(\bar{x},\bar{y})}(p(\bar{x}), q(\bar{y}))\ge n.
\]
{\noindent \underline{$n=0$.}} This happens iff both $p(\bar{x})$ and $q(\bar{y})$
are consistent.

{\noindent \underline{$n+1$.}} The rank is $\ge n+1$ iff for some $\bar{c}$
realising $q(\bar{y})$ both
\[
\mbox{rk}^1_{\varphi(\bar{x},\bar{y})}(p(\bar{x})\cup\{\varphi(\bar{x},\bar{c})\},
 q(\bar{y}))\ge n
\]
and 
\[
\mbox{rk}^1_{\varphi(\bar{x},\bar{y})}(p(\bar{x}),
 q(\bar{y})\cup\{\neg(\exists\bar{x})(\varphi(\bar{x},\bar{y})\wedge
 \varphi(\bar{x},\bar{c}))\})\ge n.
\]
If the rank is $\ge n$ for all $n$ then we say it is inifinite, otherwise we say it is finite.
\end{Definition}

\begin{Theorem}\label{ranksucess}(Shelah-Usvyatsov \cite{ShUs xx}) A theory $T$ is
NSOP${}_1$ iff
\[
{\rm rk}_{\varphi(\bar{x},\bar{y})}^1(\bar{x}=\bar{x}, \bar{y}=\bar{y})<\infty
\]
for every formula $\varphi(\bar{x},\bar{y})$.
\end{Theorem}

\section{$\initial^\ast$-maximality revisited}\label{ekvivalencija} 
In this section we come back to our main thesis, which is that properties
SOP${}_2$ and the maximality in the $\initial^\ast$-order are closely
connected.

Our main proof will use two auxiliary notions. The first is the order
$\initial^{\ast\ast}_\lambda$, which is a version
of the $\initial^\ast_\lambda$-order.

\begin{Definition}\label{red}
(1) For (complete first order theories) $T_1, T_2$
and a regular cardinal $\lambda>\card{T_1},\card{T_2}$, let
$T_1\initial^{\ast\ast}_\lambda T_2$ mean:

There is
a $\lambda$-relevant $(T_1,T_2)$-superior $(T^\ast,\bar{\varphi}, \bar{\psi})$ 
(see Definition \ref{order}) such that
$T^\ast$ has Skolem functions and if $T^{\ast\ast}\supseteq T^\ast$ is complete
with $\card{T^{\ast\ast}}<\lambda$ \underline{then}

$(\oplus)$ there is a model $M$ of
$T^{\ast\ast}$ of size $\lambda$ and an $M^{[\bar{\psi}]}$-type $p$
omitted by $M$ such that for every
elementary extension $N$ of $M$ of size
$\lambda$ which omits $p$ and a type $q$ (in one variable) over $N^{[\bar{\varphi}]}$,
there is an elementary extension of $N$ of size
$\lambda$ which realises $q$ and omits $p$.

{\noindent (2)} Let $T_1\initial^{\ast\ast} T_2$ mean that
$T_1\initial^{\ast\ast}_\lambda
 T_2$  holds for all large enough regular $\lambda$. 
 
 {\noindent (3)} $T_1$ is said to be $\initial^{\ast\ast}_\lambda$-maximal iff there is no $T_2$
 such that $T_1\initial^{\ast\ast}_\lambda T_2$. Similarly for $\initial^{\ast\ast}$.
\end{Definition}

The connection between this notion and $\initial^\ast$ is given by the
following claim:

\begin{Claim}\label{novi} Suppose that $T_1, T_2$
are theories and $\lambda>\card{T_1},
\card{T_2}$ satisfies $2^\lambda=\lambda^+$.
\underline{Then}
\[
T_1\initial^{\ast}_{\lambda^+} T_2\implies 
\neg( T_2\initial^{\ast\ast}_{\lambda} T_1).
\]
\end{Claim}

\begin{Proof} This statement is just a reformulation of the beginning
of the proof of Theorem \ref{notmax}. In other words, let
$(T,\bar{\varphi_1}, \bar{\varphi_2})$ show that $T_1\initial^{\ast}_{\lambda^+} T_2$.
This means that $\card{T}<\lambda^+$ but since 
$\lambda^{<\lambda}=\lambda$ and $\lambda>\card{T_1},
\card{T_2}$ we may assume that $\card{T^\ast}<\lambda$. Namely since there is a consistent
theory $T\supseteq \bar{\varphi_1}\cup \bar{\varphi_2}$ in which $\bar{\varphi_l}$
interprets $T_l$, and each $T_l$ has size $<\lambda$, there is a consistent
theory $T'$ of size $<\lambda$ which does the same. Without loss of generality
$T'\subseteq T$. In particular $\card{\tau(T')}
<\lambda$ so by extending $T'$ to a complete subtheory of $T$ and renaming we may
assume $T'$ is complete. Any model $M$ of $T$ has a reduct $N$ that is a model of $T'$
and that satisfies $M^{[\bar{\varphi}]}=N^{[\bar{\varphi}]}$ and similarly for $\bar{\psi}$.
Hence $(T',\bar{\varphi},\bar{\psi})$ is a $\lambda$-relevant $(T_1,T_2)$-superior
that exemplifies $T_1\initial^{\ast}_{\lambda^+} T_2$, so by renaming we may assume
$\card{T}<\lambda$.

Suppose for contradiction that $T_2\initial^{\ast\ast}_{\lambda} T_1$ and let
$(T^\ast,\bar{\varphi},\bar{\psi})$ exemplify this. Without loss of generality,
$\bar{\varphi_1}=\bar{\psi}$ and $\bar{\varphi_1}=\bar{\varphi}$ and the common vocabulary
of $T$ and $T^\ast$ is $\tau(\bar{\varphi_1})\cup \tau(\bar{\varphi_2})$. Hence
$T^{\ast\ast}=T\cup T^\ast$ is consistent by Robinson Consistency Criterium.
Without loss of generality $T^{\ast\ast}$ is complete. Hence let
$M$ be a model of $T^{\ast\ast}$ of size $\lambda$ and $p$ be a $M^{[\bar{\psi}]}$
type omitted by $M$ exemplifying the definition of $\initial_\lambda^{\ast\ast}$.
Using the assumption $2^\lambda=\lambda^+$ we
can build by induction an elementary extension $N$ of $M$ with $\card{N}=\lambda^+$, with
$N$  omitting $p$ and being $\bar{\varphi}$-saturated.
This is a contradiction with the choice of $T$.
$\eop_{\ref{novi}}$
\end{Proof}

\begin{Corollary}\label{GCH} Suppose that for all large enough regular $\lambda$ we have
$2^\lambda=\lambda^+$. \underline{Then} any $\initial^\ast$-maximal theory is also
$\initial^{\ast\ast}$-maximal.
\end{Corollary}

\begin{Proof} Suppose otherwise and let $T$ exemplify this. Hence for every $\kappa$
there is regular $\lambda>\kappa$ such that $T$ is not $\initial^{\ast\ast}$-maximal
and $2^\lambda=\lambda^+$. Hence $T$ is not $\initial^\ast_{\lambda^+}$-maximal
by Claim \ref{novi}, a contradiction.
$\eop_{\ref{GCH}}$
\end{Proof}

The next notion we need is a syntactic property.

\begin{Definition}\label{SOP"_2} Let $T$ be a theory.

{\noindent (1)} For a formula $\sigma (x,\bar{y})$
we say that $\sigma (x,\bar{y})$ has SOP${}''_{2}$ iff for some
[by compactness equivalently all] regular $\lambda>\card{T}$ there is
a sequence
\[
\left
\langle 
\bar{e}_{\bar{\eta}}:\,\bar{\eta}=\langle\eta_0,\ldots\eta_{n^\ast-1}\rangle,
\eta_0\initial\eta_1\initial\ldots\initial\eta_{n^\ast-1}\in {}^{\lambda>}
\lambda\mbox{ and }\llg(\eta_i)\mbox{ a successor}
\right\rangle
\]
such that 
\begin{description}
\item{($\alpha$)} for each $\eta\in {}^\lambda \lambda$, the set
\begin{equation*}
\left\{
\begin{split}
\sigma(x,\bar{e}_{\bar{\eta}}):\,&\bar{\eta}=
\langle\eta\rest(\alpha_0+1), 
\eta\rest(\alpha_1+1),\ldots \eta
\rest(\alpha_{n^\ast-1}+1)\rangle\\
&\mbox{ and }
\alpha_0<\alpha_1<\ldots\alpha_{n^\ast-1}<\lambda\\
\end{split}
\right\}
\end{equation*}
is consistent

\item{($\beta)$} for every large enough $m$,
\underline{if} $g:\,{}^{n^\ast\ge
}m\into{}^{\lambda>}\lambda$ satisfies \[
\rho\initial\nu\implies g(\rho)\initial g(\nu)
\]
and 
\[
\rho\in {}^{n\ge}m\implies\llg(g(\rho))\mbox{ is a successor},
\]
while
for $l<n^\ast-1$
\[
(g(\rho))\frown\langle l\rangle\vartrianglelefteq g(\rho\frown\langle
l\rangle), \]
\underline{then}
\[
\{\sigma(x,\bar{e}_{\langle g(\rho\rest 1), g(\rho\rest 2),
\ldots g(\rho)\rangle}):\,\rho\in {}^{n^\ast}m\}
\]
is inconsistent. Here $n^\ast=\llg(\bar{y})$ in $\sigma(x,\bar{y})$.
\end{description}

{\noindent (2)} $T$ is said to have SOP${}''_{2}$ iff some $\sigma (x,\bar{y})$
exemplifies it.
\end{Definition}

Our theorem \ref{glavnith} is phrased in terms of SOP${}''_{2}$. Answering
a question from an earlier version of this paper
Shelah and Usvyatsov proved in \cite{ShUs xx} the following Theorem \ref{SOP2andprime},
which then can be used together with theorem \ref{glavnith} to prove
Corollary \ref{zakljucak} which states that $\initial^\ast$-maximality implies SOP${}_2$.

\begin{Theorem}\label{SOP2andprime} (Shelah-Usvyatsov \cite{ShUs xx}) 
For any theory $T$, $T$ has SOP${}_2$ iff
it has SOP${}''_2$.
\end{Theorem}

\begin{Main Theorem}\label{glavnith} 
For any theory $T$ and regular cardinal $\lambda>\card{T}$, \underline{if}
$T$ is $\initial^{\ast\ast}_\lambda$-maximal \underline{then} $T$ has SOP${}''_2$.
\end{Main Theorem}

\begin{Proof} 
Let
$T$ be a given theory and let $\lambda=\cf(\lambda)>\card{T}$. We shall assume that
$T$ is $\initial^{\ast\ast}_\lambda$-maximal and prove that $T$ has SOP${}''_2$.
To make the reading of the proof easier we shall break it into stages.

{\bf Stage} {\bf A}.
Let $T_{\rm tree}^n\deq\mbox{ Th}({}^{n\ge}2, <_{\rm tr})$ for $n<\omega$,
where $<=<_{\rm tr}$ stands for 
the relation of ``being an initial segment of", and let $T_{\rm tree}\deq\lim
\langle T_{\rm tree}^n:\,n<\omega\rangle$, that is to say the set of all
$\psi$ which are in $T_{\rm tree}^n$ for all large enough $n$.
In order to use our assumptions at a later point, let us
fix a
theory $T^\ast$ 
which is a $\lambda$-relevant $(T_{\rm tree},T)$-superior with Skolem functions
(such a $T^\ast$ is easily seen to exist),
and let $\bar{\varphi}$,
$\bar{\psi}$ be the
interpretations of $T_{\rm tree}$ and $T$ in $T^\ast$, respectively. We can
without loss of generality, by renaming if necessary, assume that
$\LL(T)\subseteq \LL(T^\ast)$, so the interpretation $\bar{\psi}$ is trivial.

As $\card{T},\card{T^\ast}<\lambda$, we can find $A\subseteq\lambda$ which
codes $T$ and $T^\ast$. 
Working in ${\bf L}[A]$, we shall define a model $M$ of $T^\ast$ of size
$\lambda$ as follows. Let
\begin{equation*}
\begin{split}
\Gamma\deq T^\ast& \cup\{\varphi_=(x_\eta,x_\eta):\,\eta\in
{}^{\lambda>}\lambda\}\\
&\cup\{x_\eta<_\varphi x_\nu:\,\eta\initial\nu\in {}^{\lambda>}\lambda\}\\
&\{\neg(x_\eta<_\varphi x_\nu):\,\neg (\eta\initial\nu)\mbox{ for }\eta,\nu\in
{}^{\lambda>}\lambda)\}.\\
\end{split}
\end{equation*}
By a compactness argument and the fact that $\bar{\varphi}$ interprets
$T_{\rm tree}$ in $T^\ast$, we see that $\Gamma$ is consistent.
Let $M$ be a model of $\Gamma$ of size $\lambda=\lambda^{<\lambda}$ (as we are
in ${\bf L}[A]$). For $\eta\in {}^{\lambda>}\lambda$ let $a_\eta$ be the
realisation of $x_\eta$ in $M$. For $\eta\in {}^\lambda\lambda$, let
\[
p_\eta(x)\deq\{a_{\eta\rest\alpha}<_\varphi x:\,\alpha<\lambda\}
\]
By the choice of $M$ and the compactness
argument it follows that each $p_\eta$ is a (consistent) type. Note that for
$\eta_0\neq\eta_1\in {}^\lambda\lambda$, types $p_{\eta_0}$ and $p_{\eta_1}$
are contradictory. Let \[
p'_\eta(x)=\{a<_\varphi x:\,\mbox{for some }\alpha<\lambda,
a<_\varphi a_{\eta\rest\alpha}\}.
\]
By the axioms of $T_{\rm tree}$, we have that $p_\eta$ and $p'_\eta$ are
equivalent.
Now we observe
that by the size of $M$ there is $\eta^\ast\in
{}^\lambda\lambda$ such that the type $p'_{\eta^\ast}$ is omitted in $M$, and
$p'_{\eta^\ast}$ is not definable in $M$, i.e. for no formula
$\vartheta(y,\bar{z})$
and $\bar{c}\subseteq M$ do we have: for $a\in M$, the following are
equivalent:
$[a<_\varphi x]\in p'_{\eta^\ast}$ and
$M\models\vartheta[a,\bar{c}]$.
Let $p\deq
p'_{\eta^\ast}$ for such a fixed $\eta^\ast$. For $\alpha<\lambda$, let
$a_\alpha\deq a_{\eta^\ast\rest\alpha}$. We now go back to $V$
and make an observation about $M$. 

\begin{Subclaim}\label{jedan}
$T_{\rm tree}$ satisfies the following property:

for any formula $\vartheta(x,\bar{y})$ we have that $T_{\rm tree}\proves
\sigma=\sigma(\vartheta)$, where
\begin{equation*}
\begin{split}
\sigma\equiv(\forall\bar{y})[[(\forall x_1,x_2))\vartheta(x_1,\bar{y})&\wedge
\vartheta(x_2,\bar{y})\implies x_1\le_{\rm tr} x_2\vee x_2\le_{\rm tr} x_1)]\\
&\implies (\exists z)(\forall x)(\vartheta(x,\bar{y})\implies x\le_{\rm
tr}z)].\\
\end{split}
\end{equation*}
\end{Subclaim}

\begin{Proof of the Subclaim} Let $\vartheta(x,\bar{y})$ be given. By the
definition of $T_{\rm tree}$ we only need to show that $T^n_{{\rm
tree}}\proves\sigma$ for all large enough $n$, which is obvious as for every
$n$ the tree ${}^{n\ge} 2$ has the top level.
$\eop_{\ref{jedan}}$
\end{Proof of the Subclaim}

Hence the interpretation $\bar{\varphi}$ of $T_{\rm tree}$ in $T^\ast$
satisfies the same statement claimed about $T_{\rm tree}$.
We conclude:

$\tensor$ if $M\elementary N$ and $p$ is not realised in $N$, \underline{then}
there is no $\vartheta(x,\bar{c})$ with $\bar{c}\subseteq N$ such that 
$\vartheta(a_{\eta^\ast\rest\alpha},\bar{c})$ for all $\alpha<\lambda$ holds
and every two elements of $N$ satisfying $\vartheta(x,\bar{c})$ are
$<_\varphi$-comparable.

{\bf Stage} {\bf B}.
We shall choose a filtration $\bar{M}=\langle M_i:\,i<\lambda\rangle$
of $M$, and an increasing sequence $\langle \alpha_i:\,
i<\lambda\rangle$, requiring:

\begin{description}

\item{(a)} $M_i\elementary M$
and $M_i$ are $\elementary$-increasing continuous of size $<\lambda$,
with $M$ being the $\bigcup_{i<\lambda}
M_i$, 

\item{(b)} 
$a_{\alpha_i}\in M_{i+1}\setminus M_i$.

\end{description}
We may note that
the branch induced by $\{a_{\alpha_i}:\,i<\lambda\}$ is
the same as the one induced by $\{a_\alpha:\,
\alpha<\lambda\}$. Hence $p$
is realised in any model in which
$p'(x)\deq\{a_{\alpha_i}<_{\varphi} x:\,i<\lambda\}$
is realised (or even the similarly defined type using any
unbounded subset of $\{\alpha_i:\,i<\lambda\}$). Hence, by renaming,
without loss of generality we have $\alpha_i=i$ for all $i<\lambda$.

{\bf Stage} {\bf C}.
At this point we shall use the $\initial^{\ast\ast}_\lambda$-maximality
of $T$, which implies that
it is not true that $T\initial^{\ast\ast}_\lambda T_{\rm tree}$.
In particular, our $T^\ast$, $M$ and $p$ do not exemplify this, hence
there is $N$ with $M\elementary N$ and $\norm{N}=\lambda$, such
that $N$ omits $p$, but for some $N^{[\bar{\psi}]}$-type $q$ over $N$,
whenever $N\elementary N^+$ and $N^+$ realises $q$, also
$N^+$ realises $p$. By $\tensor$, the branch induced by $\{a_{\eta^\ast}\rest
\alpha:\,\alpha<\lambda\}$ is not definable in $N$, so
without loss of generality $N=M$. We can also assume that
$q$ is a complete type over $M^{[\bar{\psi}]}$.
Let us now use the choice of $q$ to define
for each club $E$ of $\lambda$ a family of
formulae associated with it, and to show that each of these
families is inconsistent.
We use the abbreviation c.d. for ``the
complete diagram of".

For any club $E$ of $\lambda$ we define
 \[
\Gamma_E\deq \mbox{c. d.}(M)\cup
q(x)\cup \{\neg(a_{i}<_\varphi\tau(x,\bar{b})):\,i \in E,
\tau\mbox{ a term of }T^\ast, \bar{b}\subseteq M_i\}. \]
Clearly, for any club $E$, if $\Gamma_E$ is consistent
then there is a model $N$ in which $\Gamma_E$ is realised. Identifying any
$b\in M$ with its interpretation in $N$ and letting $a^\ast$ be the
interpretation of $x$ from $\Gamma_E$, we can assume that $N$ is
an elementary extension of $M$ in which $q$ is realised by $a^\ast$.
As $T^\ast$ has Skolem functions, we have $M\elementary N$.
Let $N_1$ be the submodel of $N$ with universe
\[
A^\ast\deq M\cup\bigcup_{i\in E}\{\tau(a^\ast, \bar{b}):\,\bar{b}\subseteq
M_i\mbox{ and }\tau \mbox{ a term of }T^\ast\}.
\]
Note that the size of $N_1$ is $\lambda$.
Clearly, $N_1$ is closed under the functions of $T^\ast$, so $M\subseteq
N_1\subseteq N$. As $T^\ast$ has Skolem functions, we get that $M\elementary
N_1\elementary N$. By the third part of the definition of $\Gamma_E$,
$p$ is omitted in
$N_1$. This is in contradiction with our assumptions, as $a^\ast\in
N_1$ realises $q(x)$.

Hence we can conclude
\begin{equation*}
\mbox{ for every club }E\mbox{ of }
\lambda,\mbox{ the set } \Gamma_E \mbox{ is inconsistent.}
\label{gammae}
\end{equation*}

{\bf Stage} {\bf D}.
Now we start our search for a formula that exemplifies that $T$
has SOP${}''_2$. In the following
definitions, we shall use the expression ``an almost branch" or the
abbreviation a.b. to stand for a
set linearly ordered by
$<_{\varphi}$ (but not necessarily closed under $<_{\varphi}$-initial
segments and not necessarily unbounded). 
Let
\begin{equation*}
\Theta_{T^\ast}^0\deq
\left\{
\begin{split}
\vartheta(x, y, \bar{z}):\,
&\mbox{ there is }l=l_\vartheta<\omega \mbox{ such that}\\
&\mbox{for every }M^\ast\models T^\ast, a\in M^\ast, \bar{c}\subseteq M^\ast,
\mbox{ the set }\\
&  \vartheta(a,y,\bar{c})^{M^\ast} \mbox{ is the union of }
\le l \mbox{ a.b. in }
{M^\ast}^{[\bar{\varphi}]}\\
\end{split}
\right\},
\end{equation*}
and let $\Theta_{T^\ast}$ be the set of all
$\vartheta(x,\bar{y},\bar{z})$ of the form
$\bigvee_{j<n}\vartheta_j(x,y_j,\bar{z}_j)$ for some $\vartheta_0,\ldots
\vartheta_{n-1}\in \Theta_{T^\ast}^0$ (where $\bar{y}=
\langle y_j:\,j<n\rangle$ and $\bar{z}={}_{j<n}^\frown\bar{z}_j)$.
The formulae in $\Theta_{T^\ast}$ will be called candidates.
For
every candidate
\[
\vartheta(x,\bar{y},\bar{z})\equiv\bigvee_{j<n}\vartheta_j(x,y_j,\bar{z}_j)
\]
and a $\bar{\psi}$-formula $\sigma(x,\bar{t})$,
we consider the following game $\Game_{n,\sigma,\vartheta}$ (whose
definition also depends on our fixed $p$, $q$ and $\bar{M}$),
played by two players $\exists$ and $\forall$. The game starts by
$\exists$ playing 
$\bar{b}^0$
from ${}^{{\llg(\bar{z}_0)}}M$, then
$\forall$ playing $\alpha_0<\lambda$. After that
$\exists$ chooses 
$\beta_0\in (\alpha_0,\lambda)$
and
$\bar{b}^1\in {}^{\llg(\bar{z}_1)}M$ such that $\bar{b}^0\in {}^{\llg(\bar{z}_0)}
M_{\beta_0}$, 
after which
$\forall$ chooses $\alpha_1 <\lambda$ etc.,
finishing by $\exists$ choosing $\bar{b}^{n-1}\in {}^{\llg(\bar{z}_{n-1})}M$
and $\forall$ choosing $\alpha_{n-1}$, while $\exists$ chooses
$\beta_{n-1}\in
(\alpha_{n-1},\lambda)$ such that
$\bar{b}^{n-1}\in {}^{\llg(\bar{z}_{n-1})}M_{\beta_{n-1}}$. Player $\exists$
wins this game iff for some $\bar{e}\in {}^{\llg(\bar{t})}M$ we have
\begin{equation*}\label{tensor_1}\tag{$\tensor_1$}
\sigma(x,\bar{e})\in q\mbox{ and }M\models(\forall
x)[\sigma(x,\bar{e})
\implies
\vartheta(x, \langle a_{\beta_0},
\ldots, a_{\beta_{n-1}}\rangle,{}_{k<n}^\frown\bar{b}^k)].
\end{equation*}
(Note:
the constants $a_{\beta_{k}}$ are from the set $\{a_i:\,
i<\lambda\}$ we chose above.) Observe that every sequence
$\langle \alpha_0,\ldots\alpha_{n-1}\rangle\in
{}^{n}\lambda$ is an admissible sequence
of moves for $\forall$.

We shall show that for some $n\ge 1$ 
and $\sigma, \vartheta$, player $\exists$
has a winning strategy in the game $\Game_{n,\sigma,\vartheta}$,
where $\vartheta=\bigvee_{j<n}\vartheta_j$ as above.
As these are determined games, it suffices to show that
for some $n\ge 1$ and $\sigma, \vartheta$, player $\forall$
does not have a winning strategy.
Suppose
that this is not the case, arguing in 
$(\HH(\chi), \in , <^\ast_\chi, \bar{M}, p, q)$, where $\chi$ is
large enough and $<^\ast_\chi$ is a fixed well ordering of
$\HH(\chi)$. Fix
for a moment $(n,\sigma,\vartheta)$.
Player $\forall$ has a winning strategy in 
$\Game_{n,\sigma,\vartheta}$, which,
replacing the ordinals $\alpha_l$ by
constants $a_{\alpha_l}$, can be represented
by a sequence of functions $G_{n,
\sigma,\vartheta}^l$ for $l< n$ (in $(\HH(\chi),\in,<^\ast_\chi, \bar{M},
p, q)$),
where for $l<n$,
if the play up to time $l$ has been
$\bar{b}_0,\alpha_0,\beta_0,\ldots,\alpha_{l-1}, \beta_{l-1},
\bar{b}^l$, then
$G_{n,
\sigma,\vartheta}^l$ applied to this play is $a_{\alpha_l}$ for the $\alpha_l$ in
the choice of player $\forall$.
We shall assume that these functions
are the $<^\ast$-first which can act in this manner.
Using this and elementarity,
we notice that for every $n,\sigma,\vartheta$ the values of
$G_{n,\sigma,\vartheta}^l$ take place in $M$, and that
\[
 E_0\deq\{\delta<\lambda:\,(\forall \sigma, \vartheta)
(\forall n)(\forall l< n)[M\cap
\mbox{ Skolem}_{(\HH(\chi),\in,\bar{M}, G_{n, \sigma,\vartheta}^l)}(M_\delta)=
M_\delta] \}
\] is a club of $\lambda$ (as
$\card{T^\ast}$, $\norm{M_i}<\lambda$ for all $i$ and $\bar{M}$ is increasing
continuous). Let $E\deq \acc(E_0)$.
Consider now the set
$\Gamma_E$. It is contradictory, so there is a finite
subset of it which is contradictory. Hence for some $n_0, n_1, n_2<\omega$ and
formulae $\varrho_l(\bar{z}_l)\,(l<n_0)$ from the c.d.$(M)$,
formulae $\sigma_k(x,\bar{e}_k)\,(k<n_1)\in q(x)$,
ordinals $\delta_0<\ldots <\delta_{n_2-1}\in E$, 
a sequence $\langle \bar{b}_{j,l}:\,j<n_2, l<l_j\rangle$ with
$\bar{b}_{j,l}\subseteq M_{\delta_j}$ and terms $\langle\tau_{j,l}:\,j<n_2, l<l_j\rangle$
of $T^\ast$, the following is inconsistent:
\[
\bigwedge_{l<n_0}\varrho_l(\bar{z}_l)\wedge
\bigwedge_{k<n_1}\sigma_k(x,\bar{e}_k)\wedge
\bigwedge_{j<n_2, l<l_j}
\neg\left(a_{\delta_j}<_\varphi\tau_{j,l}(x,\bar{b}_{j,l})\right). \]
As $\varrho_l$ come from the c.d.$(M)$
and $q(x)$ is a complete type over $M^{[\bar{\psi}]}$,
we may assume that
$n_0=1$ and $n_1=1$.
Note that we must have $n_2\ge 1$ and that there is no loss of generality in
assuming that $\bar{b}_{j,l}=\bar{b}_j$ for all $l<l_j$ for $j<n$. We
shall omit the subscript 0 from $\varrho,\sigma, \bar{e}$. Let $n=n_2$ and
let us define $\vartheta_j(x,y_j,\bar{z}_j)$ for $j<n$ by
\[
\vartheta_j(x,y_j,\bar{z}_j)\equiv \bigvee_{l<l_j}
y_j<_\varphi\tau_{j,l}(x,\bar{z}_j), \]
and let $\vartheta=\bigvee_{j<n} \vartheta_j$.
Note that for each $j$ we have that $\vartheta_{j}\in \Theta_{T^\ast}^0$, as
$<_\varphi$ is a tree order. Hence $\vartheta$ is a candidate,
$\sigma(x,\bar{e})\in q(x)$, and since $M
\models\varrho[\bar{d}]$ for some $\bar{d}$ we have
\begin{equation*}
M\models(\forall x)[\sigma(x,\bar{e})\implies
\bigvee_{j<n}\vartheta_j(x, a_{\delta_j},
\bar{b}_j)]. \tag{$\ast$}\label{inM}
\end{equation*}
Now we consider the following play of
$\Game_{n,\sigma,\vartheta}$. Let $\exists$ choose $\bar{b}_0$. Recall that
$\bar{b}_0\subseteq M_{\delta_0}$.
The strategy
$G^0_{n,\sigma,\vartheta}$ of
$\forall$ yields an ordinal $\alpha_0$. By the choice of $E_0$
we have $\alpha_0<\delta_0$ and $\bar{b}_0\in M_{\delta_0}$, so we can let
$\exists$ choose $\beta_0=\delta_0$. Let
$\exists$ choose $\bar{b}_1$ and then let $\forall$ choose $\alpha_1$
according to the strategy, etc. At the end of the
play, player $\forall$ should have won (as he/she used the supposed winning
strategy), but clearly (\ref{inM}) implies that $\exists$ won,
a contradiction.

{\bf Stage} {\bf E}. We conclude that
(for our $\lambda,\bar{M},p, q$), for some
$\sigma,\vartheta$ and $n\ge 1$  the player $\exists$
has a winning strategy in the game $\Game_{n,\sigma,\vartheta}$, call it
{\em St}. Let us fix $n=n^\ast,\sigma,\vartheta$, and $St$ 
and
use them
to get SOP${}''_2$.

For any $\bar{\alpha}=\langle\alpha_0,\ldots,\alpha_{n-1}\rangle\in
{}^n\lambda$, we can let
$\langle 
\bar{b}^{\bar{\alpha}\rest {k}}, \beta^{\bar{\alpha}\rest ({k+1})}:\,k< n
\rangle$ be the sequence of moves that $\exists$ plays
by following the winning strategy {\em St} in a play in which
$\forall$ plays $\bar{\alpha}$, as the dependence is as marked. Let $E$ be a
club of $\lambda$ such that if
$k\le n$ and $\alpha_0<\ldots <\alpha_{k-1}<\delta\in E$,
then
$\bar{b}^{\langle\alpha_0,\ldots,\alpha_{k-1}\rangle}\in
{}^{\llg(\bar{z}_j)}M_\delta$.
(Such a club can be found by a method similar to the one used in Stage D).
Renaming the $M_i$ and $a_i$'s, we can without loss of generality
assume that $E=\lambda$.
For $\bar{\alpha}\in {}^n\lambda$ let
$\bar{e}^{\bar{\alpha}}$ be such that:
\[
M\models\forall x [\sigma(x,\bar{e}^{\bar{\alpha}})
\implies\bigvee_{j<n}\vartheta_j(x, a_{\beta^{\bar{\alpha}\rest (j+1)}},
\bar{b}_j^{\bar{\alpha}\rest (j+1)})].
\]

Notice that $\sigma$ is a formula in the language of $T$. 
We shall show that $\sigma$, together
with a conveniently chosen sequence of $\bar{e}_{\bar{\eta}}$'s,
exemplifies SOP${}''_2$.
The proof now proceeds similarly to the proof of Main Claim \ref{main}.
Namely

\begin{Lemma}\label{opetembeddings} There are sequences
\[
\langle N_\eta:\,\eta\in {}^{\lambda>}\lambda\rangle, 
\langle h_\eta:\,\eta\in {}^{\lambda>}\lambda\rangle
\]
such that
\begin{description}
\item{(i)} $h_\eta$ is an elementary embedding of $M_{\llg(\eta)}$ into
${\mathfrak C}_{T^\ast}$
with range $N_\eta$,
\item{(ii)} $\nu\initialeq\eta\implies h_\nu\subseteq h_\eta$,
\item{(iii)} for $\alpha\neq\beta<\lambda$ and $\eta\in {}^{\lambda>}\lambda$
we have
\[
h_{\eta\frown\langle\alpha\rangle}(a_{\llg(\eta)})\perp_{\varphi}
h_{\eta\frown\langle\beta\rangle}(a_{\llg(\eta)}),
\]
\item{(iv)} $N_{\eta_0}\cap N_{\eta_1}=N_{\eta_0\cap\eta_1}$ for all
$\eta_0,\eta_1$.
\end{description}
\end{Lemma}

\begin{Proof of the Lemma} This Lemma has
the same proof as that of
Main Claim \ref{main} Stage B. In the notation of that proof, ignore
$b_{\delta_i}$. When defining $\Gamma$ use
\[
\Gamma=\cup_{\alpha<\lambda}\Gamma^\alpha_0 \cup
\cup_{\alpha<\lambda}\Gamma_3^\alpha\cup \Gamma_4\cup \Gamma_2^+, \]
where $\Gamma^+_2=\{x^\alpha_0\perp_\varphi x^\beta_0:\,\alpha\neq\beta<\lambda
\}$ and $\Gamma_0^\alpha,\Gamma_3^\alpha$ and $\Gamma_4$ are defined as in
the proof of Main Claim \ref{main},
allowing for the replacement of ${}^{\lambda>}2$ by ${}^{\lambda>}\lambda$ by
using
$\{\bar{x}^\alpha:\,\alpha<\lambda\}$ in place of $\{\bar{x}^0,\bar{x}^1\}$.
Assumptions on $\Gamma_0^\alpha,\Gamma^+_2$ and $\Gamma_3^\alpha$ are
analogous to the ones we made in that proof. Fact \ref{cinjenica} still holds, except that we drop the last set from the definition of $r(\bar{x})$. The
rest
of the proof is the same, recalling
that the branch induced by $\{a_i:\,i<\lambda\}$ is undefinable in $M$. 
$\eop_{\ref{opetembeddings}}$
\end{Proof of the Lemma}

{\bf Stage} {\bf F}.
For $\eta\in {}^\lambda\lambda$, let $h_\eta\deq\cup_{\alpha<\lambda}
h_{\eta\rest\alpha}$. Let $q_\eta\deq h_\eta(q)$, hence each $q_\eta$ is
a consistent type.
For $\bar{\eta}=\langle\eta_0,\ldots,\eta_{n-1}\rangle$ and
$\eta_0\initial\ldots\initial\eta_{n-1}$ with $\llg(\eta_i)=\alpha_i+1$,
let $\bar{e}_{\bar{\eta}}\deq h_{\eta_{n-1}}
(\bar{e}^{\langle\alpha_0,\ldots\alpha_{n-1}\rangle})$.

Suppose now that $\eta\in {}^{\lambda}\lambda$ is given, and consider the
set
\[
\{\sigma(x,\bar{e
}_{\bar{\eta}}):\,\bar{\eta}=\langle\eta\rest(\alpha_0+1),  \ldots \eta\rest(\alpha_{n-1}+1)\rangle\mbox{ for some }\alpha_0<\ldots
 \alpha_{n-1}<\lambda \}.
\]
This set is a subset of $q_\eta$, and is hence consistent. This proves
property $(\alpha)$ from the definition of SOP${}''_2$.
For $(\beta)$, let $m$ be large enough and
$g:{}^{n\ge }m\into {}^{\lambda>
}\lambda$ be as in the statement of $(\beta)$.
For $\rho\in {}^{n}m$
let
$\bar{e}_{g_\rho}\deq\bar{e}_{\langle g(\rho\rest 1),\ldots
g(\rho) \rangle}$ (note that this is always defined). 
We shall now
show that the set 
\[
\{\sigma (x, \bar{e}_{g_\rho}):\,\rho\in {}^{n} m\}
\]
is inconsistent. Suppose otherwise, so let $d\in {\mathfrak C}_{T^\ast}$
realise it. For each $\rho\in {}^{n} m$, let
$\eta_\rho\in {}^\lambda\lambda \supseteq g(\rho)$ and
let $\bar{\alpha}^\rho\deq\langle\alpha^\rho_0,\ldots,\alpha^\rho_{n-1}
\rangle$ satisfy $\llg(g(\rho\rest k))=\alpha^\rho_k+1$ for $k\le n$,
so for each $k<n$ we have $g(\rho\rest (k+1))=\eta_\rho\rest (\alpha^\rho_k
+1)$. Now we have that for each $\rho\in {}^n m$
\begin{description}
\item{(i)}
$\sigma(x,\bar{e}_{g_\rho})
\equiv \sigma(x, h_{\eta_\rho\rest(\alpha^\rho_{n-1}+1)}
(\bar{e}^{\bar{\alpha}^\rho}))\in
q_{\eta_\rho}\rest\sigma_{\eta_\rho}(x)$
\item{(ii)}
$N_{\eta_\rho}\models(\forall x)[\sigma(x,\bar{e}_{g_\rho})
\implies \vartheta(x, \langle h_{\eta_\rho}(a_{\beta^{\bar{\alpha}^{\rho\rest
1} }}), \ldots h_{\eta_\rho}(a_{\beta^{\bar{\alpha}^\rho}})\rangle,
{}^\frown_{j<n} h_{\eta_\rho}(\bar{b}^{\bar{\alpha}^\rho_j}))]$
(hence the same holds in ${\mathfrak C}_{T^\ast}$),
\item{(iii)} 
\end{description}
\begin{equation*}
\begin{split}
\vartheta(x, \langle h_{\eta_\rho}(a_{\beta^{\bar{\alpha}^{\rho\rest 1}
}}), \ldots h_{\eta_\rho}(a_{\beta^{\bar{\alpha}^\rho}})\rangle,
{}^\frown_{j<n} h_{\eta_\rho}(\bar{b}^{\bar{\alpha}^\rho_j}))
&\implies \\
&\bigvee_{j<n}\vartheta_j(x,h_{\eta_\rho}(a_{\beta^{\bar{\alpha}^{\rho\rest
(j+1) }}}), h_{\eta_\rho}(\bar{b}^{\bar{\alpha}^{\rho\rest(j+1)}}_j))\\
\end{split}
\end{equation*}
for our $\vartheta_0,\ldots \vartheta_{n-1}$.

For each $\rho\in {}^n m$ let $j(\rho)<n$ be the first such
that 
\[
\vartheta_{j}(d
,h_{\eta_\rho}(a_{\beta^{\bar{\alpha}^{\rho\rest(j+1)}}}),
h_{\eta_\rho}(\bar{b}^{\bar{\alpha}^{\rho\rest(j+1)}}_j))
\]
holds. Let $l^\ast=\max\{l_0^\vartheta,\ldots,l_{n-1}^\vartheta\}$.

As $m$ is large enough, there
are $\rho_0,\ldots,\rho_{l^\ast}\in {}^n
m$ such that $j(\rho_s)=j^\ast$ for all $s\in \{0,\ldots,l^\ast\}$, while
$\rho_s\rest j^\ast$ is fixed and
$\rho_s (j^\ast)\neq \rho_t (j^\ast)$ for $s\neq t\le l^\ast$. 
(We use that there is a full ${}^{l^\ast+1\ge} n$ subtree $t^\ast$ of
${}^{n\ge}m$ such that for all $\rho\in t^\ast\cap {}^n m$ we have
$j(\rho)=j^\ast$. Choose $\rho_s$ belonging to $t^\ast$
and splitting at the level $j^\ast$).
In particular, $\alpha^{\rho_s}_0=\alpha_0,\ldots, \alpha^{\rho_s}_{j^\ast-1}=
\alpha_{j^\ast-1}$ is fixed, and so is $h_{\eta_{\rho_s}}\rest
M_{\alpha^\ast_{j-1}
+1}$,
but  \begin{equation*}
g(\rho_s)\rest(\alpha_{j^\ast-1}+2) \mbox{ for }s\le l^\ast
\mbox{ are incomparable in }{}^\lambda\lambda.\tag{$\ast\ast$}
\end{equation*}
Let $\bar{\alpha}\deq \bar{\alpha}^{\rho_0}$.

For each $\rho \in {}^n m$ and $k<n$ we have that
$\bar{b}^{\bar{\alpha}^{\rho\rest(k+1)}}\in M_{\alpha^\rho_{k+1}}$
(by the choice of $E$), so in
particular $\bar{b}^{\bar{\alpha}^{\rho\rest j^\ast}} \in
M_{\alpha^\rho_{j^\ast-1}+1}$, and hence
$h_{\eta_{\rho_s}}(\bar{b}^{\bar{\alpha}^{\rho\rest j^\ast}})$ is a fixed
$\bar{b}^\ast$. By the choice of $d$ and
definitions of $j^\ast,l^\ast$ and 
$\Theta_{T^\ast}$, there are $s\neq t<l_{\vartheta_{j^\ast}}\le l^\ast$ such
that $h_{\eta_{\rho_s}}(a_{\beta^{\bar{\alpha}^{\rho_s\rest(j^\ast+1)}}})$
and $h_{\eta_{\rho_t}}(a_{\beta^{\bar{\alpha}^{\rho_t\rest(j^\ast+1)}}})$  are
on the same almost branch. Now note that for all $\rho$
we have 
\[
a_{\beta^{\bar{\alpha}^{\rho\rest(j^\ast+1)}}}
\in M_{{\beta^{\bar{\alpha}^{\rho\rest(j^\ast+1)}}}+1}\setminus
M_{\beta^{\bar{\alpha}^{\rho\rest(j^\ast+1)}}}
\]
and $\beta^{\bar{\alpha}^{\rho\rest(j^\ast+1)}}> \alpha^\rho_{j^\ast}$.
Hence $h_{\eta_{\rho_s}}(a_{\beta^{\bar{\alpha}^{\rho_s\rest(j^\ast+1)}}})$
and $h_{\eta_{\rho_t}}(a_{\beta^{\bar{\alpha}^{\rho_t\rest(j^\ast+1)}}})$
are incomparable, by property (iii) in Lemma
\ref{opetembeddings}, a contradiction. This shows $(\beta)$ from
the definition of SOP${}''_2$, so finishing the proof.
$\eop_{\ref{glavnith}}$
\end{Proof}

Putting this together with Corollary \ref{GCH} and Shelah-Usvyatsov theorem 
\ref{SOP2andprime} above we get the following corollary \ref{zakljucak}. 

\begin{Corollary}\label{zakljucak} (1) Suppose that $T$ is a theory that is
$\initial^\ast$-maximal in some universe of set theory in which $2^\lambda=\lambda^+$ holds
for all large enough regular $\lambda$. \underline{Then} $T$ has SOP${}_2$.

{\noindent (2)} Suppose that $T$ is a theory that is
$\initial^\ast_{\lambda^+}$-maximal in some universe of set theory in which
$\lambda$ is regular and $2^\lambda=\lambda^+$.
\underline{Then} $T$ has SOP${}_2$.
\end{Corollary}

\begin{Proof} (1) Let $W$ be a universe of set theory in which $2^\lambda=\lambda^+$ holds
for all large enough regular $\lambda$ and in which $T$ is
$\initial^\ast$-maximal. Hence by Corollary \ref{GCH} $T$ is $\initial^{\ast\ast}$-maximal
in $W$ and hence by Main Theorem \ref{glavnith} in $W$ it satisfies SOP${}''_2$.
By Shelah-Usvyatsov Theorem \ref{SOP2andprime} above
$T$ satisfies SOP${}_2$ in $W$. An application of
the Compactness Theorem shows that satisfying SOP${}_2$ is absolute, hence 
$T$ satisfies SOP${}_2$ in $V$.

{\noindent (2)} This follows similarly, but more directly, from Main Theorem \ref{glavnith}
and the Shelah-Usvyatsov Theorem \ref{SOP2andprime}.
$\eop_{\ref{zakljucak}}$
\end{Proof}

This section hence provides us with the proof of one side of our thesis that
SOP${}_2$ and $\initial^\ast$-maximality are closely connected. Recall that
Shelah proved in \cite{Sh 500} that SOP${}_3$ implies $\initial^\ast$-maximality.
So an important open question (provided that SOP${}_3$ are not actually equivalent,
which we still do not know) is

\begin{Question} Does SOP${}_2$ imply $\initial^\ast$-maximality?
\end{Question}

In a partial answer to this question posed in an earlier version of the paper
Shelan and Usvyatsov in Theorem 3.12 of \cite{ShUs xx}
provided a local positive answer to this question, where by ``local" we mean
that they proved that any theory with SOP${}_2$ is $\initial^\ast$ above $T_{\rm tree}$
when only types localised by a certain formula are considered (see Definition \ref{localised}).

\eject


\begin{thebibliography}{xxxxxxx}

\bibitem[ChKe]{ChKe} C.C Chang and H.~J. Keisler, {\em Model Theory, 3rd ed.}, Studies in
Logic and Foundations of Mathematics, vol. 73, North-Holland, 1990.

\bibitem[DjSh 614]{DjSh 614} M. D\v zamonja and S. Shelah, {\em On the
existence of universals and an application to Banach spaces and triangle free
graphs}, submitted.

\bibitem[DjSh 659]{DjSh 659} M. D\v zamonja and S. Shelah, {\em
Universal graphs at the successor of a singular cardinal}, Journal of Symbolic Logic, 68 (2)
(June 2003), pg. 366-387.

\bibitem[DjSh 710]{DjSh 710} M. D\v zamonja and S. Shelah,
{\em On properties of theories which preclude the existence of universal
models},
submitted.

\bibitem[GrSh 174]{GrSh 174} R. Grossberg and S. Shelah, {\em On universal
locally finite groups}, Israel Journal of Mathematics, 44 (1983), 289-302.

\bibitem[HaLa]{HaLa} J.D. Halpern, H. Lauchli, {\em A Partition Theorem},
Transactions of the American Mathematical Society, vol. 124, Issue 2 (Aug.
1966), pg. 360-367.

\bibitem[Ke]{Keisler} H.~J. Keisler, {\em Six classes of theories}, Journal of
Australian Mathematical Society 21 (1976), 257-275.

\bibitem[Ki]{Kim} B. Kim and A. Pillay,
{\em  Simple theories}, Annals of
Pure and Applied Logic, 88 (2-3) (1997), 149-164 .

\bibitem[KjSh 409]{KjSh 409} M. Kojman and S. Shelah, {\em Nonexistence of
universal
orders in many cardinals}, Journal of Symbolic Logic 57 (1992), 875-891.

\bibitem[KjSh 447]{KjSh 447} M. Kojman and S. Shelah,
{\em The universality spectrum of stable unsuperstable theories},
Annals of Pure and Applied Logic, 1992.

\bibitem[Sh c]{Sh c} S. Shelah, {\em Classification theory and the number of
nonisomorphic models}, vol. 92 of Studies in Logic and the Foundation of
Mathematics,
North-Holland, Amsterdam, 1990.


\bibitem[Sh 457]{Sh 457} S. Shelah, {\em The Universality Spectrum: Consistency
for more classes} in {\em Combinatorics, Paul Erd\"os is Eighty},
Vol. 1, 403-420, Bolyai Society Mathematical Studies, 1993,
Proceedings of the Meeting in honor of P. Erd\"os, Keszthely, Hungary 7. 1993,
an improved version available at http://www.math.rutgers.edu/$\tilde{}$
shelaharch

\bibitem[Sh 500]{Sh 500} S. Shelah, {\em Toward classifying unstable theories},
Annals of Pure and Applied Logic 80 (1996) 229-255.

\bibitem[ShUs 844]{ShUs xx}
S. Shelah and A. Usvyatsov, {\em More on SOP${}_1$ and SOP${}_2$}, Annals of Pure and Applied Logic
155 (2008) 16-31.

\bibitem[PiHa]{PiHa} D. Pincus and J.D. Halpern, {\em Partitions of Products}, Transactions of the American Mathematical Society 267 no. 2 (1981) 549-568.

\end{thebibliography}
\end{document}